\font\headfont=cmr8
\headline={\hfil\headfont\ifodd\pageno
HOMOTOPY HYPERBOLIC 3-MANIFOLDS ARE HYPERBOLIC
\else
DAVID GABAI, G. ROBERT MEYERHOFF, AND NATHANIEL THURSTON
\fi\hfil
}

\overfullrule=0pt
\input psfig
\def\figdir{Pictures}
\def\sectionskip{\penalty-1000\vskip 18pt plus 6pt minus 6pt}
\sectionskip
\magnification=1200
\baselineskip = 0.2 true in
\centerline {\bf HOMOTOPY HYPERBOLIC 3-MANIFOLDS ARE HYPERBOLIC}

\centerline {David Gabai, G. Robert Meyerhoff, and Nathaniel Thurston}
\bigskip
\noindent {\bf Section 0: Introduction}
\medskip
This paper introduces a rigorous computer-assisted procedure for analyzing hyperbolic
3-manifolds.  This technique is used to complete the proof of several 
long-standing
rigidity conjectures in 3-manifold theory as well as to provide a new lower
bound for
the volume of a closed orientable hyperbolic 3-manifold. 
\medskip

\noindent {\bf
Theorem 0.1:} Let $N$ be a closed hyperbolic 3-manifold. Then

1) If $f\colon M \to N$ is a homotopy equivalence where $M$ is a closed
irreducible 3-manifold, then $f$ is homotopic to a homeomorphism.

2) If $f,g\colon M\to N$ are homotopic homeomorphisms, then $f$ is
isotopic to $g$.

3)  The space of hyperbolic metrics on $N$ is path connected.
\medskip

\noindent {\bf Remarks:}  Under the hypothesis that $M$ is hyperbolic,  conclusion i)
follows from
Mostow's rigidity theorem [Mo].  Under the hypothesis that N is Haken, 
conclusions i)-ii) follow from Waldhausen [W].  If $N$ is both Haken and
hyperbolic, then iii) follows by combining [Mo] and [W].  Since non-Haken
manifolds are necessarily orientable we will from now on assume that all
manifolds under discussion are orientable.
\medskip
Theorem 0.1 with the added hypothesis that a closed geodesic
$\delta \subset N$ has a
{\it non coalescable insulator family} was proven by  Gabai (see [G]).  Thus
Theorem 0.1 follows from [G] and the main technical result of this paper
which is,
\medskip
\noindent {\bf Theorem 0.2:}   If $\delta$ is a shortest geodesic in a closed orientable
hyperbolic
3-manifold, then $\delta$ has a 
non-coalescable insulator family.
\medskip

\noindent {\bf Remarks:}  If $\delta$ is the core of an embedded hyperbolic tube of radius
$(\ln 3)/2 =  0.549306...$ then $\delta$ has a non-coalescable insulator family by
Lemma 5.9 of [G].  
 In this paper we establish a
second condidition, sufficient to guarantee the existence of a non-coalescable insulator family
for $\delta.$ That is if {\it Corona} ($\delta) < 2\pi/3$. 
(Corona($\delta)<2\pi/3$ if tube radius($\delta)>(\ln 3)/2)$.
We use the expression ``$N$ satisfies the {\it insulator condition}"
when there is a geodesic $\delta$ which has a non-coalescable insulator family.

\medskip
We prove Theorem 0.2 by first showing that
all closed hyperbolic 3-manifolds, 
with seven families of exceptional cases, 
have embedded hyperbolic
tubes of radius $(\ln 3)/2$ about their shortest geodesics.  
Conjecturally, up to isometry, there are exactly six exceptional
manifolds (see Conjecture 4.2).
Second, we show that any shortest
geodesic $\delta$ in six of the seven families has
Corona($\delta)<2\pi/3$.
Finally, we show that the seventh family corresponds to Vol3, the third smallest
known hyperbolic 3-manifold, and that the insulator condition holds for Vol3.  Each of the three parts of the proof is carried out with the assistance of a rigorous computer program.

Here is a hint why Theorem 0.2 might be amenable to computer proof.
If a shortest geodesic $\delta$ in a hyperbolic 3-manifold $N$ does not have a
$(\ln 3)/2$
tube then there is a 2-generator subgroup $G$ of $\pi_1(N)=\Gamma$ which also
does not have
that property.  That is, after identifying $N={\bf H^3}/\Gamma$ and letting $Z={\bf H^3}/G$,
then a
shortest geodesic in $Z$ does not have a $(\ln3)/2$ tube.  $G$ is a group of 2
generators, generated by $f$ and $w$, where $f\in \Gamma$ is a primitive
hyperbolic isometry
whose fixed axis $\delta_0 \subset {\bf H^3}$ projects to $\delta$, and $w$ is an
isometry of
${\bf H^3}$ which takes $\delta_0$ to a nearest translate.  Here $\delta_0$ is a lift
of $\delta$
to ${\bf H^3}$.

The space of relevant 2-generator groups in Isom$({\bf H^3})$
naturally
lives in ${\bf C^3}$.  We show that except for seven small regions in ${\bf C^3}$ the shortest
geodesic in any discrete, torsion free, parabolic
free 2-generator group must have a $(\ln 3)/2$
tube.  Further, if Corona$(\delta)\ge 2\pi/3$, then there is a 2-generator
subgroup with that property.  
That is, there is a 2-generator subgroup $G$ of
$\pi_1(N)$ such that if $N_1={\bf H^3}/G,$ then Corona$(\delta_1)={\rm Corona}(\delta)\ge
2\pi/3$ for some shortest geodesic $\delta_1$ in $N_1.$  We show that away from
a single small open set in ${\bf C^3},$ every discrete, torsion free, parabolic
free 2-generator group $G$ satisfies Corona$(\delta)<2\pi/3,$ where $\delta$ is a
shortest geodesic in ${\bf H^3}/G.$
 We finally show that the exceptional open set
contains a
unique manifold which is Vol3.  In fact we show that if a shortest
geodesic $\delta$ in
$N$ satisfies Corona$(\delta)>2\pi/3$, then $N$ = Vol3.  A variant of the above
arguments shows that Vol3 satisfies the insulator condition.

This paper is organized as follows.  In Chapter 1 we describe a space 
${\cal P}^\prime \subset {\bf C^3}$ which
naturally parametrizes all relevant 2-generator groups.  We explain how a theorem of
Meyerhoff as
well as elementary hyperbolic geometry considerations imply that we need
only consider a
compact portion ${\cal P}$ of ${\bf C^3}$.  We explain
in detail the plan for proving Theorem 0.2.  We will actually
be working in the parameter space ${\cal W}=\exp({\cal P})$.  The technical reasons for
working in ${\cal W}$ is
described at the end of this section.  In Chapter 2 we
describe and prove the necessary results about the Corona function.
 In Chapter 3 we prove that the exceptional open set in $C^3$
contains only Vol3. Also if $\delta$ is a shortest geodesic in a closed
orientable
hyperbolic 3-manifold $N$ and Corona$(\delta)\ge 2\pi/3$, then $N=$Vol3.  Nonetheless, we are able to
show that 
Vol3 satisfies the insulator condition.   
In Chapter 4, we prove some applications, one of which is discussed briefly below.

In Chapters 5 through 8 we address the computer-related aspects of the proof.   Here, the method for describing the decomposition of the parameter space ${\cal W}$ into sub-regions is given, and the conditions used to eliminate all but seven of the sub-regions are discussed.  At the end of this chapter, the first part of a detailed example is given.  Eliminating a sub-region requires that a certain function is shown to be bounded appropriately over the entire sub-region.  This is carried out by using a first-order Taylor series approximation of the function together with a remainder bound.  Our computer version of such Taylor series with remainder bounds is called an {\it AffApprox} and in Chapter 6, the relevant theory is developed.  At this point, the detailed example of Chapter 5 can be completed.

Finally, in Chapters 7 and 8, round-off error analysis appropriate to our set-up is introduced.  Specifically, in Chapter 8, round-off error is incorporated into the {\it AffApprox} formulas introduced in Chapter 6.  The proofs here require an analysis of round-off error for complex numbers, which is carried out in Chapter 7.

We used two rigorous computer programs in our proofs---{\it verify} and {\it fudging}.  These programs are provided in Appendices 1 and 2.  Actually, {\it fudging} is a variation of {\it verify} and as such we only provide the sections of {\it fudging} that are changed.  The proofs amount to having {\it verify} and {\it fudging} analyze several computer files.  These computer files are available from the Geometry Center. Details about how to get them and the programs can be found at
\bigskip
{\noindent\tt http://www.geom.umn.edu:/locate/HomotopyHyperbolic}
\bigskip

A consequence of this work is that either the shortest geodesic in a closed
orientable 3-manifold $N$ has a $1.059191579962\ldots/2$ tube or $N=\ $Vol3.  The volume of Vol3 is $1.01...$ and by
[GM2] if $N$ has a $\log (3)/2 = 
0.549306\ldots$ tube about a geodesic, then the volume of $N$ is greater than $0.16668\ldots.$ 
This leaves some exceptional cases that can be analyzed using data provided by {\it verify}.  In any case, 
we obtain
\medskip
\noindent {\bf Theorem 4.5:}  If $N$ is a closed hyperbolic 3-manifold, then
${\rm vol}(N)>0.16668\ldots.$
\bigskip
\noindent {\bf Remark:}  The best published lower bound for volume is 0.001 by [GM1], which
improved the lower bound of 0.0008 of [M2].
\medskip
\noindent {\bf Acknowledgements:}  We thank The Geometry Center and
especially Al Marden and David Epstein for the vital and multifaceted roles
they played
in  this work.   Jeff Weeks and SNAPPEA provided  valuable data and
ideas.  In
fact it was the data from an undistributed version of
SNAPPEA that encouraged us to pursue a computer-assisted proof of Theorem 0.2.  

Bob Riley specially tailored his program
POINCAR\'E
to directly address the needs of our project.  His work provided
many leads in our search for
killerwords.  Further, he provided the first proofs that the six
exceptional regions (other than the Vol3 region) correspond to
closed orientable 3-manifolds.  The authors are deeply grateful for his
help. 

The first-named author thanks
the NSF for partial support.  Some of the first author's preliminary ideas
were formulated while visiting
David Epstein at the University of Warwick Mathematics Institute.  The
second-named author thanks the NSF for
partial support; the USC and Caltech Mathematics Departments for supporting
him as a visitor while much of
this work was done; and Jeff Weeks, Alan  Meyerhoff, and especially Rob
Gross for computer assistance.  The
third-named author thanks the NSF for partial support; and the Geometry Center and the Berkeley Mathematics Department for their support.

\sectionskip\magnification=1200
\baselineskip = 0.2 true in
\def\qed{{\bf \ \ \ \ \vrule height6pt width5pt depth1pt}}
\def\Relength{{\rm Relength}}
\def\length{{\rm length}}
\def\Arccosh{{\rm Arccosh}}
\def\trace{{\rm trace}}
\def\distance{{\rm distance}}
   
\noindent {\bf Chapter 1: Killer Words and the Parameter Space}
\bigskip\noindent{\bf Definition 1.1:}  
We will work in the upper-half-space model for
hyperbolic 3-space.  All isometries will be
orientation preserving.  If $f$ is an isometry, then we define
$\Relength(f)=\inf\{\rho(x,f(x))\mid x\in {\bf H^3}\}.$
Thus $\Relength(f)=0$ if and only if $f$ is either a parabolic or elliptic
isometry.  If $\Relength(f)\neq 0,$ then $f$
is hyperbolic and fixes a unique geodesic $\sigma$ in ${\bf H^3}.$  In that case
$\sigma$ is oriented (the negative end
being the repelling fixed point on $S^2_\infty$)  and the isometry $f$ is the
composition of a  rotation of $t \pmod{2\pi}$
radians along $\sigma$ (the sign of the angle of rotation is determined by the right-hand rule) followed by a pure translation of ${\bf H^3}$ along $\sigma$ of
$l = \Relength(f)$.  We define
$\length(f)=l+it.$
  If $\sigma$ is an oriented geodesic in ${\bf H^3}$, then it makes sense
to talk about an $l+it$ translation of
${\bf H^3}$ along $\sigma$, even when $l \le 0$.  

If $f$ is elliptic, then $f$ is a rotation of $t$ radians about some oriented geodesic.  If $f$ is elliptic, we define $\length(f) = |t| i,$ the absolute value accounting for the arbitrariness of the orientation of the fixed geodesic.
If $f$ is parabolic or the identity, we define $\length(f) = 0 + i0.$
So, for all isometries we have that $\Relength = {\rm Re}(\length).$

\medskip \noindent {\bf Definition 1.2:} If
$G$ is a subgroup of ${\rm Isom}({\bf H^3}),$ then  we say that $f$ is an element of {\it smallest
length} in $G$ if $\Relength(f)\le \Relength(g)$
for all $g\in G, \ g\neq {\rm id}.$

\medskip\noindent{\bf Convention 1.3:}  Let $B$ denote the oriented geodesic $t(0,0,1)$, with negative
end $(0,0,0).$  Let $C$ denote the
oriented geodesic with negative endpoint $(-1,0,0)$ and positive endpoint $(1,0,0)$.

\medskip\noindent{\bf Lemma 1.4:}
  If the isometry $f$ is represented by the matrix $A\in {\rm PSL}(2,{\bf C}),$ then
$$\length(f)=2 \Arccosh(\trace(A)/2),$$
where the branch of $\Arccosh$ with positive real values is taken, unless the real part is zero in which case the non-negative imaginary part is taken.

\medskip\noindent{\bf Proof:}  Because $\trace$ is a conjugacy invariant, we can normalize our set-up via conjugation and assume that the axis of $f$ is $B.$  As such, $A$ is a diagonal matrix, with diagonal entries $p$ and $p^{-1}.$

The action of $A$ on the bounding complex plane is simply multiplication by $p^2.$  Extending this action to upper-half-space in the natural way rotates the $z-$ axis by angle $\arg(p^2)$ and sends $(0,0,1)$ to $(0,0,|p|^2).$  Thus, ${\rm Im}(\length(f)) = \arg(p^2) = {\rm Im}(\ln(p^2))$ and, using the hyperbolic metric, 
${\rm Re}(\length(f)) = \ln(|p|^2) = {\rm Re}(\ln(p^2)).$
That is, $\length(f) = \ln(p^2).$  Thus, we need only show that $2 \ln(p) = 2 \Arccosh(\trace(A)/2).$  But this follows because $\cosh(\ln(p)) = (p + p^{-1})/2 = \trace(A)/2.$ \qed

\medskip\noindent{\bf Definition 1.5:}
  If $\sigma,\ \tau$ are disjoint oriented geodesics  in 
${\bf H^3}$ which do not meet at infinity,
then define
$\distance(\sigma,\tau)=\length(w)$ where $w\in {\rm Isom}({\bf H^3})$ is the hyperbolic
element which translates ${\bf H^3}$
along the unique common perpendicular between $\sigma$ and $\tau$ and which takes the
oriented geodesic $\sigma$ to the
oriented geodesic $\tau$.  The oriented common perpendicular from $\sigma$ to $\tau$ is
called the {\it orthocurve} between $\sigma$
and $\tau$.  The {\it ortholine} between $\sigma$ and $\tau$ is the complete oriented
geodesic in ${\bf H^3}$ which contains the
orthocurve between $\sigma$ and $\tau$.

If $\sigma,\ \tau$ intersect at one point in ${\bf H^3}$ then slight changes must be made in the above definition.  The ortholine has no natural orientation, the orthocurve is a point, and $w$ is an elliptic isometry.

If $\sigma,\ \tau$ intersect at infinity, then there is no unique common perpendicular, hence no ortholine, and 
$\distance(\sigma,\tau) = 0 + i0,$ or $0 + i\pi$ depending on whether or not $\sigma$ and $\tau$ point in the same direction at their intersection point(s) at infinity.

\medskip\noindent{\bf Lemma 1.6:}
  $\distance(\sigma,\tau) = \distance(\tau,\sigma) \qed$

\medskip\noindent{\bf Lemma 1.7:}
  If the isometry $f$ is represented by the matrix  
$$A = \left(\matrix{ a & b \cr 
                                c  & d \cr} \right) \in {\rm PSL}(2,{\bf C}),$$
then $2\distance(f(B),B) = 2\Arccosh(ad + bc).$  Again,  the branch of $\Arccosh$ with positive real values is taken, unless the real part is zero in which case the non-negative imaginary part is taken.

\medskip\noindent{\bf Proof:}  In the case where $B$ and $f(B)$ do not intersect at infinity, we will compute the length of $h$, the square of the transformation taking $B$ to $f(B)$ along their ortholine.  
$h = (f \circ \rho \circ f^{-1}) \circ \rho$ where $\rho$ is 180-degree rotation about $B$ and hence $(f \circ \rho \circ f^{-1})$ is 180-degree rotation about $f(B).$  $\rho$ and $f$ are represented by the matrices $$\pm \left(\matrix{ i & 0 \cr 
                                0  & -i \cr} \right) {\rm \ \ and\ \ } 
                         \left(\matrix{ a & b \cr 
                                c  & d \cr} \right)  \in {\rm PSL}(2,{\bf C}). $$ Hence,   $h = (f \circ \rho \circ f^{-1}) \circ \rho$ can be computed to have representation $$\left(\matrix{ ad + bc & 2ab \cr 
                                2cd & ad + bc \cr} \right).$$ 

Using Lemma 1.3, we have that $$2 \distance(f(B),B) =  \length(h) = 2\Arccosh(\trace(h)/2) = 2\Arccosh(ad + bc).$$

If $f$ fixes the point $(0,0,0)$ at infinity, then $c = 0,\ ad = 1$ and the formula holds.  Similarly for the other cases in which $B$ and $f(B)$ intersect at infinity.

\noindent Note that this formula only defines $\distance(f(B),B)$  modulo $i \pi,$ but this is sufficient for our purposes. \qed

\medskip\noindent{\bf Definition 1.8:}  Let $\delta$ be a geodesic in the hyperbolic 3-manifold $N.$
Then {\it tuberadius}$(\delta)=\sup\{r\mid$
there exists in $N$ an embedded $D^2\times S^1$ of radius $r$ centered about the
geodesic $\delta\}.$

\medskip\noindent{\bf Lemma 1.9:}  Let $\delta$ be a geodesic in the hyperbolic 3-manifold $N$ and
$\{\delta_i\}_{i\ge 0}$ be the set of its distinct lifts to ${\bf H^3},$
then tuberadius$(\delta)=\min\{{\rm Redistance}(\delta_0,\delta_i)\mid i\neq 0\}. \qed$

\medskip\noindent{\bf Definition 1.10:}  We define an open subset ${\cal P}^{\prime}$ of ${\bf C^3}$ which naturally parametrizes
the collection of conjugacy classes
of 2-generator subgroups $G \subset {\rm Isom} ({\bf H^3})$ with specified generators $f, w$ where $f$ is hyperbolic and $w$ is not parabolic.  By
conjugating $G$ we can assume that $f$ is a
positive translation of ${\bf H^3}$ along the geodesic $B,$  and that the
orthocurve from $w^{-1}(B)$ to $B$ lies on $C$
on the negative side of $C\cap B.$  Associated to $\{G,f,w\},$ 
the group $G$ with specified generators $f,w,$ is the parameter
$(L,D,R) = (l+it, d+ib, r+ia)$ where
$l+it=\length(f)$ and  $d+ib=\distance(w(B),B).$ The complex number $r+ia$ is
defined as follows.  The isometry $w$ is
the composition of two isometries, first a $d+ib$  translation of ${\bf H^3}$ along
$C,$ which takes $w^{-1}(B)$ to $B,$ as
oriented geodesics, followed by an $r+ia$ translation along $B.$

We are primarily interested in the set 
${\cal T}^{\prime}\subset {\bf C^3}$
which parametrizes all conjugacy classes of triples $\{G,f,w\}$ where $G$ is a group
generated by a {\it shortest} element $f$ which (positively) translates $B$ and $w\in G$
takes $B$ to a
{\it nearest} translate $w(B)$ such that $-\Relength(f)/2<
{\rm Redistance}($
(ortholine from $w^{-1}(B)$ to $B$),
(ortholine from $B$ to $w(B)))
\le \Relength(f)/2.$ By conjugating $G$ we can assume both 
that $f$ is a positive
translation along the geodesic $B$ and that the orthocurve from $w^{-1}(B)$ to $B$ lies in $C$ on the negative side of
$B\cap C.$ See Figure 1.1.

\midinsert
\centerline{\psfig{file=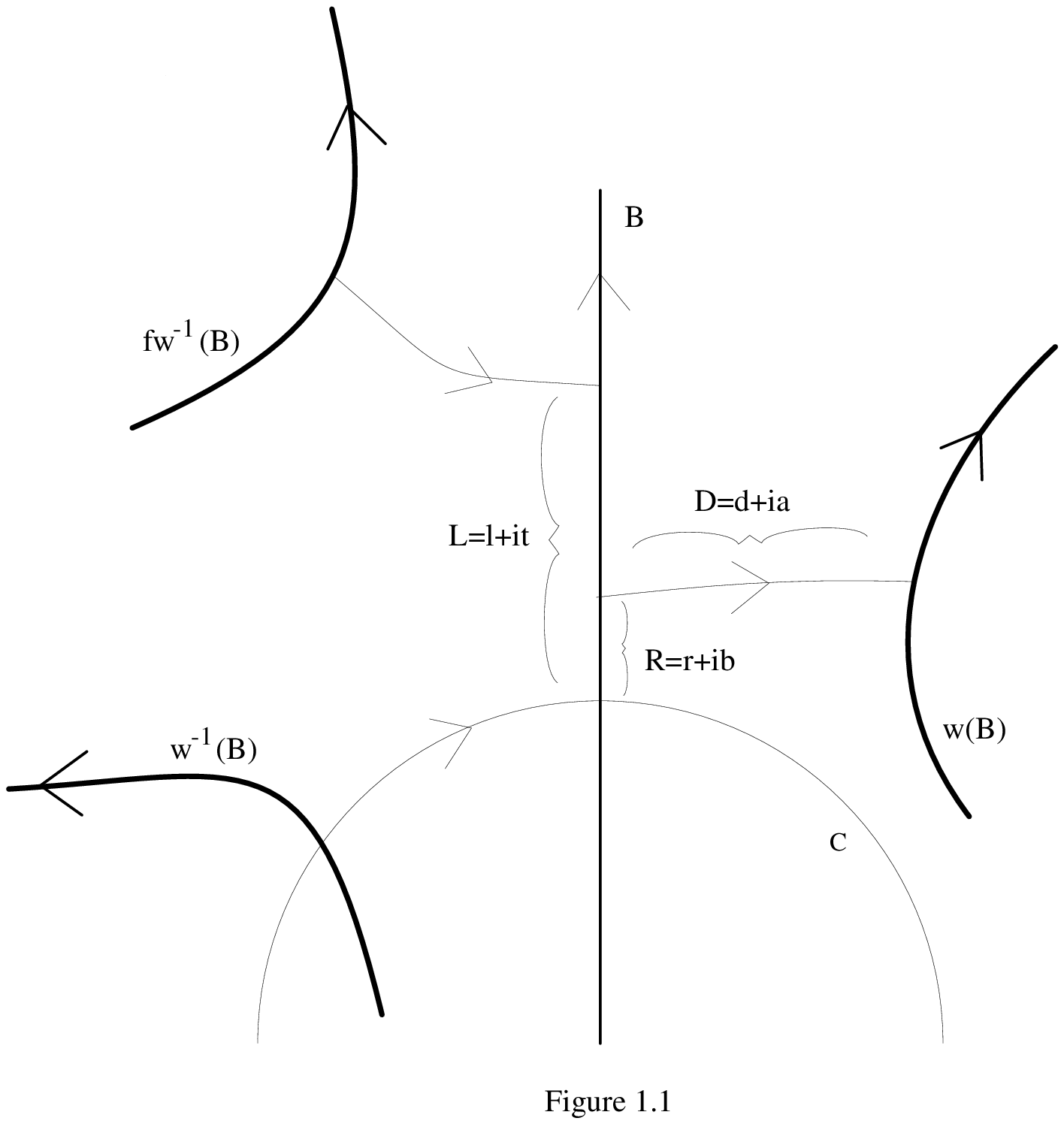,height=8cm}}
\endinsert
\bigskip
\noindent {\bf  Remark 1.11:}  Said another way, ${\cal T}^{\prime}$ corresponds to those parameters such
that $l$ is the real length of a
shortest element of $G,\  d$ is the real distance between $B$ and a nearest
translate, and $-l/2<r\le l/2.$  
In what follows, it is
essential to remember that an element $\alpha$ of ${\cal P^{\prime}}$ corresponds not only
to a group $G,$ but a group with two
special generators.  When $\alpha\in {\cal T}^{\prime},$ then two (i.e. $l$ and $d$) of
$\alpha's$ six real parameters correspond
to invariants of $\{G,f,w\}$. 

We are only interested in the subset of ${\cal T}^{\prime}$ corresponding to
parameters $\alpha$ with $d \le \ln(3).$  The
following two propositions imply this subset of ${\cal T}^{\prime}$ lives in a compact subset of ${\cal P}^{\prime}.$

\medskip\noindent {\bf Proposition 1.12:} All closed geodesics of length less than 0.0979 in
all hyperbolic 3-manifolds have embedded solid tube neighborhoods of
radius $\ln 3/2.$
\medskip\noindent {\bf Proof:}   In  [M1] it is proven that a geodesic of length $x + iy$ has an embedded solid-tube neighborhood of radius $r(x + iy)$ satisfying 
$$\sinh^2(r(x + iy)) = \max_{n \in {\bf Z_+}} {1 \over 2} ({\sqrt{1 - 2k(x,y,n)} \over k(x,y,n)} - 1) {\rm \ where\ } k(x,y,n) = \cosh(nx) - \cos(ny).$$
Of course, we restrict to $x+iy$ values which produce positive radii $r(x + iy)$ by means of this formula. It is easy to compute that for a given $x+iy$ we need to have $n$ for which $0 < k(x,y,n) < -1 + \sqrt 2$ to produce a positive radius tube by this method.

The function ${1 \over 2} ({\sqrt{1 - 2k} \over k} - 1)$ is decreasing on the interval $(0, -1 + \sqrt 2).$  It is therefore easy to solve for the range of $k$ values that produce radii $r$ greater than $\ln3/2.$  In fact, positive $k$ less than 0.3397 work.

Thus, to complete the proof of this proposition, we need to show that when a geodesic has real length $x$ less than 0.0979, that for all angles $y$ there exists a positive integer $n$ for which $k(x,y,n)$ is less than 0.3397.
Because $\cosh$ is an increasing function, we can restrict our analysis to $x = 0.0979.$  Thus,  we need only show that given any angle $y$, we can find a positive integer $n$ such that 
$\cosh(n 0.0979) - \cos(ny) < 0.3397.$  When $n > 8$ we can compute that $\cosh(n 0.0979) - \cos(ny) > 0.3397,$ and we therefore restrict to positive integers $n \leq 8.$  

We now consider angles $y$.  Because $\cos$ is an even function, we need only consider $y \in [0, \pi].$  
We will cover $[0,\pi]$ by 11 overlapping closed sub-intervals $\sigma_i$ each of which has an associated positive integer $n_i$ for which $\cosh(n_i 0.0979) - \cos(n_i y) < 0.3397$ is true for all $y \in \sigma_i.$

$$\sigma_0 = [0.000, 0.843] \ \ n_0 = 1$$
$$\sigma_1 = [0.835, 0.960] \ \ n_1 = 7$$
$$\sigma_2 = [0.951, 1.143] \ \ n_2 = 6$$
$$\sigma_3 = [1.123, 1.391] \ \ n_3 = 5$$
$$\sigma_4 = [1.386, 1.755] \ \ n_4 = 4$$
$$\sigma_5 = [1.733, 1.858] \ \ n_5 = 7$$
$$\sigma_6 = [1.832, 2.357] \ \ n_6 = 3$$
$$\sigma_7 = [2.334, 2.3792] \ \ n_7= 8$$
$$\sigma_8 = [2.3789, 2.647] \ \ n_8 = 5$$
$$\sigma_9 = [2.630, 2.755] \ \ n_9 = 7$$
$$\sigma_{10} = [2.730, \pi] \ \ n_{10} = 2$$
$\qed$

\medskip
\noindent {\bf Proposition 1.13:} If the {\it shortest} geodesic in a closed
hyperbolic 3-manifold has length greater than  or equal to 1.29, then it has
an embedded solid tube neighborhood of radius $\ln3/2.$
\medskip
\noindent {\bf Proof:} Consider the following folklore result:
The {\it shortest} geodesic in a closed hyperbolic 3-manifold has an
embedded solid tube neighborhood of radius $l/4$ where $l$ is the real
length of the shortest geodesic.  The proof is simple:
Expand a solid tube around the shortest geodesic, if it hits itself before a
radius of $l/4$ then we will construct a loop of length less than $l$, a
contradiction to ``shortestness."  Drop the two obvious perpendiculars from
the hitting point down to the core geodesic. Consider the following loop---
down one perpendicular, follow the shorter direction on the core geodesic,
up the other perpendicular.  This non-trivial loop has length less than $l/4
+ l/2 + l/4 < l.$

We now improve on this loop.  Replace the first half of the journey by the
hypotenuse of the right triangle formed by the first perpendicular and the
first half of the shorter arc along the core geodesic.  Replace the second
half of the journey by the hypotenuse of the right triangle formed by the
second perpendicular and the second half of the shorter arc along the core
geodesic.  Using the hyperbolic Pythagorean Theorem (see [F]) $\cosh c = (\cosh a)
( \cosh b)$ with $a = \ln 3/2$ and $b = l/4$ we get that the length of the
constructed loop is $2 \cosh^{-1} (\cosh(\ln 3/2) (\cosh (l/4)))$ and this is
less than $l$ when $l > 1.29$, by a calculation and the fact that $2 \cosh^{-
1} (\cosh(\ln 3/2) (\cosh (l/4)) - l$ is a decreasing function of $l$.

In fact, we can solve explicitly for the value of $l$ at which
$$2 \cosh^{-1} (\cosh(\ln 3/2) (\cosh (l/4)) - l = 0.$$  Noting that
$\cosh(\ln 3/2) = {2 \over \sqrt 3}$ we get ${2 \over \sqrt 3}(\cosh (l/4)))
= \cosh(l/2).$  Using a half-angle formula for $\cosh(l/2)$ we get ${2
\over \sqrt 3}(\cosh (l/4))) = 2\cosh^2(l/4) - 1.$  Setting $x =
\cosh^2(l/4)$ we get the quadratic ${2 \over \sqrt 3} x = 2 x^2 -1.$
Solving and substituting, we get $l = 1.289784... \qed$ 

\medskip \noindent {\bf Definition 1.14:} Let ${\cal P}\subset{\cal P}^{\prime}$ be those parameters $\alpha
= (l+it, d+ib, r+ia)$ such
that

a)      $-\pi \le t \le 0$

b)      $0\le r \le l/2$

c)       $0.0978\le l \le 1.29$

d)      $l/2 \le d\le \ln(3)$

\noindent Define ${\cal T}={\cal T}^{\prime}\cap {\cal P}.$

\medskip \noindent {\bf  Lemma 1.15:}  If $\alpha=(l+it, d+ib, r+ia)\in {\cal T}^{\prime}$ has $d\le \ln(3)$ and
corresponds to a 2-generator group
$\{G_\alpha,f_\alpha,w_\alpha\},$ then there exists a parameter $\beta\in {\cal T}$ 
with associated group $\{G_\beta,f_\beta, w_\beta\}$ such that $G_\beta$ is conjugate
to $G_\alpha.$

\medskip \noindent {\bf  Proof:} $d < l/2$ is eliminated from consideration by the first paragraph of the proof of Proposition 1.13, and the definition of ${\cal T}^\prime.$ If $d\le ln(3),$ then $0.0978\le l\le 1.29$ by Proposition 1.12.  If for
the triple $\{G, f, w\}$ we have $-l/2<r<0,$ then the 
triple $\{G,f,w^{-1}\}$ is conjugate to an element of ${\cal T}$ whose new
$r$-parameter is $-r.$  Thus we can assume that
b), c), and d) hold for the relevant $\{G,f,w\}.$  

This leaves property a).  Conjugating $G$ by a
reflection in the geodesic plane spanned
by $B$ and $C$ changes the $t$-parameter to $-t \pmod {2\pi}.$
The effect on $b$ and $a$ is irrelevant.  \qed

\medskip
By [G; Lemma 5.9] a closed orientable hyperbolic 3-manifold $N$
satisfies the insulator condition provided that
tuberadius$(\delta) > \ln(3)/2$ for some geodesic $\delta\subset N.$  Thus we
are led to ask

\medskip \noindent {\bf  Question 1.16:}   List all closed orientable hyperbolic 3-manifolds $N$
possessing a shortest geodesic $\delta$ such that tuberadius$(\delta)\le \ln(3)/2.$

\medskip \noindent {\bf  Remarks 1.17:}  i)  Using a J. Weeks-modified version of Snappea it was known
experimentally  that any shortest
geodesic in Vol3 has a .415... tube  (see [G]).  Conjecturally, up to isometry, there are a total of six manifolds in the list answering Question 1.16 (see Theorem 1.30 and Remark 1.31 iii, and Theorem 4.xx).

ii)  If a shortest geodesic $\delta$ in $N$ satisfies tuberadius$(\delta)\le
\ln(3)/2,$ then $N$ can be expressed as ${\bf H^3}/\Gamma$ where a lift of
$\delta$ is the geodesic $B,$ and $C$ is an ortholine between $B$ and a nearest $\Gamma-$translate.
Thus $N$ gives
rise to an element $\alpha\in {\cal T}.$  In fact $N$ may give rise to finitely many
different elements of ${\cal T}.$  Thus we need to
investigate.

\medskip \noindent {\bf  Question 1.18:}  Name all parameters $\alpha=(l+it, d+ib, r+ia) \in {\cal T}.$ 
\medskip
We now describe our method of (partially) answering Question 1.18. 
There is a technical point to mention: starting with Definition 1.22, we 
will work in the space ${\cal W} \supset \exp({\cal P})$, but for now we will describe the results in terms of the unexponentiated space ${\cal P}.$

We will partition ${\cal P}$ into about one billion regions
$\{{\cal P}_i\}$ and show that ${\cal T}$ is disjoint from all but seven small such regions.   Suppose that
${\cal P}_i$ is a region of this partition and $\alpha\in {\cal P}_i.$  Let $h$ be a word in the
letters $f,w$ and their inverses.  Associated to the parameter
$\alpha=(l_\alpha+it_\alpha, d_\alpha+ib_\alpha, r_\alpha
+ia_\alpha)$ there are the group elements
$f_\alpha, w_\alpha$ and hence $h_\alpha.$  Suppose that $h_\alpha\neq
f_\alpha^m.$  We ask
\medskip
a)  Is $\Relength (h_\alpha) < \Relength(f_\alpha) = l_\alpha?$

b)  Is Redistance$(h_\alpha(B),B) < $Redistance$(w_\alpha(B), B)= d_\alpha?$
\medskip
\noindent If either a) or b) is true, then $\alpha\notin {\cal T}.$  

Now let $\beta\in {\cal P}_i,$
with $f_\beta, w_\beta,$ and $h_\beta$ the associated hyperbolic isometries.  If
say a) is true for $\alpha$
then so is the statement $\Relength(h_\beta) < \Relength(f_\beta) = l_\beta$ for
$\beta$ sufficiently close to $\alpha.$
Thus we can show that ${\cal T}\cap {\cal P}_i =\emptyset$ if we can find an $\alpha$ for
which say statement a) is true, and then
use first-order Taylor approximation (with error/remainder term) to show that the corresponding
statement holds for all $\beta\in {\cal P}_i .$

\medskip \noindent {\bf  Definition 1.19:}  A word $h$ in $w,f,w^{-1},f^{-1}$ for which statement a) (resp.
b)) holds (non-trivially) for each $\beta\in {\cal P}_i$ is
called a {\it killer word} for ${\cal P}_i$ with respect to contradiction a) (resp. b)).

\medskip \noindent {\bf  Summary 1.20:}  With seven exceptions,  to each of the approximately one
billion regions partitioning ${\cal P},$ we will
associate a killerword and a contradiction.  

\medskip \noindent {\bf  Remark 1.21:}  Computers are well suited for partitioning a region such as ${\cal P}$
into many sub-regions $\{{\cal P}_i \},$ and finding a
killerword $h_i$ which eliminates $\alpha_i \in {\cal P}_i $ due to contradiction $C_i.$
Depending on the contradiction, we find
computable expressions  for approximations of the values of
$\Relength(h_\beta)$ or Redistance$(h_\beta(B),B)$ and
thus use the computer to eliminate all of ${\cal P}_i .$ 

There are a number of difficulties in executing this procedure.  First, a uniform mesh of the partition would
yield far too many sub-regions to be handled by computer.  In fact with 6
real parameters, refining a
given mesh by a factor of 10 would change the partition size by a factor of
$10^6$.  Our method for refining the
parameter space  and the way the computer keeps
track of the refinements are discussed briefly in Remark 1.27 and in more detail in
Chapter 5.

A second difficulty is finding the killerwords.  In practice, most of the
parameter space is eliminated by killerwords of length less than 7, but a number of spots need killerwords of length 10
and a few regions need killerwords of
length 35.  A brute force enumeration and testing of the various words
would take far too long.  Note that there are
more than 70,000 words of length 10 and 4$\times 3^{29}$ words of length 30.
Techniques for finding killer words can
be found in [Txxx].

Finally, there is the issue of rigor.
The main difficulty in making the plan work rigorously is that we need to bound
the difference between what the computer
thinks a value is and what the value actually is.  In particular we need to
control roundoff error, which becomes
quite significant when one does a large number of multiplications, e.g. in
the computation of $\length(h_\alpha)$ when $h$
is a 35 letter word.  Another issue is to make sure that $h_\beta$ is bounded
away from $f_\beta^m,$ and in particular
bounded away from $id,$ when computing say $\Relength(h_\beta).$   A large
portion of this paper is devoted to addressing
these  issues.  See Remark 1.30 for a more detailed discussion.

\medskip \noindent {\bf  Definition 1.22:}  Let 
$$ {\cal W} = \{ (x_0,x_1,x_2,x_3,x_4,x_5) : |x_i| \le 4 \times 2^{(5 - i) /6} {\rm \ for\ } i = 0,1,2,3,4,5 \}$$
$$\displaylines{\supset \exp({\cal P})=
\{(x_0,x_1,x_2,x_3,x_4,x_5)\mid \hfill\cr\hfill
x_0 + i x_3=\exp(e),  x_1 + i x_4 =\exp(f),  x_2 + i x_5 =\exp(g) 
{\rm where} (e,f,g)\in {\cal P}\}\cr}$$
 and let
$${\cal S}=\exp({\cal T}).$$
\noindent  Also, let $$L^\prime = \exp(L) = \exp(l+it),\ D^\prime = \exp(D) = \exp(d + ib),\ R^\prime = \exp(R) = \exp(r+ia).$$

\medskip \noindent {\bf  Remarks 1.23:} 
i) We work with ${\cal W}$ instead of $\exp({\cal P})$ because we want our initial region to be a (6-dimensional) box that is easily sub-divided.  This has the side-effect that certain sub-boxes ${\cal W}_i$ of ${\cal W}$ will be eliminated because they are outside of $\exp({\cal P})$ rather than by the analogues of conditions a) and b) above.
 The entire collection
of conditions is given in Section 5.

ii) All the ideas expressed in 1.18-1.21 will be carried out in the
parameter space ${\cal W}$ rather than the space
${\cal P}.$  That is mainly because of Lemmas 1.24 - 1.26 which demonstrate
that while working in ${\cal W}$ one need only
understand the basic arithmetic operations $+, -, \times, /, \sqrt{}$.

iii)  The reason for choosing the co-ordinates of ${\cal W}$  so that $L^\prime = x_0 + i x_3,\ D^\prime = x_1 + i x_4,\ R^\prime = x_2 + i x_5\ $ was to gain a mild computer advantage.  

\medskip \noindent {\bf  Lemma 1.24:}  If $(L^\prime, D^\prime, R^\prime)\in {\cal W}$ and  $f,w$ are the generators of the associated group $G$ then 

\noindent a) $${\rm Matrix}[f]= \left(\matrix{\sqrt{L^\prime}&0\cr 
0 & 1/\sqrt{L^\prime}\cr}\right)$$

\noindent b) $${\rm Matrix}[w]= \left(\matrix{\sqrt{R^\prime} * ch &  sh * \sqrt{R^\prime} \cr sh / \sqrt{R^\prime} & ch/\sqrt{R^\prime}\cr}\right)$$ 

where 
$ch = (\sqrt{D^\prime} + 1/\sqrt{D^\prime})/2\ \ {\rm and}\ \ 
sh = (\sqrt{D^\prime} - 1/\sqrt{D^\prime})/2$

\medskip \noindent {\bf  Proof:}  a)  By our set-up we have that the (oriented) axis of $f$ is $B$.  Following the proof of Lemma 1.3, $${\rm Matrix}[f]= \left(\matrix{p & 0 \cr 0 & 1/p \cr}\right)$$ where $p = \exp(\length(f)/2) = \sqrt{\exp(\length(f))} = \sqrt{\exp(L)} = \sqrt{L^\prime}$

b)  $w = \beta \circ \alpha$ where $\beta$ is translation of distance $R$ along $B$ and $\alpha$ is translation of distance $D$ along $C$.  Thus, $${\rm Matrix}[\beta]= \left(\matrix{\sqrt{R^\prime} & 0 \cr 0 & 1/\sqrt{R^\prime}\cr}\right)$$ and Matrix$[\alpha]$ can be computed to be $$\left(\matrix{\cosh(D/2) & \sinh(D/2) \cr \sinh(D/2) & \cosh(D/2)\cr}\right).$$ But $\cosh(D/2) = (\exp(D/2) + \exp(-D/2))/2 = 
(\sqrt{D^\prime} + 1/\sqrt{D^\prime})/2 = ch$ and similarly for $sh.$
Thus, $${\rm Matrix}[\alpha]= \left(\matrix{ch & sh \cr sh & ch\cr}\right)$$ and b) follows by matrix multiplication.
\qed
\medskip \noindent {\bf  Lemma 1.25:} If $h \in {\rm Isom}({\bf H^3})$ is represented by the matrix $$A = \left(\matrix{a &  b \cr c & d\cr}\right)\in {\rm PSL}(2,{\bf C}),$$ then 

\noindent a) $\exp(\Relength (h))= |\trace(A) \pm \sqrt{(\trace(A)/2)^2 - 1}|^2$

\noindent b) $\exp({\rm Redistance} (h(B),B))=|{\rm orthotrace}(A) \pm \sqrt{({\rm orthotrace}(A)/2)^2 - 1}|$ 

where 
orthotrace($A) = ad + bc.$

\noindent Here the $+,\ -$ produce reciprocal values for $\exp(\Relength (h)),$ and we take the one producing the larger value, unless the value is 1, in which case there is no choice.

\medskip \noindent {\bf Proof:}  Because $\cosh(x) = (\exp(x) + \exp(-x))/2$ it is easy to compute that $\cosh^{-1}(x) = \log(x \pm \sqrt{x^2 - 1}).$  Of course, $x - \sqrt{x^2 - 1}$ and $x + \sqrt{x^2 - 1}$ are inverses, which corresponds to the fact that $\cosh^{-1}(x)$ 
for $x \ne 1$ consists of values differing by a factor of $-1.$
\medskip
\noindent
a)  $\exp(\Relength (h)) = |\exp(\length(h))| = |\exp(2 \Arccosh(\trace (A)/2))| = |(\trace(A)/2) \pm \sqrt{(\trace(A)/2)^2 - 1}|^2$ where the second equality follows from Lemma 1.2.
\medskip
\noindent
b) $\exp({\rm Redistance} (h(B),B))= |\exp(\distance(h(B), B))| = |\exp(\Arccosh(ad + bc))| = |(ad + bc) \pm \sqrt{(ad + bc)^2 - 1}|$ where the second equality follows from Lemma 1.7.
\qed
\medskip \noindent {\bf Remarks 1.26:} i)  It follows from Lemma 1.25 that if $h$ is a word in $f,w$ and their inverses, then for any
parameter value $\alpha\in {\cal W},\ \ $
$\exp(\Relength(h_\alpha)),\  $ and 
$\exp({\rm Redistance}(h_\alpha(B), B))$ can be computed using only the
operations $+, -, \times, /, \sqrt{}. \qed$

 ii) During the course of the computer work needed to prove the main theorems, the parameter space
${\cal W}$ was decomposed into sub-boxes by computer via a recursive sub-division process:
Given a sub-box that is being analyzed, either it can be handled, or it
cannot.  If it cannot be handled, it is sub-divided in half by a hyper-plane
$\{x_i = c \}$ (where $i$ runs through the various co-ordinate dimensions
cyclically) and the two pieces are analyzed separately.  And so on. 

As such, a sub-box of ${\cal W}$ can be described by a sequence of 0's and 1's where 0 means ``take the lesser $x_i$ values" and 1 means ``take the greater $x_i$ values."  
Remarkably, for the decomposition of ${\cal W}$ into sub-boxes, all the sub-
box descriptions could be neatly encoded into one tree (although in practice we found it preferable to use several trees to describe the entire decomposition).  This is described in
Chapter 5.

iii) In the following proposition, seven exceptional boxes are described as sequences of 0's and 1's.  Four of the boxes---$X_0, X_4, X_5, X_6$--- are each the union of two abutting sub-boxes, $X_0 = X_{0a} \cup X_{0b}$ and so on.  It is a pleasant exercise to work through the fact that they abut.  

It is also a pleasant exercise to calculate by hand the co-ordinate ranges of the various sub-boxes.  For example, the range of the last co-ordinate (i.e., $x_5$) of the sub-box 
$$\displaylines{X_{6a} = 
111000000001000111\ 
111111110101001111\ 
011111010111111111\hfill\cr\hfill  
110001001011000111\ 0\cr}$$
 \noindent is found by taking the 6th entry, the 12th, entry, the 18th entry and so on.  These entries are 011111111111.  The first entry (0) means take the lesser $x_5$ values, and produces the interval $[-4,0].$  The second entry (1) means take the greater $x_5$ values, and produces the interval $[-2,0].$  The third entry (1) produces $[-1,0].$  Continuing, we see that $X_{6a}$ has $-2^{-9} \le x_5 = {\rm Im}(R') \le 0.$  The other co-ordinates can be computed in the same fashion, although they must at the end be multiplied by the facter $2^{(5 - i)/6}$ (see the definition of the initial box ${\cal W}$).  The range of co-ordinate values is given for each of the seven boxes.  Finally, two {\it quasi-relators} are given for each sub-box $X_0, X_1, \ldots, X_6.$

\medskip\noindent {\bf Definition 1.27:} A {\it quasi-relator} in a sub-box $X$ of ${\cal W}$ is a word in $f,w,f^{-1},w^{-1}$ that is close to the identity throughout $X$ and experimentally appears to be converging to the identity at some point in $X.$  In particular, a quasi-relator rigorously has Relength less than that of $f$ at all points in $X.$

\medskip \noindent {\bf  Proposition 1.28:}  $S\cap ({\cal W} - \bigcup_{n = 1,\dots, 7} X_n)=\emptyset$ where the $X_n$ are the exceptional sub-boxes

$$\displaylines{X_{0a} = 
001000110111110001\ 
101001010101011001\ 
011011010111101101\hfill\cr\hfill  
100001101101000111\ 
010001110101100101\ 
1101110111110100\cr}$$
$$\displaylines{X_{0b} = 
001001110110110000\ 
101000010100011000\ 
011010010110101100\hfill\cr\hfill  
100000101100000110\ 
010000110100100100\ 
1101100111100100\cr}$$
$$X_0 = \left(\matrix{
-0.840655162503\ldots  \le  Re(L')  \le  -0.840600786360\ldots\cr
-0.840642408899\ldots  \le  Re(D')  \le  -0.840593965263\ldots\cr
0.999979499517\ldots  \le  Re(R')  \le  1.000022657890\ldots\cr
-2.137267196028\ldots  \le  Im(L')  \le  -2.137228746289\ldots\cr
-2.137295441962\ldots  \le  Im(D')  \le  -2.137226932315\ldots\cr
-0.000061035156\ldots  \le  Im(R')  \le  0.000061035156\ldots\cr
}\right)$$

\noindent $X_0\ \ $ quasi-relators:

$r_1 = fwFwwFwfww$
 
$r_2 = FwfwfWfwfw$

$$\displaylines{X_1 = 
001000110001110110\ 
011101000110111110\ 
100010110000100011\hfill\cr\hfill  
101101001101001000\ 
110101011000000100\ 
000\cr}$$

$$X_1 = \left(\matrix{
-1.348528333122\ldots  \le  Re(L')  \le  -1.348310828552\ldots\cr
-0.543343817104\ldots  \le  Re(D')  \le  -0.543150042561\ldots\cr
0.903908961497\ldots  \le Re(R')  \le  0.904081594988\ldots\cr
-2.661029541660\ldots  \le  Im(L')  \le  -2.660721943747\ldots\cr
-2.858770529287\ldots  \le  Im(D')  \le  -2.858496490701\ldots\cr
-1.471679687500\ldots  \le  Im(R')  \le  -1.471435546875\ldots\cr
}\right)$$

\noindent $X_1\ \ $ quasi-relators:

$r_1 = FFwFWFWfWFWFwFFww$
 
$r_2 = FFwwFwfwfWfwfwFww$

$$\displaylines{X_2 = 
001000110101010010\ 
101010110001100101\ 
110111100001101010\hfill\cr\hfill  
111100100000010001\ 
111100\cr}$$

$$X_2 = \left(\matrix{
-1.787017545957\ldots  \le  Re(L')  \le  -1.785277509398\ldots\cr
-1.074286063490\ldots  \le  Re(D')  \le  -1.072735867150\ldots\cr
0.741633479486\ldots  \le Re(R')  \le  0.743014547418\ldots\cr
-2.272533378081\ldots  \le  Im(L')  \le  -2.271302986431\ldots\cr
-2.718462773249\ldots  \le  Im(D')  \le  -2.717366618905\ldots\cr
-1.529296875000\ldots  \le  Im(R')  \le  -1.528320312500\ldots\cr
}\right)$$

\noindent $X_2\ \ $ quasi-relators:

$r_1 = FwfwfWffWfwfwFww$

$r_2 = FFwFFwwFwfwfwFww$

$$\displaylines{X_3 = 
111000000001000110\ 
011011101101011000\ 
111101011110001100\hfill\cr\hfill  
111111100110110000\ 
0000100010100010\cr}$$

$$X_3 = \left(\matrix{
0.581172210661\ldots  \le  Re(L')  \le  0.581607219801\ldots\cr
1.156446469500\ldots  \le  Re(D')  \le  1.156834018585\ldots\cr
1.404200819866\ldots  \le  Re(R')  \le  1.404546086849\ldots\cr
-3.312214322575\ldots  \le  Im(L')  \le  -3.311906724662\ldots\cr
-2.756280098119\ldots  \le  Im(D')  \le  -2.755732020947\ldots\cr
-1.179687500000\ldots  \le Im(R')  \le  -1.179199218750\ldots\cr
}\right)$$

\noindent $X_3\ \ $ quasi-relators:

$r_1 = FFwfwFFwwFWFwFWfWFWffWFWfWFwFWFww$ 

$r_2 = FFwfwFwfWfwfWWfwfWfwFwfwFFwwFWFww$

$$\displaylines{X_{4a} = 
111000000001000110\ 
011001001111101010\ 
011110110110111101\hfill\cr\hfill  
100011111110110110\ 
10000111101\cr}$$

$$\displaylines{X_{4b} = 
111000000001000110\ 
011001001111101010\ 
111110010110011101\hfill\cr\hfill  
000011011110010110\ 
00000101101\cr}$$

$$X_4 = \left(\matrix{
0.333217001023\ldots  \le Re(L')  \le 0.334957037582\ldots\cr
0.977398792251\ldots  \le  Re(D')  \le 0.978173890421\ldots\cr
1.354137107330\ldots  \le  Re(R')  \le  1.354827641296\ldots\cr
-3.319596672476\ldots  \le  Im(L')  \le  -3.318981476651\ldots\cr
-2.825337821794\ldots  \le  Im(D')  \le  -2.824789744622\ldots\cr
-1.225585937500\ldots  \le  Im(R')  \le  -1.224609375000\ldots\cr
}\right)$$

\noindent $X_4\ \ $ quasi-relators:

$r_1 = FFwfwFwfWfwfWfwFwfwFFwwFWFwFWFww$

$r_2 = FFwfwFwfwFFwwFWFwFWfWFWfWFwFWFww$

$$\displaylines{X_{5a} = 
001000110001110111\ 
001111000101111111\ 
101111100111001111\hfill\cr\hfill 
000001111011110111\ 1\cr}$$

$$\displaylines{X_{5b} = 
001001110000110110\ 
001110000100111110\ 
101110100110001110\hfill\cr\hfill  
000000111010110110\ 1\cr}$$

$$X_5 = \left(\matrix{
-1.379848991182\ldots  \le  Re(L')  \le  -1.378108954623\ldots\cr
-1.379674742433\ldots  \le  Re(D')  \le  -1.376574349753\ldots\cr
0.999893182771\ldots  \le  Re(R')  \le  1.002655318635\ldots\cr
-2.537067582893\ldots  \le  Im(L')  \le  -2.534606799593\ldots\cr
-2.536501152136\ldots  \le  Im(D')  \le  -2.534308843448\ldots\cr
-0.001953125000\ldots  \le  Im(R')  \le  0.001953125000\ldots\cr
}\right)$$

\noindent $X_5\ \ $ quasi-relators:

$r_1 = FwFWFwFwfwfWfwfw$

$r_2 = FwfwfWfWFWfWfwfw$

$$\displaylines{X_{6a} = 
111000000001000111\ 
111111110101001111\ 
011111010111111111\hfill\cr\hfill  
110001001011000111\ 0\cr}$$

$$\displaylines{X_{6b} = 
111001000000000110\ 
111110110100001110\ 
011110010110111110\hfill\cr\hfill  
110000001010000110\ 0\cr}$$

$$X_6 = \left(\matrix{
1.378108954623\ldots  \le  Re(L')  \le  1.379848991182\ldots\cr
1.376574349753\ldots  \le  Re(D')  \le  1.379674742433\ldots\cr
0.999893182771\ldots  \le  Re(R')  \le  1.002655318635\ldots\cr
-2.537067582893\ldots  \le  Im(L')  \le  -2.534606799593\ldots\cr
-2.536501152136\ldots  \le  Im(D')  \le  -2.534308843448\ldots\cr
-0.001953125000\ldots  \le  Im(R')  \le  0.001953125000\ldots\cr
}\right)$$

\noindent $X_6\ \ $ quasi-relators:

$r_1 = FWFwFWfWFwFWFwfw$

$r_2 = FWFwfwFwfWfwFwfw$

\medskip \noindent {\bf  Proof:}  
Two
computer files contain the data needed for the proof.   The first computer file describes the partition of ${\cal W}$ into
sub-boxes, and the second describes the
killerwords and contradictions associated with each sub-box (other than the $X_i$).  We have a
computer program {\it verify}  which shows
that the killerwords in question actually do kill off their associated sub-
boxes.  This computer program addresses
the issues of Remark 1.21.  The code for {\it
verify} is given in Appendix 1, although we
encourage readers to produce their own verification programs. 

In addition, {\it verify} showed that the listed words were quasi-relators for the given sub-boxes. 
\qed

\medskip \noindent {\bf Corollary 1.29:} If $\delta$ is a shortest geodesic in $N,$ a closed orientable
hyperbolic 3-manifold, then 

i)  either tuberadius($\delta) > \ln(3)/2$ or
$\exp(\length(\delta))\in {\cal L}(X_k)$ for some $k\in {0,\ldots,6}$
 where ${\cal L}(X_k)$
denotes the range of $L^\prime$ values in the sub-box $X_k.\qed$

ii) either tuberadius($\delta) > \ln(3)/2$ or
${\rm tuberadius}(\delta) = {\rm Re}(D)/2$ where $\exp(D)
\in {\cal D}(X_k)$ for some $k\in {0,\ldots,6}$
 where ${\cal D}(X_k)$
denotes the range of $D^\prime$ values in the sub-box $X_k.\qed$

\medskip \noindent {\bf Experimental Theorem 1.30:} Associated to each of the sub-boxes $X_0, X_1, \ldots, X_6$ of ${\cal W}$ is a closed orientable 3-manifold of Heegaard genus 2 with fundamental group generated by $f,w$ with 2 relators $r_1,\ r_2$ as given in Proposition 1.28.

\medskip \noindent {\bf Experimental Proof:}  Experimentally, at some point in the sub-box $X_i$ under consideration the quasi-relators are actually relators.  Applying Berge [Be], it follows that the  2-generator 2-relator presentation $<f,w;\ r_1,r_2>$ is the group presentation of a Heegaard genus 2 closed orientable 3-manifold. \qed

\medskip \noindent {\bf Remarks 1.31:}  i) In Chapter 3 it is proven rigorously that Vol3, the hyperbolic 3-manifold with the third smallest volume known, is the unique manifold associated with sub-box $X_0.$  
As of this writing its not known if the other boxes have unique manifolds.

ii) Except for Vol3, Riley's program POINCAR\'E was the first to show (experimentally) that there is a closed orientable (hyperbolic) 3-manifold associated to each box.  It provided a group presentation, which presumably could have been shown to be those of the above.

         iii) Experimental evidence suggests that the manifolds associated to $X_5$ and $X_6$ are isometric
hyperbolic 3-manifolds, isometric to the Weeks census manifold $s479(-3,1).$
Also experimentally,   the manifold associated with $X_2$ is $s778(-3,1)$ and the manifold associated with $X_1$ is $v2678(2,1).$

        iv) Berge [Be] provides the explicit Heegaard genus 2 diagrams for
each of the manifolds in Experimental Theorem 1.30, so that the dilligent reader can recover Dehn
surgery descriptions of these manifolds.

\medskip \noindent {\bf  Remark 1.32:}  It might be best to say that the list of killerwords and
contradictions used to prove Proposition 1.28 was generated in an artistic
rather than systematic mathematical way.  Nevertheless, to have a rigorous
mathematical proof one need only prove that
these words work.

\medskip \noindent {\bf  Remark 1.33:}  Recall that we use first-order Taylor approximations (with Remainder term) which we
denote {\it AffApprox's} to show that a killer word which eliminates a point $x\in {\cal W}_i$
eliminates all of ${\cal W}_i.$  In setting up a Taylor
approximation  system we avoided a considerable amount of calculation by
using the parameter space ${\cal W}$ instead of working with ${\cal P}.$

In the parameter space ${\cal W},$ all functions analyzed via Taylor
approximations are built up from the
operations $ +,\ -,\ \times, \  /,\ \sqrt{} \ .$  We prove combination
formulas for these functions, which show how the Taylor approximations
(including the Remainder term) change when one of these operations is
applied to two AffApprox's.   This is carried out in Chapter 6.

To ensure that all of our computer calculations are rigorous, we use a round-off error analysis.  Typically, this is done by using interval
arithmetic on floating-point numbers. This is too slow for our purposes,
and so we introduce round-off error at the level of AffApprox's and
incorporate the round-off error into the Remainder term.  This also
requires developing round-off error for complex numbers.  This is all
done in Chapters 7 and 8.

\sectionskip
\magnification=1200
\baselineskip = 0.2 true in
\def\qed{{\bf \ \ \ \ \vrule height6pt width5pt depth1pt}}
\def\Relength{{\rm Relength}}
\def\length{{\rm length}}
\def\Arccosh{{\rm Arccosh}}
\def\trace{{\rm trace}}
\def\distance{{\rm distance}}
\def\Redistance{{\rm Redistance}}
\def\va{{\rm visualangle}}

\noindent {\bf Chapter 2:  The Corona Insulator Family}
\medskip
The upshot of Proposition 1.28 is that if a closed orientable hyperbolic 3-manifold has  a shortest geodesic which does not have an embedded $\ln3/2$ tube then the parameters for its associated 2-generator subgroup(s) $(G,f,w)$  must be in one of the sub-boxes $X_0, X_1, \ldots, X_6$ listed in that proposition.  
Nonetheless---as we shall see in this Chapter and the next--- such manifolds have non-coalescable insulator families about their shortest geodesics.  However, they might not be Dirichlet insulator families.

In this Chapter, we describe a new insulator  family $\{\kappa_{ij}\}$  called the Corona insulator, and   
we describe a condition sufficient for this family to be
non-coalescable---a condition which is weaker than the
${\rm tuberadius}(\delta) > \ln(3)/2$ sufficient condition for the Dirichlet insulator family.

The reason Dirichlet insulator families for geodesics with solid tubes of radius greater than $\ln3/2$ are non-coalescable is that the amount of visual angle taken up by the various insulators is less than 120 degrees, and thus there is no chance for tri-linking to occur.   The visual angle (measured at one axis) for a member of the Dirichlet insulator family associated to two axes depends only on the real distance between the two axes.  

In contrast, the visual angle for a member of the Corona insulator family associated to two axes depends on the (complex) distance between the two axes.  We now give this function,
${\cal C},$ and name it the {\it Corona} function.  After that we give a precise definition of the visual angle function, and prove that the Corona function is the proper visual angle function for the Corona insulator family.

 \medskip \noindent {\bf Definition 2.1:} Let
${\cal C}: (0,\infty)\times S^1\to {\bf R}$ be defined by

$${\cal C}(u,v) = |{\rm Im}
(\Arccosh(1-{4 \over 1 \pm \cosh(u+iv)} ))|$$  
\noindent where $\pm$ is positive for $-\pi/2 \le v \le \pi/2$ and negative otherwise. 

\midinsert
\centerline{\psfig{file=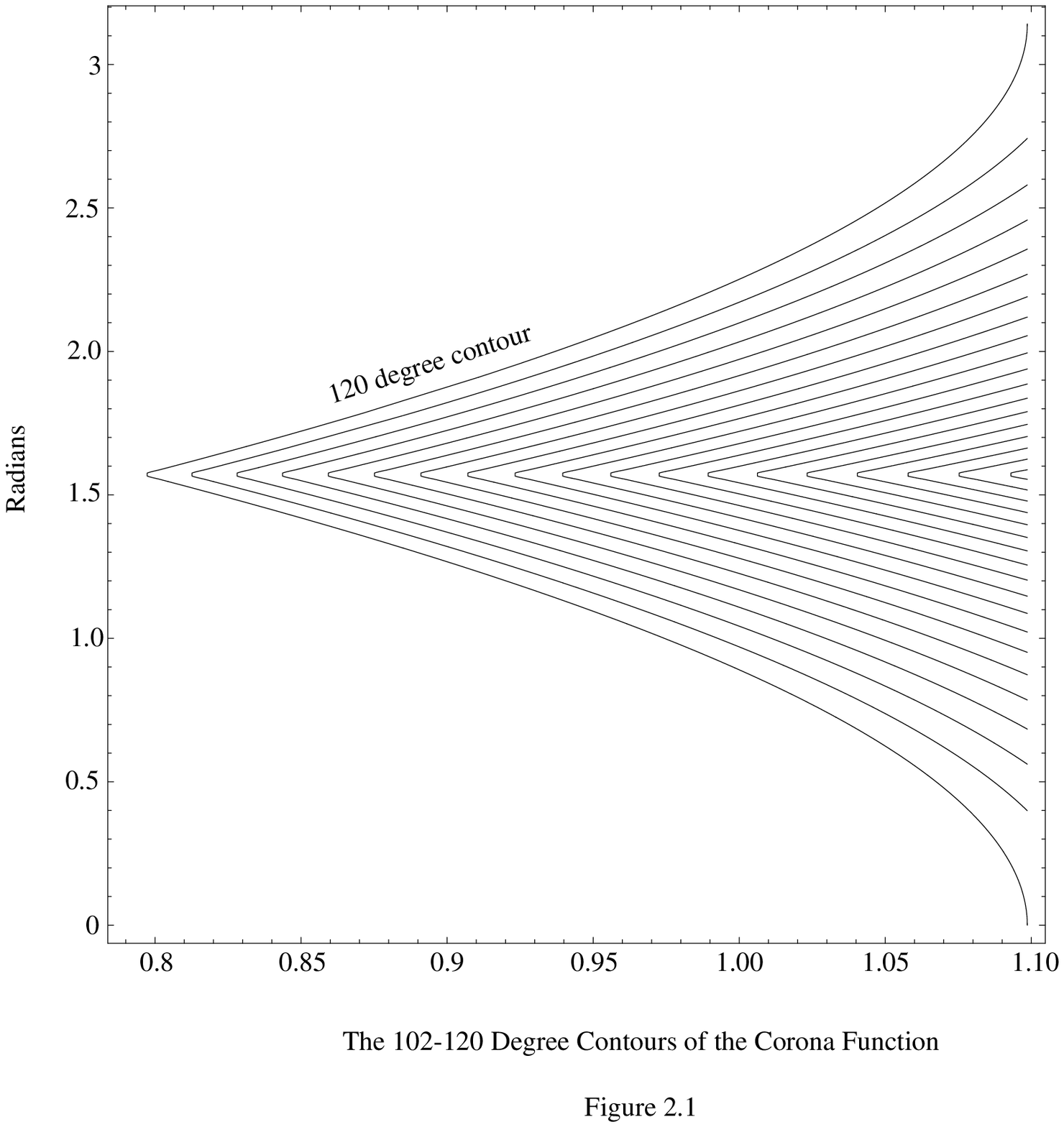,height=10cm}}
\endinsert  

In the following definition, it is helpful to imagine the geodesic $\sigma$ as being the $z-$axis in the upper-half-space model of ${\bf H^3}.$ 
\medskip \noindent {\bf Definition 2.2:}  If $\sigma\subset {\bf H^3}$ is a geodesic, then $S^2_\infty-\partial\sigma$ is
parametrized by $S^1\times {\bf R}$ (these are sometimes called Steiner circles) where each $x\times {\bf R}$ lies in the ideal boundary of
a hyperbolic halfplane bounded by $\sigma,$ two such lines $x\times {\bf R},\ y\times {\bf R}$
are at distance $\theta$ in the $S^1$ factor if they meet at $\partial \sigma$ at
angle $\theta.$  If $R \subset S^2-\partial \sigma,$ then define $\va_\sigma(R)=\theta\in [0,2\pi],$ where $\theta=\inf \{\theta_2-\theta_1 \bmod 2\pi\mid R
\subset [\theta_1,\theta_2]\times {\bf R}\}.$  
The possible choice of 0 or $2\pi$ is made
in the obvious manner.

\medskip \noindent {\bf Proposition 2.3:}  Let $\delta_i, \delta_j$ be disjoint oriented geodesics in ${\bf H^3}$.
Then there exists a smooth simple closed curve $\kappa_{ij}$ in $S^2_{\infty}$
separating $\partial\delta_i$ from $\partial\delta_j$ such that for $k\in \{i,j\},\ 
\va_{\delta_k}(\kappa_{ij})={\cal C}(\distance(\delta_i, \delta_j)).$  

\medskip \noindent {\bf Proof:} 
Let $P$ be the orthocurve from  $\delta_i$ to  $\delta_j.$
Consider the half-plane with boundary $\delta_i$ determined by $\delta_i$ and $P,$ and the half-plane with boundary $\delta_j$ determined by $\delta_j$ and $P.$   Allow these half-planes to expand into {\it solid angles} at the same rate. (A solid angle is a closed set in $B^3 = {\bf H^3} \cup S^2_\infty$ bounded by two hyperboic halfplanes which meet along a common geodesic.)  At first, the four half-planes that bound these solid angles intersect in ${\bf H^3},$ but at some angle $\theta$ these half-planes intersect only at infinity (that is, in $S_\infty^2$). By reasons of symmetry they intersect in two or four points (four points of intersection occur when ${\rm Im}(\distance(\delta_i, \delta_j))$ is  $\pi/2$ or $3 \pi/2$).

Let $S_i, S_j$ be the solid angles which exist at angle $\theta.$ Let
$T_k=S_k\cap S^2_\infty$ for $k\in \{i,j\}.$  Let $\kappa_{ij}$ be a simple closed curve in
$T_i\cap T_j$ which separates $\partial \delta_i$ from $\partial \delta_j.$  By
construction, for $k\in \{i,j\}\ \va_{\delta_k}(\kappa_{ij})=\theta.$ See
Figure 2.2.

\midinsert
\centerline{\psfig{file=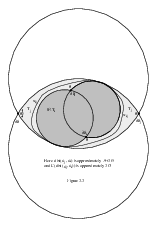,height=12cm}}
\endinsert 

To complete the proof of the Proposition, we now show that $\theta = {\cal C}(\distance(\delta_i,\delta_j)).$
To do this, we use hyperbolic trigonometry on a degenerate right-angled hexagon in ${\bf H^3}.$  
 Following [F], a degenerate right-angled
hexagon is a 5-tuple of oriented geodesics $S_1,\cdots, S_5$ in ${\bf H^3}$ such that
$S_i$ is orthogonal to $S_{i+1}$ and $S_1$ and $S_5$ limit at a common point $S_0$ at
infinity.   These oriented geodesics give rise to complex numbers 
$\sigma_0, \sigma_2,\sigma_3, \sigma_4.$   $\sigma_0 = 0$ if the axes $S_1$ and $S_5$ either both point
 into $S_0$ or both point out of $S_0.$  Otherwise $\sigma_0=\pi i.$  For
$k\in \{2,3,4\} \ \sigma_k=e$ if an $e$-translation of the oriented geodesic $S_k$
takes the oriented geodesic $S_{k-1}$ to the oriented geodesic $S_{k+1}.$  By [F; pg. 83] we have the following Hyperbolic Law of Cosines.
$$\cosh(\sigma_0) = \cosh(\sigma_{2}) \cosh(\sigma_{4})
+ \sinh(\sigma_{2}) \sinh(\sigma_{4}) \cosh(\sigma_{3}).$$

We work in the upper-half-space model of hyperbolic 3-space, and normalize so that the ortholine from $\delta_j$ to $\delta_i$ is B (thus $\delta_i$ intersects $B$ above $\delta_j$), while the oriented axis $\delta_i$ is $C.$ $B$ will be $S_3,$ while the oriented geodesics $\delta_i$ and $\delta_j$ will be $S_2$ and $S_4$, respectively.  

Of course, $u+iv = \distance(\delta_i,\delta_j).$
  If $-\pi/2 \le v \le 0$ then the intersection points at infinity occur in the second quadrant and the fourth quadrant (see Figure 2.3a).  For convenience, we work with the point in the second quadrant and send (unique) perpendiculars from it to the geodesics $\delta_i$ and $\delta_j$.  The perpendicular to $\delta_i$ will be oriented towards $\delta_i$ and then denoted $S_1$, while the perpendicular to $\delta_j$ will be oriented away from $\delta_j$ and then denoted $S_5.$  The  intersection point at infinity (in the second quadrant) is  $S_0.$

\midinsert
\centerline{\psfig{file=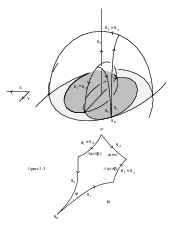,height=12cm}}
\endinsert  

This is the proper set-up for applying the (degenerate) Hyperbolic Law of Cosines (see Figure 2.3b).  Note that $\sigma_3 = -(u+iv),$ and $\sigma_0 = i \pi.$ By symmetry $\sigma_2 = \sigma_4 = (\alpha + i \beta)/2$ where $(\alpha + i \beta)/2$ is $\distance(S_1, S_3).$  Plugging into the Law of Cosines, using a half-angle formula ($\cosh(2z) = 2 \cosh^2(z) - 1 = 2 \sinh^2(z) + 1$), solving for $\cosh(\alpha + i \beta),$ and taking the \Arccosh, we get the desired result.  Note that the visual angle in this set-up is $-\beta,$ thus necessitating taking the absolute value.

When $0 \le v \le \pi/2$ our 2 intersection points occur in the first and third quadrants, and we carry out the same procedure.  This
time $S_2$ and $S_4$ are traversed in the direction opposite to their orientations (the attendant changes in sign drop out though). In this case,  the visual angle is $\beta.$

The cases $-\pi \le v \le -\pi/2$ and $\pi/2 \le v \le \pi$ reduce to the previous cases after adding or subtracting $\pi.$  The formula in Definition 2.1 is then obtained after noting that $\cosh(z \pm i\pi) = - \cosh(z). \qed$ 

\medskip \noindent {\bf Definition 2.4:}  Let $\delta$ be a simple closed geodesic in the closed
orientable hyperbolic 3-manifold $N.$  Let $\{\delta_i\}_{i\ge 0}$ be the lifts of
$\delta$ to ${\bf H^3}.$  For each $\pi_1(N)$-orbit of unordered pairs $(\delta_i,
\delta_j)$ choose a representative where $i=0.$ 
If $\Redistance(\delta_0, \delta_i)\le \ln(3)/2,$ then let
$\kappa_{0j}$ be a smooth simple closed
curve in $S^2$ separating $\partial \delta_0$ from $\partial \delta_j$ such that for
$k\in \{0,j\},\ \  \va_{\delta_k}(\kappa_{0j})={\cal C} (\distance(\delta_0, \delta_j)).$
 If
$\Redistance(\delta_0,\delta_i) > \ln(3)/2,$ then let $\kappa_{0j}$ be the Dirichlet
insulater, i.e. the boundary of the geodesic plane orthogonally bisecting
the orthocurve between $\delta_0, \delta_j.$ 
 
In either case, extend the collection $\pi_1(N)$-equivariantly to a family $\{\kappa_{ij}\}$ defined for all
$i,j.$  This is the Corona family for $\delta.$

\medskip \noindent {\bf Lemma 2.5:} The Corona family $\{\kappa_{ij}\}$ is

i)  an insulator family for $\delta$

ii)  noncoalescable if $\max\{{\cal C}(\delta_0, \delta_j)\mid j>0\}<2\pi/3.$

\medskip \noindent {\bf Proof:} We check that $\{\kappa_{ij}\}$ satisfies the various conditions of
Definitions 0.4-0.5 of [G].

i)  By construction, $\kappa_{ij}$ separates $\partial \delta_i$
from $\partial \delta_j$ and $\{\kappa_{ij}\}$ is $\pi_1(N)$-equivariant.  Because for $k\in \{i,j\},\ $
$\delta_k-{\rm visual angle}(\kappa_{ij})<\pi,\ \  \{\kappa_{ij}\}$ satisfies the convexity
condition.  

Modulo the natural action of $\pi_1(N)$ on  $\kappa_{ij},$ 
there are only finitely many insulators $\kappa_{ij}$ which are not Dirichlet insulators.  Therefore, for fixed $i,$ there exist
only finitely many $\kappa_{ij}$ such that ${\rm diam}(\kappa_{ij})>\epsilon.$  This
establishes local finiteness.\qed

ii) No trilinking follows immediately from ii).\qed

\medskip \noindent {\bf Definition 2.6:} If $\delta$ is a simple closed geodesic in the hyperbolic
3-manifold $N,$ define ${\rm maxcorona}(\delta)=\max\{{\cal C}(\delta_0, \delta_j)\mid j>0\}$

\medskip \noindent {\bf Remark 2.7:} a)  It seems possible that the Dirichlet insulator family
associated to the geodesic $\delta\in N$ may be non-coalescable, while the
Corona insulator family is coalescable and conversely.

b) If ${\rm tuberadius}(\delta)>\ln(3)/2,$ then both the Dirichlet insulator and Corona
insulator families coincide and by Lemma 5.9 of [G] they are non-coalescable.

\medskip \noindent {\bf Proposition 2.8:} Let $\delta$ be a shortest geodesic in the closed orientable
hyperbolic 3-manifold $N.$  Then either the Corona insulator family is
noncoalescable or there exists a 2-generator subgroup $G$ of $\pi_1(N)$
with generators $f,w$ such that the parameter associated to
$\{G,f,w\}$ lies in the sub-box  $X_0= X_{0a} \cup X_{0b} \subset {\cal W}.$

\medskip\noindent {\bf Proof:}  Let $\delta$ be a shortest geodesic in $N,$ a closed orientable hyperbolic
3-manifold.  By Corollary 1.29 and Lemma 2.5 either $\delta$ has a non-coalescable insulator family 
or $\length(\delta)$ lies in ${\cal L}(X_k)$ for some $k\in
\{0,1,\ldots, 6\}$ and ${\rm maxcorona}(\delta)\ge 2\pi/3.$  Let $G$ be a 2-parameter
subgroup of $\pi_1(N)$ generated by $f$ and $g,$ where $f$ corresponds to $\delta$
and has axis $B$ and $g$ maximizes ${\cal C}(g(B),B).$  We will show that if
${\cal C}(g(B),B)\ge 2\pi/3,$ then the parameter associated to $\{G,f,g\}$ lies in the
sub-box $X_0.$

 In fact we will show that if $\{G,f,g\}$ is any torsion-free subgroup of
${\rm Isom}({\bf H^3})$ where $f$ is a length-minimizing loxodromic element with axis $B$ and
$g$ maximizes ${\cal C}(g(B),B)$ where ${\cal C}(g(B),B)\ge 2\pi/3,$ then the parameter
$\beta=(L_\beta, D_\beta, R_\beta)\in \exp^{-1}(X_0).$

It follows as in the first paragraph that $L_\beta\in {\cal L}(X_k)$ for some $k\in
\{0,1,2,3,4,5,6\}.$    Also $D_\beta$ is subject to nontrivial constraint.  For
example if $L_\beta\in {\cal L}(X_1),$ then $D_\beta$ must lie in the decorated region of Figure 2.4, because $\Redistance(g(B),B)\ge {\rm min}_d(X_k)> 1.059$ where
${\rm min}_d(X_k)$ is the minimal $d$ value in $X_k.$  
In fact the disjointness of the
$L(X_k)$'s implies that if $\beta\in S \cap X_k,$ then 
tuberadius$(\delta)\ge {\rm min}_d(X_k)/2,$ 
where $\delta$ corresponds to the element $f.$
Finally ${\cal C}(g(B),B)\ge 2\pi/3$ implies that $D_\beta$ cannot lie
in the decorated region of Figure 2.4.

\midinsert
\centerline{\psfig{file=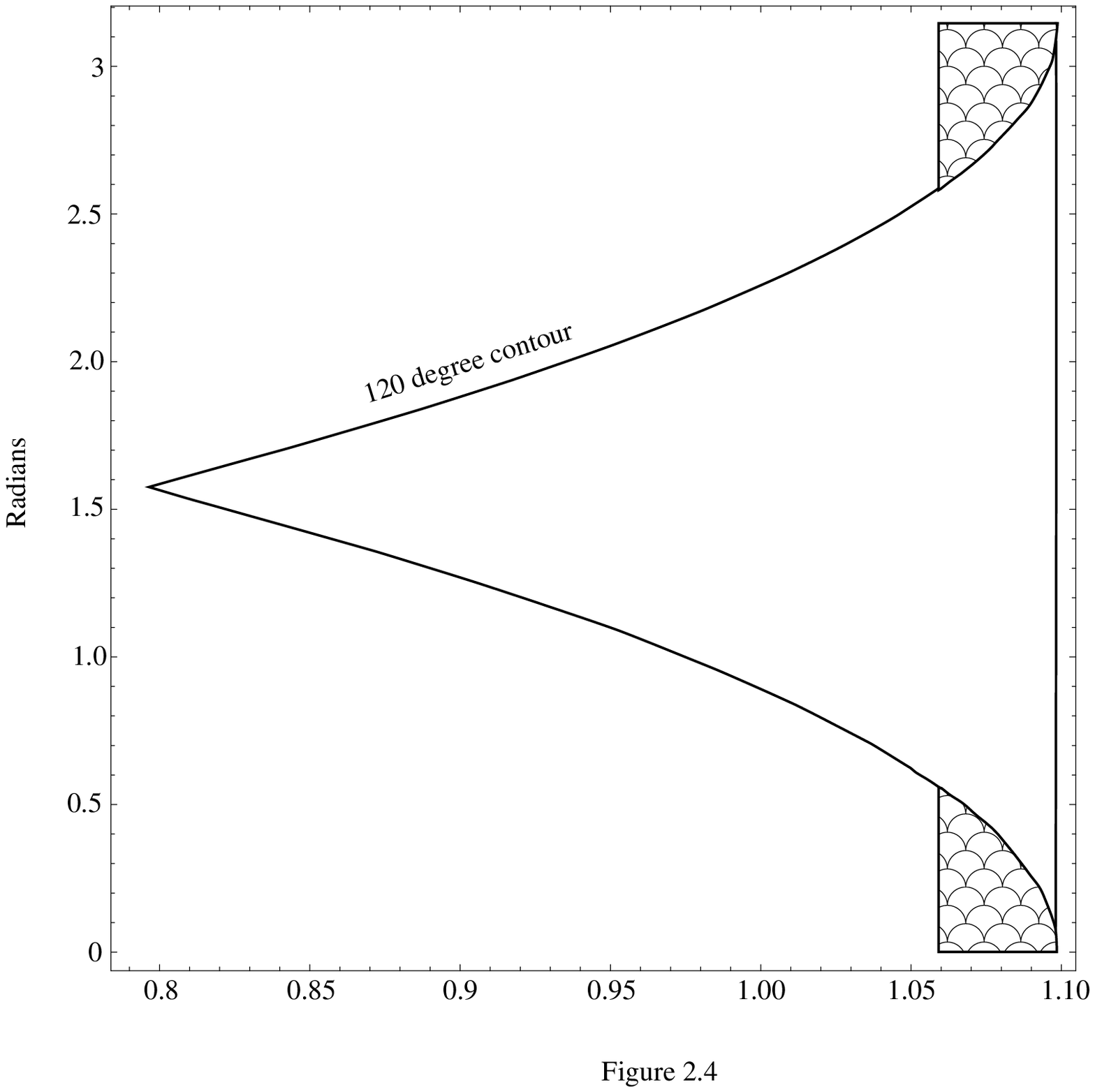,height=8cm}}
\endinsert

Our proof is now similar to the proof of Proposition 1.28.
We partition the
initial box ${\cal W}$ into sub-boxes ${\cal W}_i$ and eliminate ${\cal W}_i$ if
any of the following conditions hold (for clarity, the conditions have been translated into ``pre-exponentiated" form).

a)  There exists no $\beta\in {\cal W}_i$ such that $\length(f_\beta) \in {\cal L}(X_k).$ 

b)       ${\cal W}_i$ has some  $L$ values in ${\cal L}(X_k)$ but
there exists no $\beta\in {\cal W}_i$ such that distance $(w_\beta(B),B)\in\ $ decorated region for that $k.$

c)  There exists a killerword $h$ in $f,w,f^{-1},w^{-1},$ such that
$\Relength(h_\beta) < \Relength(f_\beta)$  and $h_\beta \neq {\rm id}$
for all $\beta\in {\cal W}_i.$

d)  There exists a killer word $h$ in $f,w,f^{-1},w^{-1}$ such that
$${\cal C}(\distance(h_\beta(B),B)) >{\cal C}(\distance(w_\beta(B),B))$$
\noindent and $h_\beta(B) \neq B$ for all $\beta\in {\cal W}_i.$ 

We have two files that contain the decomposition of ${\cal W}$ into  sub-boxes and associated conditions/(killer words).  The program {\it fudging} checks that these files do indeed eliminate all of 
${\cal W} - X_0.\ $  {\it Fudging} analyzes the cases $k = 1,\ldots, 6$ all at once. \qed

\medskip\noindent {\bf Remarks 2.9:}  i)  Note that in practice the proof of Proposition 2.8 requires working in a
considerably smaller parameter space than that of Proposition 1.28. Condition a)
implies that the parameter space is ``$(2+\epsilon)$-complex dimensional" and
condition b) implies that one of these parameters is greatly constrained.   This
suggests why it took so much longer to come up with the partition and the associated killer words for Proposition 1.28.  
In fact, it took roughly 1500 CPU days to find the partition and the associated killer words for Proposition 1.28, versus roughly 2 CPU days for Proposition 2.8.  Here, the term ``CPU day''  refers to 24 hours of running an SGI Indigo 2 workstation with an R4400 chip, and the estimate of 1500 CPU days refers to 15 to 20 such machines running 80 to 90 percent of the time for 3 to 4 months.

ii) We took pains to make {\it fudging} as similar to {\it verify} as we could, thereby lessening the amount of analysis needed to show the veracity of {\it fudging}.   Appendix 2 contains those sections of {\it fudging} that differ from corresponding sections of {\it verify}.

iii)  When working with exponentiated co-ordinates (that is, in ${\cal W}$ rather than ${\cal P}$) the Corona function changes as follows.  Let $X = \exp(\alpha + i \beta)$ and 
$ U = \exp(u + iv),$  then the Corona function formula 
$$\cosh(\alpha + i \beta) = 1 - {4 \over 1 + \cosh(u+iv)}$$
becomes
$${X + X^{-1} \over 2} = 1 - {4 \over 1 + (U + U^{-1})/2}$$

It is a pleasant exercise to solve this, and we find that 
$$X = {(U^2 -6U + 1) \pm 4(U - 1)\sqrt{-U} \over (U+1)^2}$$
the two answers so gotten are reciprocals, which means their associated arguments are opposites.  
In {\it fudging}, the exponentiated version of the Corona function is the function {\it horizon(ortho)}, which  takes in $U = ${\it ortho} and computes the associated $X$ value.  $\beta,$ the argument of $X,$ is implicitly gotten in the function {\it larger-angle}. 

iv) It is possible that by working purely in the context of the Corona
function, rather than first working with $\Redistance$ and attempting to prove
Proposition 1.28, the computer proof can be simplified.   We started this
project with the naive idea that perhaps Vol3 was the only manifold whose
shortest geodesic did not have a $\ln(3)/2$ tube.  The remarkable fact that this
naive idea is almost correct accounts for the fact that a proof of Theorem 0.2
can be obtained with only the mild extra effort detailed in this chapter and the next.

For $0 \le v \le \pi/2\ \ \ {\cal C}(u,v) = \beta$ where 
$\beta$ satisfies 
$$\cosh(\alpha + i \beta) = 1 - {4 \over 1 + \cosh(u + iv)}.$$

For $-\pi/2 \le v \le 0$ ${\cal C}(u,v) = -\beta$ where 
$\beta$ satisfies 
$$\cosh(\alpha + i \beta) = 1 - {4 \over 1 + \cosh(u + iv)}.$$

For $\pi/2 \le v \le \pi$ ${\cal C}(u,v) = -\beta$ where 
$\beta$ satisfies 
$$\cosh(\alpha + i \beta) = 1 - {4 \over 1 + \cosh(u+i(v - \pi))}.$$

For $-\pi \le v \le -\pi/2$ ${\cal C}(u,v) = \beta$ where 
$\beta$ satisfies 
$$\cosh(\alpha + i \beta) = 1 - {4 \over 1 + \cosh(u+i(v + \pi))}.$$
In all of these cases, we normalize the sign of $\beta$ by requiring that $\alpha$ be positive.



%
%
\medskip \noindent {\bf Lemma 2.10:} 
i) if $h =wFwfwfWfwf$   then  ${\rm Relength}(h_\alpha) < {\rm Relength}(f_\alpha)$ for all parameters $\alpha$ in $\exp^{-1}(X_0)$

ii) if  $h_2 = wFwfwwfwFw$ then  ${\rm Relength}(h_\alpha) < {\rm Relength}(f_\alpha)$ for all parameters $\alpha$ in $\exp^{-1}(X_0)$

iii) $|R| < |L/2|$ throughout $\exp^{-1}(X_0).$
\medskip
\noindent {\bf Proof:} This easy calculation was carried out by {\it verify}. \qed

\sectionskip\magnification=1200
\baselineskip = 0.2 true in
\def\qed{{\bf \ \ \ \ \vrule height6pt width5pt depth1pt}}
\def\distance{{\rm distance}}
\def\length{{\rm length}}

\noindent {\bf Chapter 3: Vol3}
\medskip
\noindent The two main results of this chapter are

\medskip \noindent {\bf Proposition 3.1:}  If $N$ is a closed orientable hyperbolic 3-manifold, then either $N$ is 
Vol3,  the third smallest known closed orientable hyperbolic 3-manifold, or maxcorona$(\delta) < 2 \pi/3$ for $\delta$ a shortest geodesic in $N.$

\medskip \noindent {\bf Proposition 3.2:} Any shortest geodesic in  Vol3 satisfies the insulator condition.

\medskip \noindent {\bf Remark 3.3:}  
Vol3 is the third smallest known closed orientable
hyperbolic 3-manifold.  Topologically Vol3 is (3,1)  surgery on manifold m007 in the census of
cusped hyperbolic 3-manifolds (see [W]).  It is also (-3,2) (-6,1) surgery on the ``left-handed Whitehead link", 
link $5_2^2$ in the standard knot tables.  Snappea gives an
experimental proof that Vol3 is hyperbolic and that its volume is that of the
regular ideal 3-simplex.  A rigorous proof can be found in [JR].  Previously,
Hodgson-Weeks [HW1] had found an exact Dirichlet domain for Vol3, that is,  the face
pairings were expressible as explicit matrices with coefficients in a finite
extension $F$ of ${\bf Q}$ and they obtained equations in $F$ for the various faces.  See
Remark 3.14.

\medskip \noindent {\bf Remarks 3.4:}
i) Idea of Proof of Proposition 3.1: If $N$ has maxcorona$(\delta) \ge 2\pi/3$, then it must have an $(L,D,R)$ parameter in the region ${\cal R} = \exp^{-1}(X_0),$ ($X_0$ is defined in Proposition 1.28).   A geometric argument (Lemma 3.7) which utilizes Lemmas 3.5 and 3.6 shows that $R = 0$ and an algebraic argument (Lemmas 3.8, 3.9, 3.11, 3.12) shows that $L = D = \omega,$ where $\exp(\omega)$ is a root of the polynomial $1 + 2d + 6d^2 + 2d^3 + d^4.$  This implies ([HW1] or [JR]) that $g = \pi_1$(Vol3) and hence $N$ is covered by Vol3.  By [JR], $N = {\rm Vol3}.$

 ii)  The proof of Proposition 3.2 will be somewhat similar to
the proofs of Propositions 1.28, 2.8.

\medskip \noindent {\bf Lemma 3.5:}  Let ${\cal R} = \exp^{-1}(X_0)$ If
$\alpha= (L,D,R)=(l+it, d+ib, r+ia)  \in {\cal R} \cap {\cal T},$ then $f_\alpha,
w_\alpha$ satisfy the relations

i)  $wFwfwfWfwf$

ii) $wFwfwwfwFw$

\medskip \noindent {\bf Proof:}
  In i), ii) above and what follows below we supress the subscripts
$\alpha.$  Also $W$ (resp. $F$) denotes $w^{-1}$ (resp. $f^{-1}$).  because i), ii) are cyclic permutations of the quasi-relators $r_2,\ r_1$ corresponding to the $X_0$ sub-box of Proposition 1.28, it follows that if 
$h=wFwfwfWfwf$     or $wFwfwwfwFw,$ then ${\rm Relength}(h) < {\rm Relength}(f)$ throughout ${\cal R}.$  Since
$\alpha\in {\cal T},$  $f$ is a
shortest element and so $h={\rm id}.$  \qed

\medskip \noindent {\bf Lemma 3.6:}  The following substitutions give rise to three sets of new relators.

a) In i), ii) replace $f$ by $w$  and replace $w$ by $f$.

b) In i), ii) replace $f$ by $F$.

c) In i), ii) replace $w$ by $W$.

\medskip \noindent {\bf Proof:}  Note that replacing $f$ by $w$ necessitates replacing $F$ by $W$, and so on.  

a)  First, a cyclic permutation of relator i) gives relator i) with $f$
replaced by
$w$ and $w$ replaced by $f$.  Second, one readily obtains the relator $fWfwffwfWf$ from i), ii), because $fWfwf = (wFwfw)^{-1} = wfwFw = (fwfWf)^{-1}$ where the first and third equalities follow from i) and the second from ii).

b)  Again it is routine to obtain b) from i), ii).

c)  Conclusion c) follows from a) and b). \qed

\medskip \noindent {\bf Lemma 3.7:} If $(L,D,R) \in \cal R \cap {\cal T},$ then $R=0$.

\medskip \noindent {\bf Proof:}  With respect to the parameter $(L, D, R)$ let $G$ be the group
generated by $f, w$.  Figure 3.1  shows a schematic picture of geodesics
$B, W(B), w(B), f(w(B)), f(W(B)). $ Also, it shows the images of the orthocurve $O$ from
$W(B)$  to $B$ after translation by $w, f$ and $fw.$  Finally, it shows the
orthocurves  $O_1$ from $fW(B)$ to $W(B)$ and $O_2$ from $w(B)$ to $fw(B).$ Note that
Figure 3.1 displays the situation where ${\rm Re}(R)>0.$  It is also apriori
possible that $O_2$  might intersect $w(B)$ on the positive side of $w(O).$
There are other similar  possible inaccuracies.

\midinsert
\centerline{\psfig{file=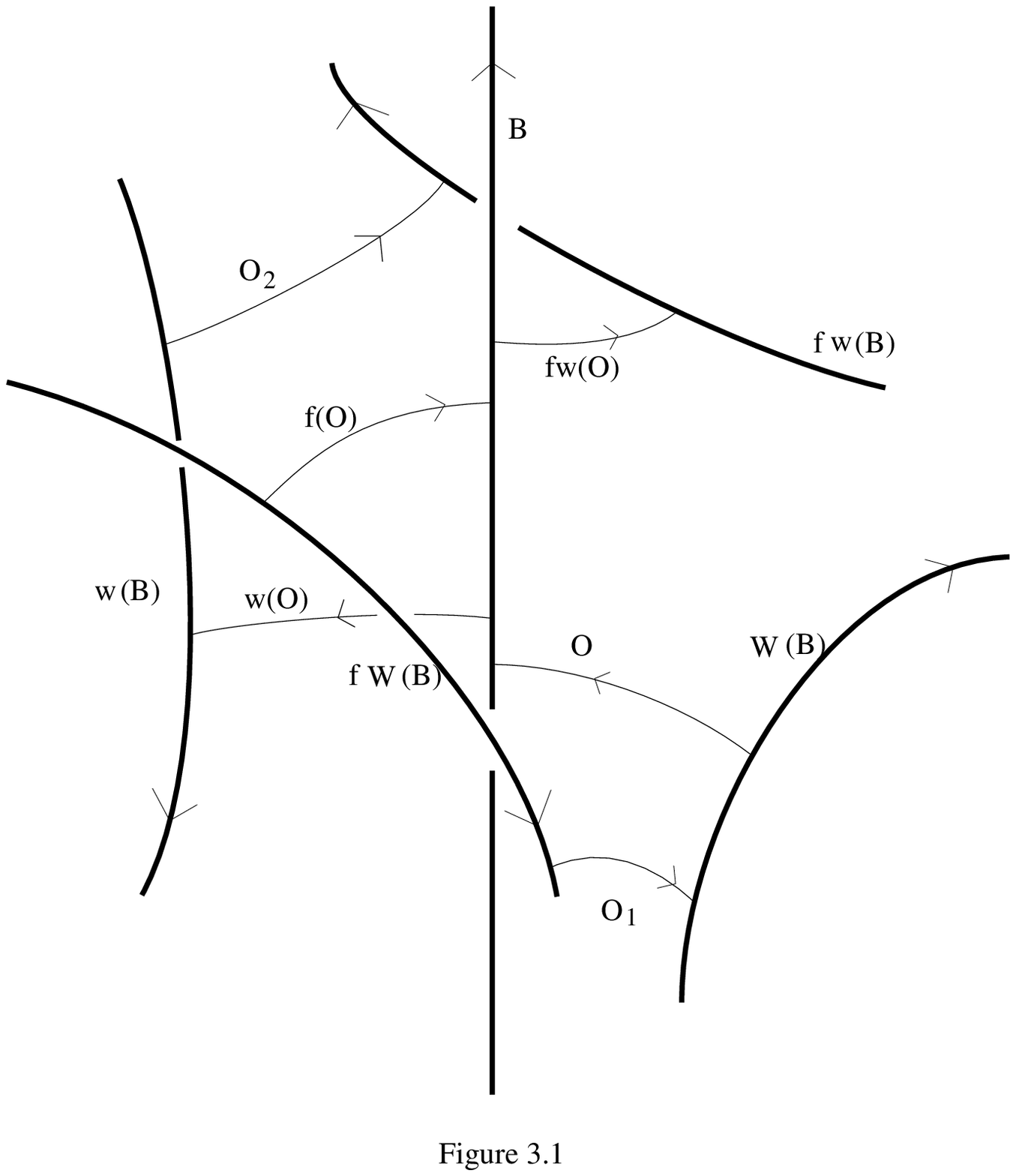,height=10cm}}
\endinsert 

$\sigma_1=fwffw\in G$ sends
geodesic $W(B)$ to $fw(B)$ and $\sigma_2=wFFwF\in G$ sends the
geodesic $fW(B)$ to $w(B)$.  
Now $\sigma_2^{-1}\sigma_1=fWffWfwffw$ is a relator
of $G$, since it is a cyclic permutation of relator ii) with $f$ replaced by $w$
and $w$ replaced by $f.$  Thus $\sigma_1=\sigma_2$ and hence
$\sigma_1(O_1)=O_2.$

As in Fenchel [F] a right angled hexagon consists of a cyclically
ordered 6-tuple of oriented geodesics $\lambda_1,\cdots,\lambda_6$
in ${\bf H^3}$ such that $\lambda_i$ intersects $\lambda_{i+1} \pmod 6$
orthogonally.  Each ``edge" of the hexagon is labelled by the
complex number $e_i$ where one obtains $\lambda_{i+1}$ from $\lambda_{i-1}$
by $e_i$ translation of ${\bf H^3}$ along the oriented geodesic
$\lambda_i.$

(3.1) The effect of reversing the orientation of
$\lambda_i$ is to change $e_i$ to $-e_i, e_{i+1}$ to $e_{i+1} + \pi i,$ and
$e_{i-1}$ to $e_{i-1}+\pi i.$

Figure 3.1 gives rise to the two right-angled hexagons drawn in Figure 3.2.  (Figure 3.2a may be inaccurate for
the following reason.  It is not clear
 whether the head of $O_1$ should be placed in  front of the tail
of $O$ or behind the tail of $O.$  A similar statement holds for
the tail of $O_1$ and for Figure 3.2b.) Assume that in Figure
3.2a, $\lambda_1$  corresponds to $B$ and the edges are
cyclically  ordered counterclockwise.  Then $e_6=D,
e_1=L, e_2=-D.$     We now show that if $e_5$ has value $c,$ then
$e_3$ has value $c+\pi i.$ Observe that there is an order-2
rotation $\tau$ of ${\bf H^3}$ about an axis orthogonal to $B$ which
reverses the orientation on $B$ and takes the oriented
orthocurve $O$ to the oriented orthocurve $f(O).$  Since 
distance($W(B),B)=\ $distance$(fW(B),B)$ it follows that
$\tau(fW(B))=-W(B)$ and $\tau(W(B))=-fW(B)$ where the - sign
indicates that the orientation has been reversed.   This in
turn imples   that $\tau(O_1)=-O_1$ and therefore using (3.1)
that  $e_3=e_5+\pi i.$

\midinsert
\centerline{\psfig{file=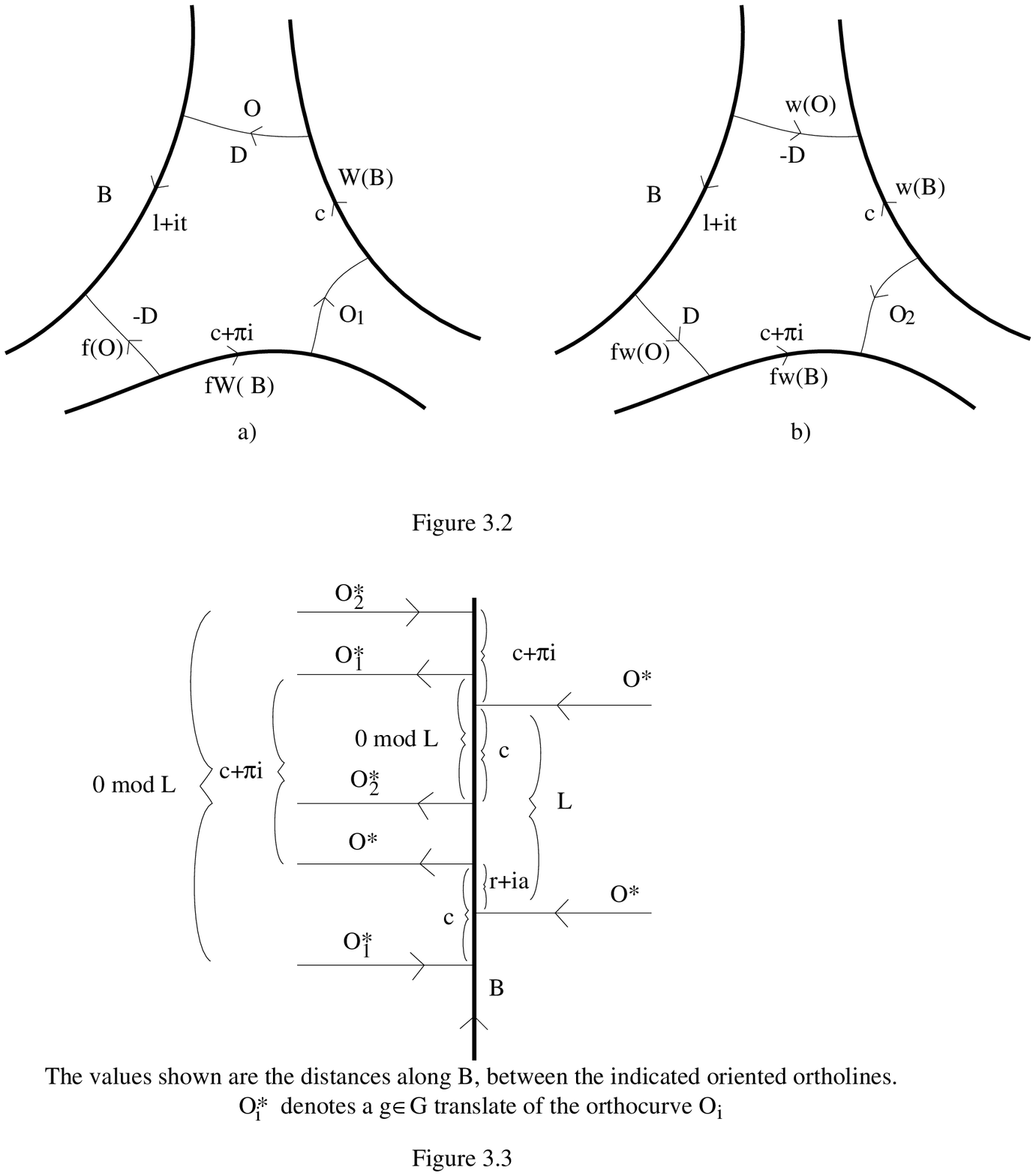,height=15cm}}
\endinsert 

Let $\phi$ be the isometry of ${\bf H^3}$ which is an  $r+i(\pi+a)$
translation of $B.$  Thus $\phi(B)=B, \ \phi(O) = -w(O),\ 
\phi(f(O)) = -fw(O).$  This implies that $\phi(W(B))=w(B)$
and $\phi(fW(B))=fw(B)$ which in turn implies that
$\phi(O_1)=-O_2.$  If the hexagon of Figure 3.2b with
edges $\lambda^\prime_1,\cdots,\lambda^\prime_6$ is
counterclockwise cyclically oriented so that 
$\lambda^\prime_1$ denotes the oriented geodesic $B$, then again
using (3.1) it follows that $e_5^\prime=c$ and
$e_3^\prime=c+\pi i.$

Via elements of $G,$ we translate each of $W(B), \ fW(B), \ w(B),
\ fw(B)$ to $B$ and after translation we obtain from Figure 3.2
the various distance relations schematically indicated on
Figure 3.3.  Here $O_i^*$ is $G$ translation of $O_i$ with an
endpoint on $B.$  Two such translates appear in Figure
3.3, one where $O_i^*$ points into $B$ and one where
$O_i^*$ points out.  Call the former (resp. latter) the
pointing in (resp. out) $O_i^*.$  Actually Figure 3.3 includes
3 more relations.   Because the oriented $O_1$ is a $G-$translate
of the oriented $O_2$ and $f$ is a primitive element of $G$ which
fixes $B,$ it follows that 
$${\rm distance}(({\rm pointing\ in}\ O_1^*),
({\rm pointing\ in }\ O_2^*))=0 \pmod L$$  and 
$${\rm distance}(({\rm pointing\ out}\ 
O_1^*), ({\rm pointing\ out}\ O_2^*))=0 \pmod L.$$  
Finally
distance($w(O),O)=R.$  

We therefore obtain the following two
equations

(3.2) $c-R+L+c+\pi i = 0 \pmod L\ \ \ \ \ $
$c-L+R+c+\pi i = 0 \pmod L$

and hence

(3.3) $2R=0 \pmod L\ \ \ \ \ $
$4 c =0 \pmod L.$

\medskip \noindent Using the $L^\prime$ and $R^\prime$ ranges for $X_0$ provided in Proposition 1.28, it is easy to compute that  for each element of ${\cal R} = \exp^{-1}(X_0),\  |{\rm Re}(R)|
< |{\rm Re}(L)/2|.$   It then  follows  that $R=0$ for $(L, D, R) \in {\cal R} \cap {\cal T}.$   \qed

\medskip \noindent {\bf Lemma 3.8:}  If $(L,D,R) \in \cal R \cap {\cal T},$ then $L=D.$

\medskip \noindent {\bf Proof:}  We will now  use the exponential coordinates
$l=\exp(L),\ d=\exp(D).$  These $l,d$ should not be confused with
the $l+it, d+ib$ used above.  In the following calculations
Mathematica [Math] was used to perform matrix
multiplication of $2\times 2$ matrices with coeficients
rational functions in the variables $\sqrt l, \ \sqrt d$.  We will follow Mathematica's notation (for example, in Mathematica's notation the $2 \times 2$ matrix $(a_{ij})$ is
$\{\{a_{11}, a_{12}\}, \{a_{21}, a_{22}\}\}$).
In particular, plugging in $R = 0$ in Lemma 1.24 we get the following matrix representations of $f,F,w,W.$  Because the $R$ term drops out, we can express the matrices of $w$ and $W$ as functions of $d$ alone.

$$f[l\_]=\{\{{\rm Sqrt}[l],0\},\{0,1/{\rm Sqrt}[l]\}\}$$
$$F[l\_]=\{\{1/{\rm Sqrt}[l],0\},\{0,{\rm Sqrt}[l]\}\}$$
$$\displaylines{w[d\_]=\{\{1/2({\rm Sqrt}[d]+1/{\rm Sqrt}[d]),1/2({\rm Sqrt}[d]-1/{\rm Sqrt}[d])\},\cr
\{1/2({\rm Sqrt}[d]-1/{\rm Sqrt}[d]),1/2({\rm Sqrt}[d]+1/{\rm Sqrt}[d])\}\}\cr}$$
$$\displaylines{W[d\_]=\{\{1/2({\rm Sqrt}[d]+1/{\rm Sqrt}[d]),1/2(-{\rm Sqrt}[d]+1/{\rm Sqrt}[d])\},\cr
\{1/2(-{\rm Sqrt}[d]+1/{\rm Sqrt}[d]),1/2({\rm Sqrt}[d]+1/{\rm Sqrt}[d])\}\}\cr}$$

  Let $Y$ be the relator
$F.w.f.w.f.W.f.w.f.w,$ which is a cyclic permutation of
relator i) of Lemma 3.4.   Multiplying
this product of 10 matrices we obtain the
following matrix for $Y$ which we know is $I$.

$$\displaylines{Y[l\_,d\_]=\{\{((1 + d)*(1 - 2*d^2 + d^4 + 8*d*l - 16*d^2*l \cr
                + 8*d^3*l - 2*l^2 + 4*d*l^2 - 4*d^2*l^2 + 4*d^3*l^2 -\cr
        2*d^4*l^2 + l^4 + 4*d*l^4 + 6*d^2*l^4 + 4*d^3*l^4\cr
        + d^4*l^4))/(32*d^{(5/2)}*l^{(3/2)}),\cr}$$

   $$\displaylines{((-1 + d)*(1 + 4*d + 6*d^2 + 4*d^3 + d^4 + 4*d*l +\cr
                8*d^2*l + 4*d^3*l - 2*l^2 - 12*d^2*l^2 - 2*d^4*l^2 +\cr
        4*d*l^3 + 8*d^2*l^3 + 4*d^3*l^3 + l^4 + 4*d*l^4 +\cr
        6*d^2*l^4 + 4*d^3*l^4 + d^4*l^4))/(32*d^{(5/2)}*l^{(3/2)})\},\cr}$$

  $$\displaylines{\{((-1 + d)*(1 + 4*d + 6*d^2 + 4*d^3 + d^4 + 4*d*l +\cr
                8*d^2*l + 4*d^3*l - 2*l^2 - 12*d^2*l^2 - 2*d^4*l^2 +\cr
        4*d*l^3 + 8*d^2*l^3 + 4*d^3*l^3 + l^4 + 4*d*l^4 +\cr
        6*d^2*l^4 + 4*d^3*l^4 + d^4*l^4))/(32*d^{(5/2)}*l^{(3/2)}),\cr}$$

   $$\displaylines{((1 + d)*(1 + 4*d + 6*d^2 + 4*d^3 + d^4 - 2*l^2 + 4*d*l^2 -\cr
                4*d^2*l^2 + 4*d^3*l^2 - 2*d^4*l^2 + 8*d*l^3 -\cr
        16*d^2*l^3 + 8*d^3*l^3 + l^4 - 2*d^2*l^4 + d^4*l^4))\cr
        /(32*d^{(5/2)}*l^{(3/2)})\}\}\cr}$$

Since $G$ is generated by $f,$ an $L$ translation  along $B,$
and $w,$ a $D$ translation along $C,$ it follows from 
Lemma 3.6 that the relation $Y=I$ holds with $l$
and $d$ switched.  Thus $0=Y_{12},\ $ and $0=Y_{12}$ (with $l,d$ switched)
which implies that

$$\displaylines{0=(1 + 4*d + 6*d^2 + 4*d^3 + d^4 + 4*d*l + 8*d^2*l +\cr
                4*d^3*l - 2*l^2 - 12*d^2*l^2 - 2*d^4*l^2 +\cr
        4*d*l^3 + 8*d^2*l^3 + 4*d^3*l^3 + l^4 + 4*d*l^4 +\cr
        6*d^2*l^4 + 4*d^3*l^4 + d^4*l^4)-\cr
        (1 + 4*l + 6*l^2 + 4*l^3 + l^4 + 4*l*d + 8*l^2*d +\cr
        4*l^3*d - 2*d^2 - 12*l^2*d^2 - 2*l^4*d^2 +\cr
        4*l*d^3 + 8*l^2*d^3 + 4*l^3*d^3 + d^4 + 4*l*d^4 +\cr
        6*l^2*d^4 + 4*l^3*d^4 + l^4*d^4)=\cr}$$

                $$4*(1 + d)^2*(1 + l)^2*(-d + l)*(-1 + d*l).$$

This implies that $d=l$ and hence $D=L$ or we  obtain one
of the following solutions which contradicts the
condition  ${\rm Re}(D)>0.$  The solution $d=-1$ implies
$D=\ln(d)=\ln(-1)=\pi i.$ The solution $d=1$ implies $D=0.$  The
solution $d*l=1$ implies that $d=1/l$ and hence
$D=\ln(d)=\ln(1/l)=-\ln(l)=-L$ and hence ${\rm Re}(D)<0.\qed$

\medskip \noindent {\bf Lemma 3.9:}  If $(L,D,R) \in \cal R \cap {\cal T},$ then $d = \exp(D)$ is a root of the polynomial
$$1+2*d+6*d^2+2*d^3+d^4.$$

\medskip \noindent {\bf Proof:}  The equation $Y_{12}=0$ yields

$$\displaylines{0= 1 + 4*d + 6*d^2 + 4*d^3 + d^4 + 4*d*l + 8*d^2*l + 4*d^3*l - 2*l^2 -\cr
12*d^2*l^2 - 2*d^4*l^2 + 4*d*l^3 + 8*d^2*l^3 + 4*d^3*l^3 + l^4 + 4*d*l^4 +\cr
6*d^2*l^4 + 4*d^3*l^4 + d^4*l^4.\cr}$$

Setting $l=d$ we obtain

$$\displaylines{0=\cr
1 + 4*d + 8*d^2 + 12*d^3 - 2*d^4 + 12*d^5 + 8*d^6 + 4*d^7 + d^8=\cr
(1 + 2*d - 2*d^2 + 2*d^3 + d^4)*(1 + 2*d + 6*d^2 + 2*d^3 + d^4).\cr}$$

On the other hand setting $l=d$ in the equation $Y_{11}=1$ we obtain
$$32d^5=(1+d)(1 + 4*d^2 - 12*d^3 + 6*d^4 + 8*d^5 + 4*d^6 + 4*d^7 + d^8)]$$
and so $$0=(-1 + d)*(1 + 2*d + 6*d^2 + 2*d^3 + d^4)*(-1 + 4*d^3 + d^4)$$

The only common solutions to these equations are the roots of the equation
\medskip
\line{(3.4) \ \ \hfil $(1 + 2*d + 6*d^2 + 2*d^3 + d^4)\qed$ \hfil}

\medskip \noindent {\bf Remark 3.10:}  The  two facts $Y=I,$ and $l, d$ can be switched in $Y$ follow from
relator i) of Lemma 3.4 and $R=0.$   Relator ii)  was used in the proof that
$R=0.$

\medskip \noindent {\bf Lemma 3.11:}  The roots of $1+2*d+6*d^2+2*d^3+d^4$ are
$$(-1 + i*\sqrt3)/2 - (-6 - 2*i*\sqrt3)^{(1/2)}/2$$
$$(-1 + i*\sqrt3)/2 + (-6 - 2*i*\sqrt3)^{(1/2)}/2$$
$$(-1 - i*\sqrt3)/2 - (-6 + 2*i*\sqrt3)^{(1/2)}/2$$
$$(-1 - i*\sqrt3)/2 + (-6 + 2*i*\sqrt3)^{(1/2)}/2$$
\qed

\medskip \noindent {\bf Remark 3.12:}  Note that if $x$ is a root of $1+2*d+6*d^2+2*d^3+d^4,\ $then so are
$\bar x,\  1/x\ $ and $1/\bar x.$

\medskip \noindent {\bf Lemma 3.13:}  If $(L,D,R) \in \cal R \cap {\cal T},$ then $$\displaylines{D=L=\ln(
(-1 + i*\sqrt3)/2 - (-6 - 2*i*\sqrt3)^{(1/2)}/2)\cr=\omega
\approx 0.83144294552931 - 1.945530759503636*i.\cr}$$

\medskip \noindent {\bf Proof:}  The other 3 solutions are $-L, Re(L)-Im(L), -Re(L)+Im(L)$ and
lie outside ${\cal R}.$  The above solution lies in ${\cal R}.\qed$

\medskip \noindent {\bf Remark 3.14:} In our language, Hodgson and Weeks [HW1] knew that $\pi_1$(Vol3)
was generated by $f,w$ with $D=L=\omega$  and $R=0.$  Also
that the various solutions of (3.4) corresponded to symmetries of Vol3.

\medskip \noindent {\bf Corollary 3.15:} Vol3 is the unique hyperbolic 3-manifold with
>associated parameter values in  ${\cal S} \cap X_0.$

\medskip \noindent {\bf Proof of Proposition 3.1:}  By Proposition 2.8 if
maxcorona$(\delta)\ge \ln(3)/2$, then the parameter for a 2-generator subgroup $G$ of $\pi_1(N)$ lies in
$S\cap X_0. $  Since $G= \pi_1$(Vol3), $N$ is covered by Vol3.  By Reid (see [JR]) Vol3
only covers Vol3.  Therefore $N=$Vol3. $\qed$

\medskip \noindent {\bf Proof Proposition 3.2:} We  show that
the Dirichlet insulator associated to a shortest geodesic $\delta$ in $N$ satisfies the insulator condition if  ${\cal  C}(O(2)) <  113.16$ degrees.  Here $O(2) = \distance(B, {\rm second\ closest\ lift\ of\ } \delta)$ and  $O(1) = \omega$ corresponds to the closest lift.  This actually requires that we be a bit more precise, so we insert the following definition.

\medskip \noindent {\bf Definition 3.16:} Let $\delta$ be a geodesic in the hyperbolic 3-manifold $N$
and $\{\delta_i\}_{i\ge 0}$ be the preimages of $\delta$ in ${\bf H^3}$.  Define an
equivalence relation on $\{\delta_i\}_{i\ge 1}$ by saying that $\delta_i =\delta_j$ if
either

i) there exists a $g\in\pi_1(N)$ such that $g(\delta_0)=\delta_0$ and
$g(\delta_i)=\delta_j$

\noindent or

ii) there exists $g\in \pi_1(N)$ such that $g(\delta_i)=\delta_0$ and
$g(\delta_0)=\delta_j.$

Call an equivalence class an {\it orthoequivalence class}.   Consider the collection
$\cal O$ of complex numbers 
$\{\distance(\delta_0, \delta_i)\mid i>0\ $ and only one
$\delta_i\ $ is represented in each equivalence class $\}$  Now order $\cal O$ to obtain
the {\it ortholength spectrum} of $\delta, \ \ \{O(1), O(2), \cdots\}$ where $i\le j$ implies
${\rm Re}(O(i))\le {\rm Re}(O(j)).$

The {\it based ortholength spectrum} consists of all distinct pairs of complex
numbers of the form $\{(\distance(C, $ ortholine from $B$ to
$g(B)), \distance(B,g(B))\mid g\in \pi_1(N)\}.$
\bigskip
 
   The geodesic plane midway between $B$ and
$w^{-1}(B)$ intersects the $B-C$ plane in a geodesic $E$ at distance  ${\rm Re}(\omega)/2 = {\rm Re}({\rm length}f)/2$
from $B.$  Conjugate $\pi_1(N)$ so that $f$ is still an $\omega$ translation of $B$ and
 so that one endpoint of $E$ is at (1,0) and the other endpoint is at $(x,0),$
where $x >1.$  By [Beardon; p. 166 ] $x < 23.815.$ Thus the Dirichlet insulator
$\lambda$ (resp. $w(\lambda$) between $B,\ w^{-1}(B)\ $ (resp. $B,\ w(B))\ $ is symmetric about the $x-$axis and lies
within the circle passing through (1,0), (23.815,0) (resp. (-1,0), (-23.815,
0)).  By  [G; Lemma 4.7] this circle takes up a
visual angle of less than 133.68 degrees.  Now $f$ is the composition of an
$\exp({\rm Re}(\omega))$ homothety centered about the origin and an ${\rm Im}(\omega)$ radian $\approx
-111.4707$ degree rotation. Because $\exp(4* {\rm Re}(\omega)) > 27.82 > 23.815$ it follows that the geodesic $E$ is taken ``beyond" $E$ by $w^4.$  In fact, we see
that $\lambda\cap (f^n(\lambda)\cup f^nw(\lambda))=\emptyset$ if $|n| \ge 4.$
Therefore if there was a
trilinking among three insulators associated to the orthoclass of $w(B),$ there would be a tri-linking involoving 3 circles from the collection  $\{f^nw(\lambda),  f^n(\lambda)\mid -3\le
n\le 0\}.$

Since ${\bf S^2_{\infty}}$ is rotated by at most 111.48 degrees
under $f^{\pm1},$ and $\lambda$ takes up less than 133.68 degrees $B$-visual angle it follows that
$f^{-1}(\lambda)\cup\lambda$ take up less than $133.68+111.48=245.06$
degrees $B-$visual angle.  Similar arguments show that if $\lambda$ and one
of $f^n(\lambda)$ or $f^mw(\lambda)$ nontrivially intersect, then the
union cannot take up more $B$-visual angle, in fact except for
$f^\pm(\lambda)$, it takes up less. See Figure 3.4 which shows the union of
$\{f^nw(\lambda),  f^n(\lambda)\mid -3\le n\le 0\}.$   Therefore, if the union of
three such circles was connected,  they would take up at most
$133.68+2(111.48)=356.64$ degrees of $B$-visual angle and hence would not create a
trilinking.

\midinsert
\centerline{\psfig{file=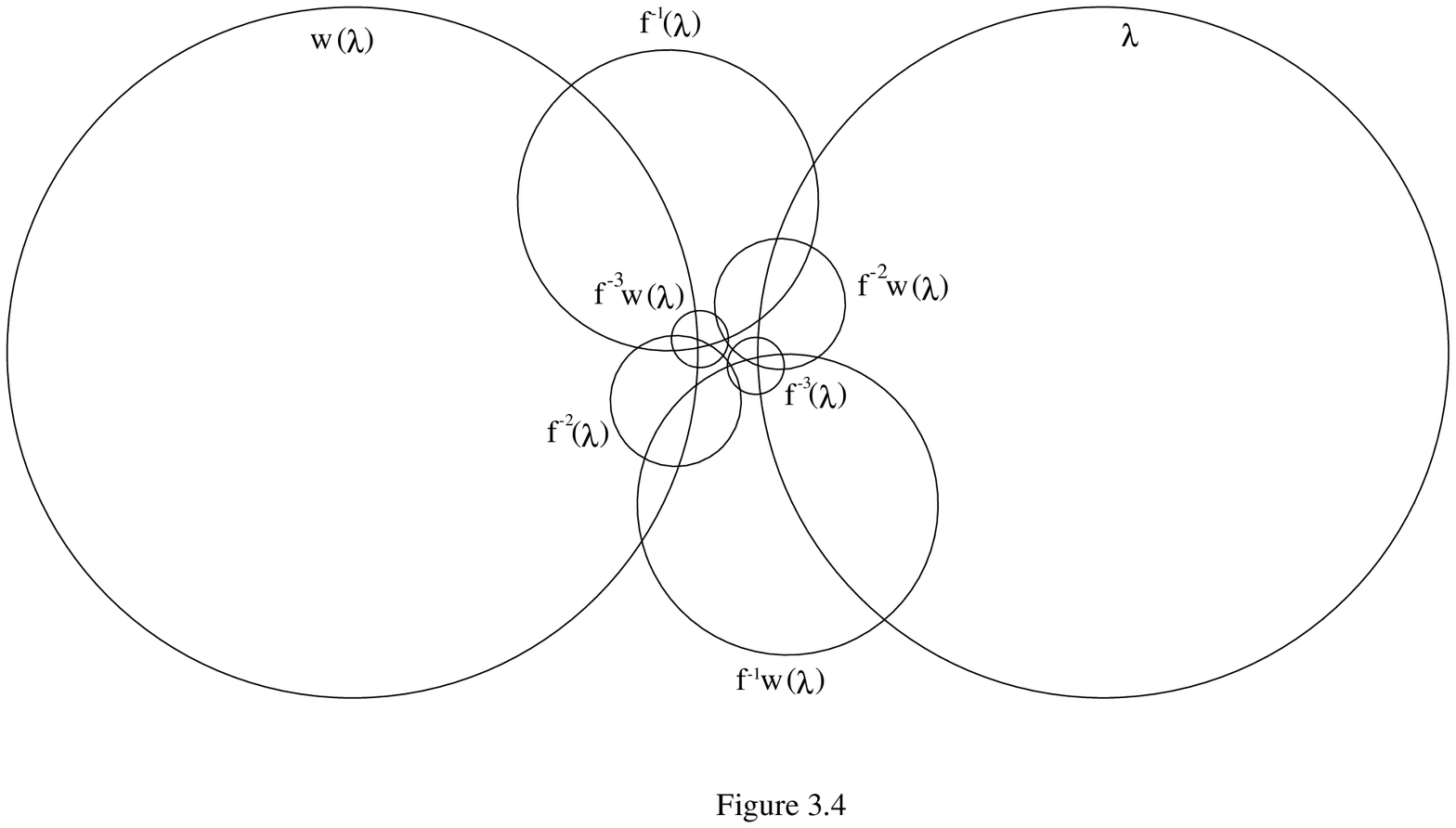,height=6cm}}
\endinsert 

Now suppose that $C(O(2))< 113.16.$  Thus any insulator
associated to a translate of $B$ not in the $O(1)$ class would take less than
113.16 degrees visual angle.  
There are three cases to consider: trilinking involving no $O(1)$ insulators; trilinking involving exactly one $O(1)$ insulator; 
trilinking involving exactly two $O(1)$ insulators.  The case of no $O(1)$ insulators cannot occur because  $3(113.16)<360$.  The case of exactly one $O(1)$ insulator cannot occur because
$B$-visual angle($O(1))+2(C(O(2))<133.68 + 2(113.16)=360.$ 
Finally, the case of exactly two $O(1)$ insulators cannot occur, because the two $O(1)$ insulators would take up at most 245.06 degrees of visual angle (see the above paragraph) and adding to this the fact that $C(O(2)) < 113.16$ produces less than 360 degrees. 

The proof of Proposition 3.2 is now completed by the following

\medskip \noindent {\bf Lemma 3.17:}  If $\delta$ is a shortest geodesic in Vol 3, then $C(O(2))<113.16.$

\medskip \noindent {\bf Proof:}  
{\bf Step 1:}  Re$(O(2))>{\rm Re}(O(1))$

Proof of Step 1.  If Re$(O(2))=$Re$(O(1)),$ then it follows by Proposition 3.1
that $O(2)=O(1)=\omega.$  Furthermore the collection of orthocurves must appear
symmetrically  on the geodesic $B,$ i.e. if distance$(O(2)-{\rm ortholine},\ C)=x+yi,$ with
$0\le x, \ x$ minimal, then there is an ortholine at distance $n(x+yi)$ with $ n\in {\bf Z}$ along $B$
from $C$ whose corresponding ortholength is $\omega.$  Use the fact that if $v\in
\pi_1$(Vol3) is the element with distance$(v(B),B)=O(2)$ as above, then the group
generated by $v,f$ is conjugate to the group generated by $w,f.$

If $x=0,$ and $y>0$ is minimal, then $y=\pi/m.$  Using Fenchel's law of cosines
[F; p. 83 ], it follows that if $m=2,$ then there exists an $\bar \omega$
ortholength thereby contradicting Proposition 3.1. See
Figure 3.5. If $m>2,$ then a similar
argument shows that there exists a real ortholength less than Re$(\omega).$

\midinsert
\centerline{\psfig{file=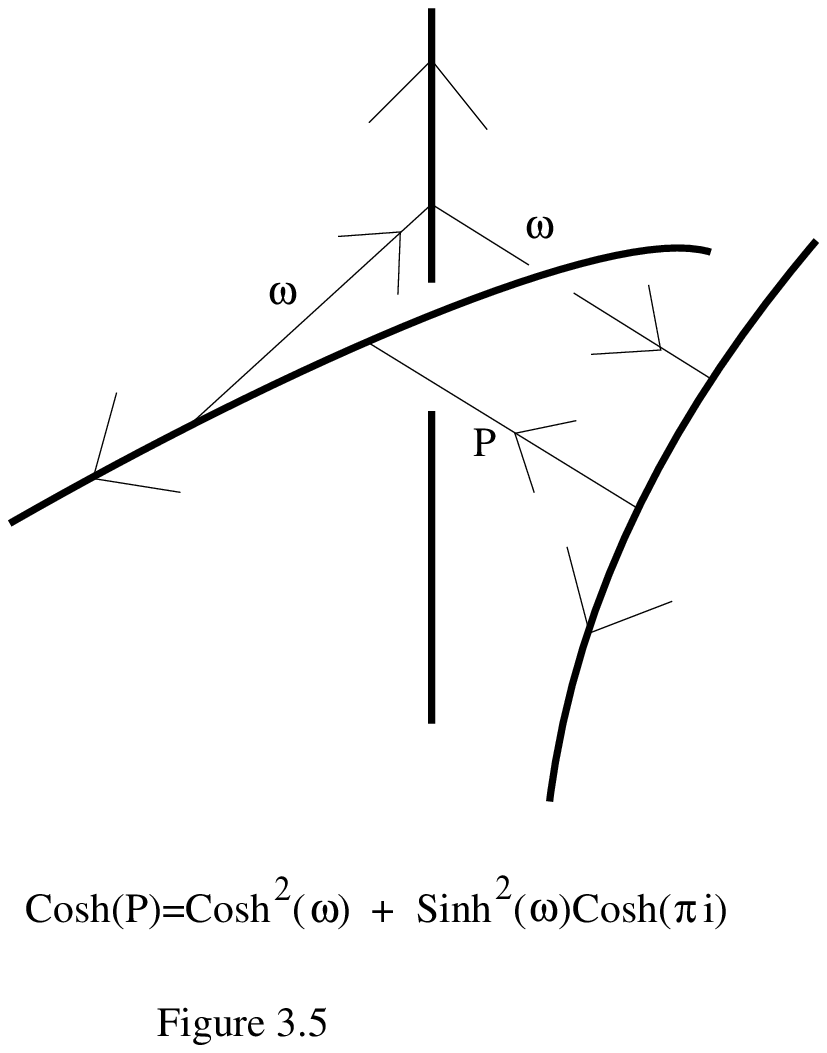,height=6cm}}
\endinsert 

If $x>0$ choose $m$ minimal so that $mx=$Re$(\omega).$  By replacing $y$ by $y+\pi i,$ if
necessary, we can assume that $m(x+yi)=\omega$.  Therefore an $x+iy$
translation $\tau$ along $B$ descends to a ${\bf Z}/m$ action $\phi$ on Vol3. Any
lift of $\phi^n$ is a conjugate of $\tau^n$ which is fixed-point free or the identity.
This  contradicts the fact that Vol3 only covers Vol3.
$\qed$

\bigskip
{\bf Step 2:} $C(O(2))<113.16$
\medskip
Proof of Step 2.  We obtain this result with computer assistance in a manner similar to
the proofs of Propositions 1.28 and 2.8.  Our parameter space ${\cal W}$ is the usual initial box.   As before, a parameter gives rise to a
group $G$ with generators $f, v.$  (Here we reserve the letter $w$ to denote one of
the generators of $\pi_1$(Vol3).)  We consider the parameters $U\subset {\cal W}$ such that
$f, v$ generate a group $G$ where $f$ is of minimal length, and maxcorona$(\delta_f)=
C(v) \ge 113.16.$ If Step 2 was false then $U\neq\emptyset$ for it contains the
2-generator subgroup of $\pi_1$(Vol3) generated by $f, v$ where
distance$(v(B),B)=O(2).$  We now show that $U=\emptyset$ as follows.

We partition ${\cal W}$ into regions ${\cal W}_i$ such that each ${\cal W}_i$ can be eliminated for one of
the following reasons.

z)  ${\cal W}_i = X_0.$

a)  There exists no $\beta\in {\cal W}_i$ such that $\length(f_\beta) \in {\cal L}(X_0).$  In particular, $\length(f_\beta) \neq \omega$ throughout ${\cal W}_i.$

b)       ${\cal W}_i$ has some $L$ values in  ${\cal L}(X_0)$ but
there exists no $\beta\in {\cal W}_i$ such that the real part of 
${\rm distance}(v_\beta(B),B)$ is greater than the minimum $d$ value for $X_0.$
In particular,  ${\rm Re(distance}(v_\beta(B),B)) \le {\rm Re}(\omega)$ throughout ${\cal W}_i.$

c)  There exists a killerword $h$ in $f,w,f^{-1},w^{-1},$ such that
$h_\beta \neq {\rm id}$ and
${\rm Relength}(h_\beta) < {\rm Relength}(f_\beta)$  
for all $\beta\in {\cal W}_i.$

\medskip
 \noindent It turns out that a Corona condition is not needed---c) is sufficient to generate all the killer words we need!

There are files containing the partition of ${\cal W}$ and the associated conditions/killerwords, and {\it fudging} verifies that they
indeed work.  {\it Fudging} actually takes care of the sub-boxes $X_0, \ldots, X_6$ all at once, and the associated files reflect that fact.

As noted earlier, {\it fudging} works with the exponentiated versions of the above conditions.
\qed

\bigskip

Here is an experimental ``proof" of Lemma 3.17.  In [HW2] an algorithm is given to
compute, with multiplicities,  the length spectrum of a hyperbolic 3-manifold $M,$
given a Dirichlet domain for $M.$   Weeks has observed [Weeks, personal
communication] that a very similar argument gives an algorithm to compute the
based ortholength spectrum.  In fact an analogue to Proposition 1.6.2 [HW2],
(with an analogous proof) is the following

\medskip \noindent {\bf Lemma 3.18 (Weeks):} Let $M$ be a closed orientable 3-manifold having a Dirichlet
domain $\cal D$ with basepoint $x$ and with spine radius $r.$  Let $\delta$ be a
geodesic
of length $l+it.$  To compute all the based ortholengths of real length less than or equal to
$\lambda$
with basing less than or equal to $l/2$ from some point on $\delta_0$ (a preimage of $\delta$) it
suffices to find all translates $g\cal D$ satisfying $\rho(x,gx)\le
2r+2{\rm Arccosh}(\cosh(l/2)\cosh(\lambda/2)).$

\medskip \noindent {\bf Proof:}  As in [HW2] we can assume that $\rho(\delta_0,x)\le r.$   The 0-basing on
$\delta_0$ will be given by the oriented perpendicular $P$ from $\delta_0$ to $x.$
Figure 3.6 shows that if there is a translate $\delta_i=g(\delta_0)$ based at
distance $ \le l/2,$ at real distance $\le \lambda$ from $\delta_0,$ then
$\rho(g(x),x)\le \lambda +l+2r$.  As in [HW2] an application of the hyperbolic
cosine law yields the better estimate of the Lemma.    $\qed$

\midinsert
\centerline{\psfig{file=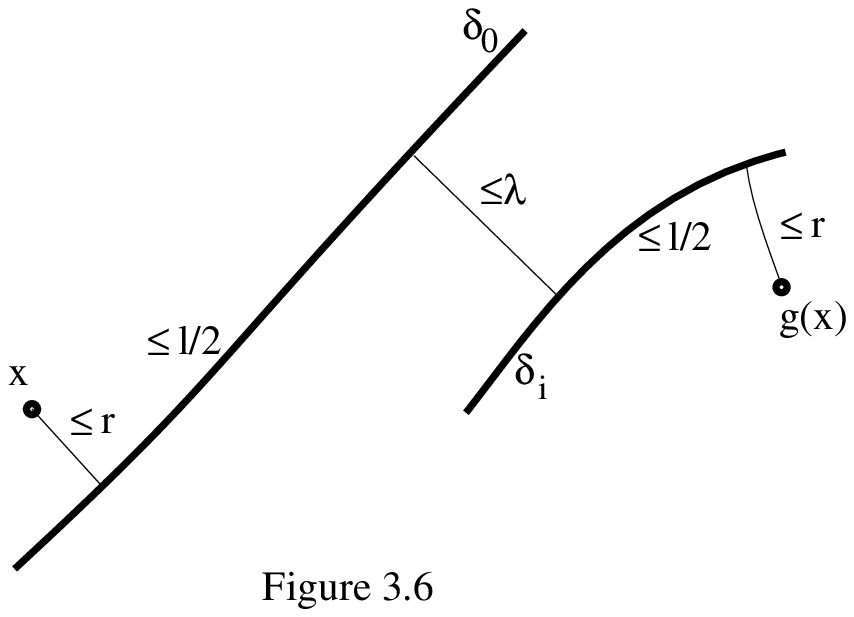,height=4.5cm}}
\endinsert 

Provided one is given a Dirichlet domain, [HW2] gives an efficient algorithm to find these $g$'s.   Finally to each such $g$ one
computes the basing and distance from $g(\delta_0)\  {\rm to}\  \delta_0.$

The collection of ortholengths with basings $\le l/2$ contains at least 2
representatives for each ortholength class,  thus the Weeks algorithm can
be used
to give lower bounds on the various $O(i)$'s.  Snappea computes a Dirichlet domain
for Vol3 with spine radius $\le 0.68.$  
Because $l< 0.83145,$ taking $\lambda=1.24$ we
obtain $2r+2{\rm arccosh}(\cosh(l/2)\cosh(\lambda/2))<2.89.$  This algorithm has been
implemented on an undistributed version of Snappea, and provided the following estimates.  
Note that ${\rm Re}(O(2)) \ge 1.24$ is sufficient to guarantee that ${\cal C}(O(2)) < 113.16$ degrees.
$$O(1)\approx .83144 - 1.94553 I,$$
$$O(2)\approx 1.3170 - \pi I, $$
$$O(3)=O(4) \approx 1.4197 + 1.0963 I, $$ 
$$O(5) \approx 1.9769 -1.2995. $$ 
These
estimates were found using a ``tiling" of radius $3.00>2.89,$ which is
sufficient for ``proving" that for Vol3, ${\rm Re}(O(2)) > 1.24.\qed$

\medskip \noindent {\bf Remark 3.19:} This experimental  proof should be easily rigorizable by
implementing Weeks's algorithm using exact arithmetic.

\sectionskip\magnification=1200
\baselineskip = 0.2 true in
\def\qed{{\bf \ \ \ \ \vrule height6pt width5pt depth1pt}}

\noindent {\bf Chapter 4: Applications}
\bigskip
We provide some applications.  First, we give a partial answer to Question 1.16.  Second a  lower bound for the volume of hyperbolic 3-manifolds is produced.  Finally, we give a relationship between isotopic closed curves in a hyperbolic 3-manifold and essential links in $B^3.$

\medskip \noindent {\bf Theorem 4.1:} If $\delta$ is a shortest geodesic in a closed
orientable hyperbolic 3-manifold $N,$ then either 
tuberadius($\delta) \ge 1.059191579962\ldots/2$
or $N =$Vol3.
        
\medskip \noindent {\bf Proof:}  If tuberadius$(\delta)=d/2 < \ln(3)/2,$ then by Corollary 1.29 ii,
$d\in$
$$ \displaylines {[0.831426508686\ldots, 0.831461989726\ldots] \cup
[1.068029907104\ldots, 1.068134862048\ldots] \cup \hfill \cr
\hfil [1.072078514724\ldots, 1.072622271888\ldots] \cup
[1.094788113123\ldots, 1.095007375882\ldots] \cup \hfil \cr
\hfill [1.094973832380\ldots, 1.095231885426\ldots] \cup
[1.059191579962\ldots, 1.060372139694\ldots] .\cr}$$ 
\noindent If $d<1.059191579962\ldots,$ 
then  by Proposition  1.28 the  parameter
$(L,D,R)\in {\cal P}^\prime$ 
associated to $N$ must be in ${\cal R} \cap {\cal T} = \exp^{-1}(X_0 \cap {\cal S}).$ 
It then 
follows by Corollary 3.15 that $N={\rm Vol3}.$\qed


\medskip \noindent  {\bf Conjecture 4.2:}  If $\delta$ is a shortest geodesic in a closed
orientable hyperbolic 3-manifold $N,$ then either 
tuberadius$(\delta) > \ln(3)/2$
or $N$ is one of six exceptional manifolds.

\medskip \noindent {\bf Remark 4.3:}  Vol3 is one of these manifolds.  As in Remark 1.32,
three others are conjecturally, $s479(-3,1),\ s778(-3,1),\ v2678(2,1).$  These correspond to sub-boxes $X_5, X_2, X_1$ respectively.

The conjectured fundamental
groups of the six manifolds are $<f,w;\ r_1(X_k),r_2(X_k)>$ for $k\in
{0,\ldots,6}.$  (The groups for $k=5,\ 6$ are isomorphic).  Using [Be]
one can get explicit Heegaard genus 2 descriptions of all the conjectured manifolds.

One could prove Conjecture 4.2 by first showing that for each $k,$
${\cal T} \cap X_k$ is a point $T_k.$  Let $M_k$ be the hyperbolic 3-manifold which
corresponds to $T_k.$  Second show that $M_k$ nontrivially covers no
3-manifold.  This procedure was carried out in Chapter 3 for the subbox $X_0.$

\medskip \noindent {\bf Remarks 4.4:}    The previous best lower bound for the volume of hyperbolic 3-manifolds was on the order of 0.001 (see [GM1] and [M2]).  Using the results of this paper and the method of [M1] it is easy to improve this to 0.1.  However, Gehring and Martin provide an improved tube-volume formula in [GM2] and we use their formula to get a lower bound of $0.16668\ldots.$
The Gehring-Martin tube-volume formula for manifolds (as opposed to orbifolds) is 
$${\cal V}(t) = \sqrt 3 \tanh(t) \cosh(2t) {\rm Arcsinh}^2(\sinh(t)/\cosh(2 t))$$
where $t$ is the radius of the embedded solid tube.  Note that the length of the core geodesic is irrelevant. 

\medskip \noindent {\bf Theorem 4.5:} ${5 \over 2 \sqrt 3} {\rm Arcsinh}^2({ \sqrt 3 \over 5}) = 0.16668\ldots$ is a lower bound for the volume of closed hyperbolic 3-manifolds.

\medskip \noindent {\bf Proof:}   [GM2; Corollary 1.7] applies the tube-volume formula ${\cal V}(t)$ to tubes of radius at least $\ln(3)/2$ and produces 
$${\rm Vol}(N) \ge {\cal V}(\ln(3)/2) = {5 \over 2 \sqrt 3} {\rm Arcsinh}^2({ \sqrt 3 \over 5}) = 0.16668\ldots$$

Corollaries 1.29 and 3.15 imply that if the tube
radius of a shortest geodesic in $N$ is less than $\ln(3)/2,$ then either $N={\rm Vol3}$
or $N$ contains an embedded open tube of radius $d/2,$ where $d\in {\rm Re}(\ln({\cal
D}(X_k))),$ about a geodesic of length $l\in {\rm Re}(\ln({\cal L}(X_k)))$ for some $k\in
{1,\ldots,6}.$  As in the Theorem 4.1 we compute that $l\ge 1.059536368901\ldots$

 In the first case, ${\rm Vol}(N)  = {\rm Vol}({\rm Vol3}) = 1.01\ldots,$ while in the second case $l$ and $d$ are bounded as follows: $l \ge 
1.059536368901\ldots$  and $d \ge 1.059191579962\ldots.$
Plugging these into the tube-volume formula $\pi l \sinh^2(d/2)$ we get ${\rm Vol}(N) \ge 1.02419\ldots.$\qed 
 
\medskip \noindent {\bf Theorem 4.6:}  Let $k_1, k_2$  be simple closed curves in $N$ such that $k_1$ is
a geodesic.  Then $k_1$ is isotopic to $k_2$ if and only if as $B^3$-links
$q^{-1}(k_1)$ is equivalent to $q^{-1}(k_2)$ where $q : {\bf H^3} \to N$ is the universal covering
projection.
\medskip \noindent {\bf Proof:}  Apply Corollary 5.6 of [G]. 
Recall that an equivalence between
$B^3$-links $\Gamma,\ \Delta$ is a homeomorphism of $B^3$ which takes $\Gamma$ to
$\Delta$ and fixes $S^2$ pointwise.\qed

\medskip \noindent {\bf Remark 4.7:}  A similar argument extends Theorem 4.6 to homotopy essential
links which lift to trivial $B^3$-links.  The general case of Conjecture 5.5A of [G] is still open.

\sectionskip\magnification=1200
\baselineskip = 0.2 true in
\noindent {\bf Chapter 5: Conditions and Sub-Boxes}
\medskip
As described in Section 1, our goal is to kill off all points in $S \subset {\bf C^3},$   but for computational reasons we will actually work with the box $B \supset S$.  We will decompose $B$ into a collection of sub-boxes such that each sub-box has an associated ``condition" that will describe how to kill off that entire sub-box.  The set-up for describing these sub-boxes will be given in Section 3. 

We now list the conditions used to kill off the sub-boxes.  There are two types of conditions: the trivial and the interesting.  The trivial conditions kill off sub-boxes in $W$ by noting that the sub-box in question misses $S.$  The interesting conditions are where the real work is done, and they require a killer word in $f, w, f^{-1}, w^{-1}$ to work their magic.

To be consistent with the computer program {\it verify} we use the following notation. $L^{\prime} = z_0 + i z_3$, $D^{\prime} = z_1 + i z_4,$ and $R^{\prime} = z_2 + i z_5.$ Here $(L^{\prime}, D^{\prime}, R^{\prime}) \in W$ and $L^{\prime} = \exp(L) = \exp(l+it),\  D^{\prime} = \exp(D) = \exp(d+ib),\  R^{\prime} = \exp(R) = \exp(r+ia).$
\bigskip
\noindent {\bf The Trivial Conditions 5.1:}
\medskip
{\bf Condition `s' (short):}  Tests that all points in the sub-box have $|z_0 + i z_3| < 1.10274.$  This ensures that  $\exp(l) = |\exp(L)| = |L^{\prime}| = |z_0 + i z_3| < 1.10274 < \exp(0.0978),$  and Definition 1.14 tells us that we are outside of $S = \exp(T)$
\medskip
{\bf Condition `l' (long):} Tests that all points in the sub-box have $|z_0 + i z_3| > 3.63201.$  This ensures that  $\exp(l) = |\exp(L)| = |L^{\prime}| = |z_0 + i z_3| > 3.63201 > \exp(1.2897845)$ and we are outside of $S.$ Here we have improved $S$ to allow for the value computed in the proof of Proposition 1.13, rather than the cruder value 1.29.
\medskip
{\bf Condition `n' (near):} Tests that all points in the sub-box have $|z_1 + i z_4| < 1.$  This ensures that  $\exp(d) = |\exp(D)| = |D^{\prime}| = |z_1 + i z_4| <1= \exp(0)$ and we are outside of $S.$  Actually, we used a stronger condition in Definition 1.14 ($l/2 \le d$).
\medskip
{\bf Condition `f' (far):} Tests that all points in the sub-box have $|z_1 + i z_4| > 3.$  This ensures that $\exp(d) = |\exp(D)| = |D^{\prime}| = |z_1 + i z_4|  > 3= \exp(\ln 3)$ and we are outside of $S.$  
\medskip
{\bf Condition `w' (whirle big):} Tests that all points in the sub-box have $|z_2 + i z_5|^2 > |z_0 + i z_3| $  This ensures that $\exp(r) = |\exp(R)| = |R^{\prime}| = |z_2 + i z_5|  > \sqrt{|z_0 + i z_3|} =  \sqrt {\exp(l)} = \exp(l/2)$ and we are outside of $S.$  In the computer program, this test requires a round-off error analysis.  As such, the formula in the program is more complicated than expected.  See Chapters 7 and 8 for a discussion of round-off error.
\medskip
{\bf Condition `W' (whirle small):} Tests that all points in the sub-box have $|z_2 + i z_5| < 1. $  This ensures that $\exp(r) = |\exp(R)| = |R^{\prime}| = |z_2 + i z_5|  < 1 = \exp(0)$ and we are outside of $S.$  \bigskip
\noindent {\bf The Interesting Conditions 5.2:}
\medskip
{\bf Condition `L':} This condition comes equipped with a killer word $k$ in $f$ and $w,$ and tests that all points in the sub-box have $|\exp({\rm length}(k))| <  |L^{\prime}| = |\exp(L)|,$ where ${\rm length}(k)$ means the length of the isometry determined by $k.$   This, of course, contradicts the fact that $L$ is the length of the shortest geodesic.

It is easy to carry out the test $|\exp({\rm length}(k))| <  |L^{\prime}|$ because
Lemma 1.25a) can be used.
Note that in {\it verify} the function which computes $\exp({\rm length})$ is called {\it length}.  

Of course, Condition `L' also checks that the isometry corresponding to the word $k$ is not the identity.
\medskip
{\bf Condition `O':}  This condition comes equipped with a killerword $k$ in $f$ and $w,$ and tests that all points in the sub-box have $|\exp({\rm distance}(k(B), B))| <  |D^{\prime}| = |\exp(D)|.$    This, of course, contradicts the ``nearest" condition.

It is easy to carry out the test $|\exp({\rm distance}(k(B), B))| <  |D^{\prime}|$ because
Lemma 1.25b) can be used. Note that in {\it verify} the function which computes $\exp({\rm distance}(k(B), B))$is called {\it orthodist}.  

Also, Condition `O' checks that the isometry corresponding to the word $k$ does not take the axis of $f$ to itself.
 bigskip
{\bf Condition `2':} This is just the `L' condition without the ``not-the-identity" check, but with the additional proviso that $k$ is of the form $f^p w^q.$ This ensures that $k$ is not the identity, because for $k$ to be the identity $f$ and $w$ would have to have the same axis, which contradicts the fact that $d$ can be taken to be greater than or equal to $l/4.$
\bigskip
{\bf Condition `conjugate':} There is one other condition that is used to eliminate points in $W.$  Following Definition 1.14 (and Lemma 1.15) we eliminate all boxes with $0 < t \le \pi.$  Of course, after exponentiating $L = l+it,$ this corresponds to eliminating all boxes with $z_3 > 0.$
Specifically, we toss all sub-boxes of $W$ whose fourth entry is a 1.
This condition does not appear in {\it verify} and {\it fudging} because it is applied ``outside" of these programs, as described 
at the end of this Chapter.
\bigskip
\noindent {\bf Construction 5.3:}
We now give the method for describing the roughly 930 million sub-boxes that the initial box $W$ is sub-divided into.

All sub-boxes are gotten by sub-division of a previous sub-box along a real
hyper-plane mid-way between parallel faces of the sub-box before
sub-division.  Of course, these midway planes are of the form $x_i = $ a
constant.   We use 0's and 1's to describe which half of a sub-divided
sub-box to take (0 corresponds to lesser $x_i$ values).  For example, 0
describes 
the box $W \cap \{(x_0,x_1,x_2,x_3,x_4,x_5): x_0 \le 0 \}$, 
010 describes the box 
$W \cap \{(x_0,x_1,x_2,x_3,x_4,x_5) : x_0 \le 0,\ x_1 \ge 0,\ x_2 \le 0 \}$,
and so on.

In this way, we get a 1-to-1 correspondence
between strings and sub-boxes.
If $s$ is a string of 0's and 1's, then let $Z(s)$ denote
the box corresponding to $s$.
The range of values for the $i$-th coordinate in the sub-box $Z(s)$ is related
to the binary fraction $0.s_{i}s_{i+6}\ldots s_{i+6k}$.
The two sub-boxes gotten from
subdividing $Z(s)$ are $Z(s0)$ and $Z(s1)$.

The directions of sub-division cycle among the various coordinate axes:
the $n$-th sub-division is across the $(n \bmod 6)$-th axis.
The dimensions of the top-level box $W$ were chosen so that sub-division
is always done across the longest dimension of the box,
and so that all of the sub-boxes are similar.  This explains the factor of $2^{(5-i)/6}$ in Definition 1.22.

To kill a sub-box $Z(s)$, the checker program has two (recursive) options:
use a condition (and, if necessary, an associated killer word)  to kill $Z(s)$ directly, or
first kill $Z(s0)$ and then kill $Z(s1)$.

There is also a third option: don't kill $Z(s)$, and instead
mark $s$ as omitted.
Any omitted sub-boxes are checked with another instance
of the checker program,
unless the sub-box is one of the 7 exceptional
sub-boxes (11 before joining abutters).

Thus, a typical output from {\it verify} would be 
$${\rm verified} 000000111101111111 - \{ 0000001111011111110 000000111101111111110 \}$$
Which means that the sub-box $Z(000000111101111111)$ was killed except for its sub-boxes $Z(0000001111011111110)$  and $Z(000000111101111111110).$  The ouput 
$${\rm verified} 0000001111011111110 - \{\ \}.$$ and 
$${\rm verified} 000000111101111111110 - \{\ \}.$$
\noindent shows that these boxes were subsequently killed as well, and thus the entire sub-box $Z(000000111101111111)$ has been killed.

Instead of immediately working on killing the top-level box, we subdivide in the six co-ordinate directions to get the 64 sub-boxes $$Z(000000),Z(000001),Z(000010),Z(000011),\ldots, Z(111111).$$ 
\noindent We then throw out
the ones with fourth co-ordinate equal to 1 (thanks to Definition 1.14 and Lemma 1.15), leaving us with the 32 sub-boxes $$Z(000000),Z(000001),Z(000010),Z(000011),\ldots, Z(111011).$$  \noindent We then use {\it verify} to kill these, with the exception of the 7 exceptional boxes (11 before joining abutters together) listed in Proposition 1.28.

The choices in {\it verify} are made for it by a sequence
of integers given as input.  The sequence of integers containing the directions for killing $Z(000000)$ is contained in the file data6/000000.
In such a sequence, 
$0$ tells {\it verify} to sub-divide the present box (by $x_i = c$),  to position itself on the "left-hand" box ($x_i \le c$) created by that sub-division, and to read in the next integer in the sequence.
A  positive integer $n$ tells {\it verify} to kill the sub-box it
is positioned at directly, using the condition (and killer word, if necessary) on line $n$ in the 'conditionlist' file,  and to then position itself at the ``next" natural sub-box. 
$-1$ tells {\it verify} to omit the sub-box, and mark it as skipped
(the sequence of integers used in killing the skipped box $Z(s)$ is contained in a file data6/s)..

The checker program {\it verify}, its inputs, and the list of conditions
will be available from the Geometry Center.
Details about how to get them can be found at
\bigskip
{\noindent\tt http://www.geom.umn.edu:/locate/HomotopyHyperbolic}
\bigskip
Similarly for the program {\it fudging} which is used on the 7 exceptional boxes.
\bigskip
\noindent {\bf Example 5.4:} 
To illustrate the checking in action, this is a (non-representative)
example, which shows how the sub-box $Z(s)$ (minus a hole) is killed,
where
$$s=001000110001110111001111000101111111101111100111001111000001111011110111.$$

The input associated with this sub-box is
$$(0, 0, 0, 1929, 12304, 0, 0, 7, 0, 1965, 0, 1929, 1929, 1996, -1),$$
which causes the program to kill $Z(s)$ in the following fashion:
{\noindent\obeylines
\def\>{\hskip 0.1in}
kill $Z(s)$:
\>kill $Z(s0)$:
\>\>kill $Z(s00)$:
\>\>\>kill $Z(s000)$ with condition 1929 = ``L(FwFWFWfWFWFwFwfww)''
\>\>\>kill $Z(s001)$ with condition 12304 = ``L(FwfWFFWFwFwfwfWfwfw)''
\>\>kill $Z(s01)$:
\>\>\>kill $Z(s010)$:
\>\>\>\>kill $Z(s0100)$ with condition 7 = ``L(w)''
\>\>\>\>kill $Z(s0101)$:
\>\>\>\>\>kill $Z(s01010)$ with condition 1965 = ``L(fwFwFWFFWFwFwfwww)''
\>\>\>\>\>kill $Z(s01011)$:
\>\>\>\>\>\>kill $Z(s010110)$ with condition 1929
\>\>\>\>\>\>kill $Z(s010110)$ with condition 1929
\>\>\>kill $Z(s011)$ with condition 1996 = ``L(FwFwFWFWfWFWFwFwfww)''
\>omit $Z(s1)$
as shown in figure 5.1.
}

\pageinsert
\centerline{\psfig{file=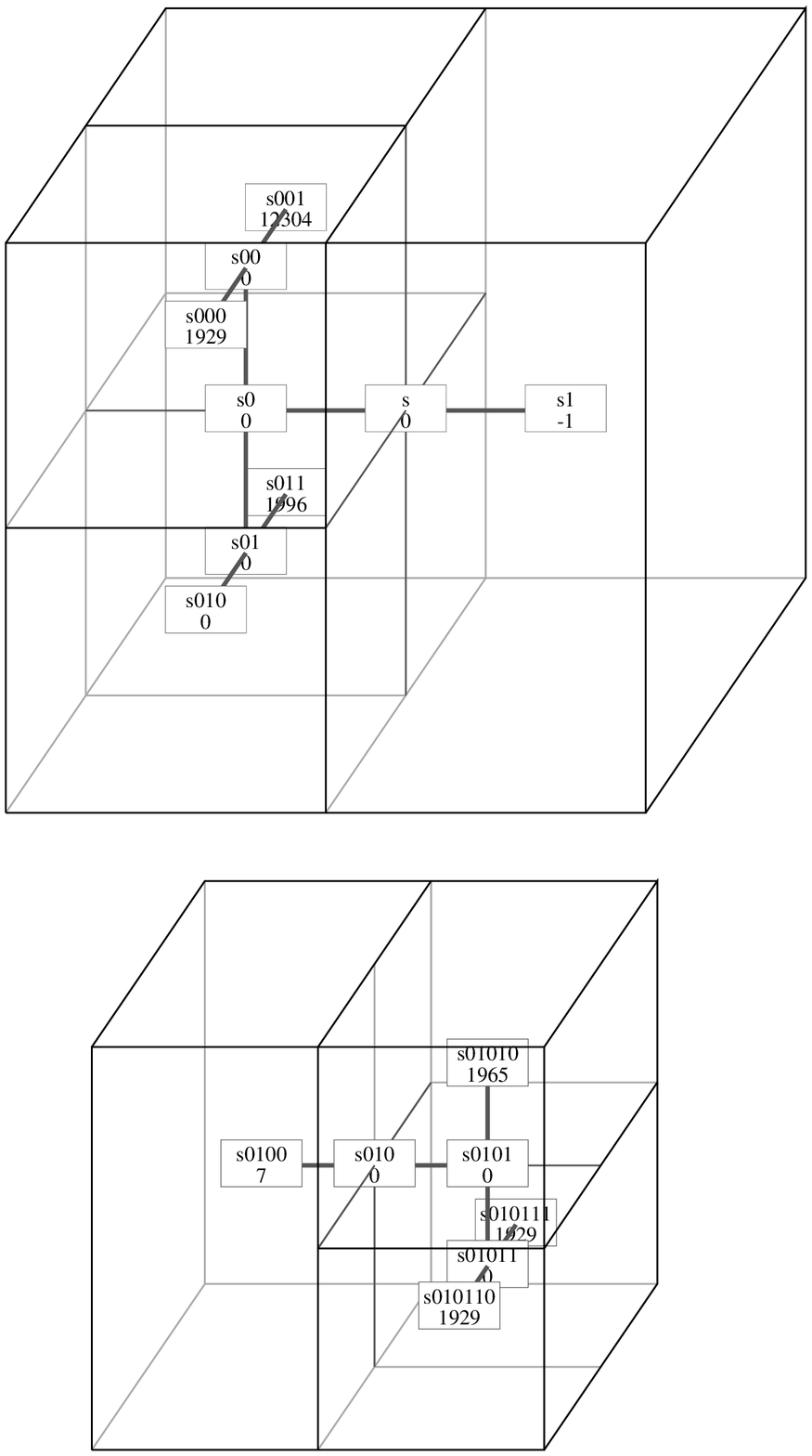,height=15cm}}
\centerline{Figure 5.1: six levels of subdivision, in two projections,
with all the trimmings}
\endinsert

$Z(s1)$ is ignored, so the checker would indicate this omission
in its report.  In fact, $Z(s1)$ is omitted entirely, since it is
one of the 11 exceptional boxes.

The use of condition ``L(w)'' so deep in the tree is unusual.
In this case, it's because the manifold in the exceptional sub-box
has $length(f) = length(w)$, so that the program will frequently
come to places where it can bound $length(f) > length(w)$ nearby.

One might wonder why the checker subdivides $Z(s01011)$,
since it's going to use the same condition to kill both halves.
The reason is the error bound for $Z(s01011)$ wasn't good enough
to prove that the sub-box is killed directly.

The binary numbers used by the computer require too much space to print.
In the example calculation which follows, we instead use a decimal 
representation.  Also, only 10 decimals are printed, less accurate
than the 53 binary digits used for the actual calculations.

The sub-box $Z(s01011)$ is the region where
$$\left(\matrix{
-1.381589027741\ldots  \le  Re(L')  \le  -1.379848991182\ldots\cr
-1.378124546093\ldots  \le  Re(D')  \le  -1.376574349753\ldots\cr
0.999893182771\ldots  \le Re(R')  \le  1.001274250703\ldots\cr
-2.535837191243\ldots  \le  Im(L')  \le  -2.534606799593\ldots\cr
2.535404997792\ldots  \le  Im(D')  \le  -2.534308843448\ldots\cr
-0.001953125000\ldots  \le  Im(R')  \le  0.000000000000\ldots\cr
}\right)$$

\bigskip
At this point, we would like to compute $$f,\ w,\ g = f^{-1}wf^{-1}w^{-1}f^{-1}w^{-1}fw^{-1}f^{-1}w^{-1}f^{-1}wf^{-1}wfww, \ {\rm length}(g),$$ 
\noindent 
and so on.  However, these items take on values over an entire sub-box and thus are computed via AffApprox's (first-order Taylor Approximations with remainder bounds), which are not formally defined until the next Chapter.  As such, we complete Example 5.4 at the end of Chapter 6.

\sectionskip\magnification=1200
\baselineskip = 0.2 true in
\def\qed{{\bf \ \ \ \ \vrule height6pt width5pt depth1pt}}

\noindent {\bf Chapter 6: Affine Approximations}
\medskip \noindent {\bf Remark 6.1:}
To show that a sub-box of the parameter box ${\cal W}$ is killed by one of the interesting conditions (plus associated killerword) we need to show that at each point in the sub-box, the killerword evaluated at that point satisfies the given condition.  That is, we are simply analyzing a certain function from the sub-box to ${\bf C}$.  

As described in Remark 6.5, this analysis can be pulled back from the sub-box in question to the unit complex 3-disc $A$, where
$A = \{(z_0,z_1,z_2) \in {\bf C}^3 : |z_k| \le 1 \ {\rm for} \ k = 1,2,3\}.$  
Loosely, we will analyze such a function on $A$ by using  Taylor series approximations consisting of an affine approximating function together with a bound on the ``error" in the approximation (this could also be described as a ``remainder bound").

\medskip \noindent {\bf Problems 6.2:}
There are two immediate problems likely to arise from this Taylor approximation approach.
The first problem is the appearance of unpleasant functions such as 
$\cosh^{-1}$.  We have already taken care of this problem by ``exponentiating'' our preliminary parameter space ${\cal P}$.  This resulted in all functions under consideration being built up from the co-ordinate functions $ L^{\prime}, D^{\prime},$ and $ R^{\prime}$ on ${\cal P}$ by means of the elementary operations $+,\ -,\ \times,\ /,\ \sqrt.$ 

Second,  for a given ``built-up function'' the computer needs to be able to compute the Taylor approximation, and the error term.  This will be handled in this section by developing combination formulas for elementary operations (see the Propositions below).  Specifically, given two Taylor approximations with error terms representing functions $g$ and $h$ and an elementary operation on $g$ and $h$, we will show how to get the Taylor approximation with error term for the resultant function from the two original Taylor approximations. 
\medskip \noindent {\bf Remark 6.3:}
We set up the Taylor approximation approach rigorously as follows.  The notation will be a bit unusual, but we are motivated by a desire to stay close to the notation used in the checker computer programs, {\it verify} and {\it fudging}.  However, it should be pointed out that the formulas in this Chapter will be superceded by the ones in Chapter 8, which incorporate a round-off error analysis.  It is the Chapter 8 formulas that are used in {\it verify} and {\it fudging}.
\medskip 
\noindent {\bf Definition 6.4:}  An {\it AffApprox}  $x$ is a five-tuple
$(x.f;\ x.f_0,\ x.f_1,\ x.f_2;\ x.e),$
consisting of four complex numbers $x.f,\ x.f_0,\ x.f_1,\ x.f_2$ and one real number $x.e,$  which represents all functions 
$g: A \rightarrow {\bf C}$ such that
$$|g(z_0,z_1,z_2) - (x.f + x.f_0 z_0 + x.f_1 z_1 + x.f _2 z_2)| < x.e$$
for all $(z_0,z_1,z_2) \in A.$  That is, $x$ represents all functions from $A$ to ${\bf C}$ that are $x.e-$well-approximated by the affine function $x.f + x.f_0 z_0 + x.f_1 z_1 + x.f _2 z_2$.  We will denote this set of functions associated with $x$ by $r(x)$.
\medskip \noindent {\bf Remark 6.5:}
As mentioned in Remark 6.1, given a sub-box to analyze, instead of working with functions defined on the sub-box, we will work with corresponding functions defined on $A.$  Specifically, rather than build up a function by elementary operations performed on the co-ordinate functions $L^{\prime}, D^{\prime}                                                                                               , R^{\prime}$ restricted to the given sub-box, we will perform the elementary operations on the following functions defined on $A$, 
$$(p_0 + i p_3; s_0 + i s_3,0,0; 0)\ \ \ (p_1 + i p_4; 0, s_1 + i s_4,0;0) \ \ \ (p_2 + i p_5; 0,0, s_2 + i s_5; 0)$$
where $(p_0 + i p_3, p_1 + i p_4, p_2 + i p_5)$ is the center of the sub-box in question, and the $s_i$ describe the six dimensions of the box. In the computer programs, these three functions are called {\it along}, {\it ortho}, and {\it whirle,} respectively, and $p_i$ and $s_i$ are denoted {\it pos}[i] and {\it size}[i], respectively.
\medskip
After the following Remarks, we state and prove the combination formulas. 
\medskip \noindent {\bf Remarks 6.6:} 
i)  We will break with the convention used previously in this paper and start the numbering of the Propositions with 6.1.  However, we will end this Chapter with Example 6.7.

ii)  The negation of a set of functions is the set consisting of the negatives of the original functions, and similarly for other operations.
\bigskip
\noindent {\bf Proposition 6.1 (unary minus):} If $x$ is the AffApprox $(x.f;\ x.f_0,\ x.f_1,\ x.f_2;\ x.e)$ then $-(r(x)) = r(-x)$ where 
$$-x = (-x.f;\ -x.f_0,\ -x.f_1,\ -x.f_2;\ x.e).$$ \medskip
\noindent {\bf Proof:}   
$$|g(z_0,z_1,z_2) - (x.f + x.f_0 z_0 + x.f_1 z_1 + x.f_2 z_2)| < e$$
if and only if 
$$|-g(z_0,z_1,z_2) - (-x.f - x.f_0 z_0 - x.f_1 z_1 - x.f_2 z_2)| < e$$ \qed
\bigskip
\noindent {\bf Proposition 6.2 (addition):}  If $x$ is the AffApprox $(x.f;\ x.f_0,\ x.f_1,\ x.f_2;\ x.e)$ and $y$ is the AffApprox $(y.f;\ y.f_0,\ yf_1,\ y.f_2;\ y.e)$, then $r(x + y) \supseteq r(x) + r(y)$, where
$$x + y = (x.f + y.f; x.f_0 + y.f_0, x.f_1 + y.f_1, x.f_2 + y.f_2; x.e + y.e)$$
\medskip
\noindent {\bf Proof:}  If $g\in r(x)$ and $h \in r(y)$ then we must show that $g + h \in r(x + y).$
$$|(g + h)(z_0,z_1,z_2) - ((x.f + y.f) + (x.f_0 + y.f_0) z_0 + (x.f_1 + y.f_1) z_1 + (x.f_2 + y.f_2) z_2)| $$
$$\displaylines{\quad
\le |g(z_0,z_1,z_2) - 
(x.f + (x.f_0) z_0 + (x.f_1) z_1 + (x.f_2) z_2) + 
\hfill\cr
\hfill{}  
h(z_0,z_1,z_2) - (y.f + (y.f_0) z_0 + (y.f_1) z_1 + (y.f_2) z_2)| 
                                               \quad\cr}$$
$$\displaylines{\quad
\le |g(z_0,z_1,z_2) - 
(x.f + (x.f_0) z_0 + (x.f_1) z_1 + (x.f_2) z_2)| + 
\hfill\cr
\hfill{}  
|h(z_0,z_1,z_2) - (y.f + (y.f_0) z_0 + (y.f_1) z_1 + (y.f_2) z_2)| 
                                              \quad\cr}$$
$$\le x.e + y.e$$ \qed
\bigskip
We now do subtraction.  The statement and proof are essentially the same as for addition.   The only thing to note is that the errors add.
\medskip
\noindent {\bf Proposition 6.3 (subtraction):}  If $x$ is the AffApprox $(x.f;\ x.f_0,\ x.f_1,\ x.f_2;\ x.e)$ and $y$ is the AffApprox $(y.f;\ y.f_0,\ yf_1,\ y.f_2;\ y.e)$, then $r(x - y) \supseteq r(x) - r(y)$, where
$$x - y = (x.f - y.f; x.f_0 - y.f_0, x.f_1 - y.f_1, x.f_2 - y.f_2; x.e + y.e)$$ \qed
\bigskip
We now state variations on Propositions  6.2 and 6.3 whose usefulness will not be apparent until Chapter 8, when we incorporate round-off error into these formulas.  In what follows, a ``double" corresponds to a real number, and has an associated AffApprox, with last four entries zero.
\medskip
\noindent {\bf Proposition 6.4 (addition of an AffApprox and a double):} If $x$ is the AffApprox $(x.f;\ x.f_0,\ x.f_1,\ x.f_2;\ x.e)$ and $y$ is a double, then $r(x + y) \supseteq r(x) + r(y)$, where
$$x + y = (x.f + y; x.f_0, x.f_1 , x.f_2 ; x.e).$$ \qed
\bigskip
\noindent {\bf Proposition 6.5 (subtraction of a double from an AffApprox):}  If $x$ is the AffApprox $(x.f;\ x.f_0,\ x.f_1,\ x.f_2;\ x.e)$ and $y$ is a double, then $r(x - y) \supseteq r(x) - r(y)$, where
$$x - y = (x.f - y; x.f_0, x.f_1 , x.f_2 ; x.e).$$ \qed
\bigskip
\noindent {\bf Proposition 6.6 (multiplication):}  If $x$ is the AffApprox $(x.f;\ x.f_0,\ x.f_1,\ x.f_2;\ x.e)$ and $y$ is the AffApprox $(y.f;\ y.f_0,\ yf_1,\ y.f_2;\ y.e)$, then $r(x \times y) \supseteq r(x) \times r(y)$, where
$$\displaylines{\quad
x \times y = (x.f \times y.f; x.f \times y.f_0 + x.f_0 \times y.f, 
x.f \times y.f_1 + x.f_1 \times y.f, x.f \times y.f_2 + x.f_2 \times y.f; 
\hfill\cr
\hfill{}  
(size(x) + x.e) \times (size(y) + y.e) + (|x.f| \times y.e + x.e \times |y.f|))
                                              \quad\cr}$$
 with $size(x) = |x.f_0| + |x.f_1| + |x.f_2|$ and $size(y) = |y.f_0| + |y.f_1| + |y.f_2|$
\medskip
\noindent {\bf Proof:}  If $g \in r(x)$ and $h \in r(y)$ then we must show that $g\times h \in r(x \times y).$  That is, we need to show
$$\displaylines{\quad
|(g \times h)(z_0,z_1,z_2) - 
((x.f \times y.f) + 
\hfill\cr
\hfill{}  
(x.f \times y.f_0 + x.f_0 \times y.f) z_0 + (x.f \times y.f_1 + x.f_1 \times y.f) z_1 + (x.f \times y.f_2 + x.f_2 \times y.f) z_2)| 
\hfill\cr
\hfill{}  
\le (size(x) + x.e) \times (size(y) + y.e) + (|x.f| \times y.e + x.e \times |y.f|)
                                              \quad\cr}$$
Note that for any point $(z_0, z_1, z_2) \in A$ and any functions $g \in r(x)$ and $h \in r(y)$ we can find complex numbers $u, v$ with $|u| \le 1$ and $|v| \le 1$, such that
$$ g(z_0, z_1, z_2) = x.f + (x.f_0 z_0 + x.f_1 z_1 + x.f_2 z_2) + (x.e) u$$ and 
$$ h(z_0, z_1, z_2) = y.f + (y.f_0 z_0 + y.f_1 z_1 + y.f_2 z_2) + (y.e) v.$$
Multiplying out, we see that 
$$\displaylines{\quad
(g \times h) (z_0, z_1, z_2) = 
(x.f \times y.f) + 
\hfill\cr
\hfill{}  
(x.f \times y.f_0 + x.f_0 \times y.f) z_0 + (x.f \times y.f_1 + x.f_1 \times y.f) z_1 +
(x.f \times y.f_2 + x.f_2 \times y.f) z_2 + 
\hfill\cr
\hfill{}  
(x.f \times y.e) v + (x.e \times y.f) u + 
\hfill\cr
\hfill{}  
((x.f_0 z_0 + x.f_1 z_1 + x.f_2 z_2) + (x.e) u) \times
((y.f_0 z_0 + y.f_1 z_1 + y.f_2 z_2) + (y.e) v)
                                              \quad\cr}$$
Hence,
$$\displaylines{\quad
|(g \times h) (z_0, z_1, z_2) - ((x.f \times y.f) + 
\hfill\cr
\hfill{}  
((x.f \times y.f_0 + x.f_0 \times y.f) z_0 + (x.f \times y.f_1 + x.f_1 \times y.f) z_1 +
(x.f \times y.f_2 + x.f_2 \times y.f) z_2))|
                                              \quad\cr}$$
$$\le (|x.f| y.e + x.e |y.f|) + 
(size(x) + x.e) \times (size(y) + y.e).$$ \qed
\bigskip
\noindent {\bf Proposition 6.7 (an AffApprox multiplied by a double):}
If $x$ is the AffApprox $(x.f;\ x.f_0,\ x.f_1,\ x.f_2;\ x.e)$ and $y$ is a double, then $r(x \times y) \supseteq r(x) \times r(y)$, where
$$
x \times y = (x.f \times y;  x.f_0 \times y, 
x.f_1 \times y,  x.f_2 \times y; 
x.e \times |y|) $$ \qed
\bigskip
\noindent {\bf Proposition 6.8 (division):} If $x$ is the AffApprox $(x.f;\ x.f_0,\ x.f_1,\ x.f_2;\ x.e)$ and $y$ is the AffApprox $(y.f;\ y.f_0,\ yf_1,\ y.f_2;\ y.e)$, then $r(x / y) \supseteq r(x) / r(y)$, where
$$\displaylines{\quad
x / y = (x.f / y.f; (-x.f \times y.f_0 + x.f_0 \times y.f)/(y.f^2), 
(-x.f \times y.f_1 + x.f_1 \times y.f)/(y.f^2), 
\hfill\cr
\hfill{}  
(-x.f \times y.f_2 + x.f_2 \times y.f)/(y.f^2); 
\hfill\cr
\hfill{}  
(|x.f| + size(x) + x.e) /(|y.f| - (size(y) + y.e)) - 
\hfill\cr
\hfill{}  
((|x.f|/|y.f| + size(x)/|y.f|) + |x.f| size(y)/(|y.f||y.f|)))
                                              \quad\cr}$$
Of course, we have to avoid division by zero.  That is, we demand that $|y.f| > size(y) + y.e.$
\medskip
\noindent {\bf Proof:} 
For notational convenience, denote $(x.f_0 z_0 + x.f_1 z_1 + x.f_2 z_2)$ by $x.f_k z_k$ and similarly for $y.f_k z_k$ and so on.
As above, note that for any point $(z_0, z_1, z_2) \in A$ and any functions $g \in r(x)$ and $h \in r(y)$ we can find complex numbers $u, v$ with $|u| \le 1$ and $|v| \le 1$, such that
$$ g(z_0, z_1, z_2) = x.f + (x.f_k z_k) + (x.e) u$$ and 
$$ h(z_0, z_1, z_2) = y.f + (y.f_k z_k) + (y.e) v.$$

We compare $(g/h)(z_0, z_1, z_2)$ with its putative affine approximation.  That is, we analyze
$$\displaylines{\quad
|(x.f + (x.f_k z_k) + (x.e) u)/(y.f + (y.f_k z_k) + (y.e) v) - 
\hfill\cr
\hfill{}  
((x.f / y.f) + {(x.f_k) y.f - x.f (y.f_k) \over y.f^2} z_k)|
                                             \quad\cr}$$
Putting this over a common denominator  of 
$|(y.f^2)(y.f + (y.f_k z_k) + (y.e) v)|$
and cancelling equal terms (in the numerator) we are left with a quotient whose numerator is 
$$\displaylines{\quad
|x.e (y.f^2) u - (x.f_k) y.f (y.f_k) z_k - x.f (y.f_k^2) z_k + 
\hfill\cr
\hfill{}  
(x.f) y.f (y.e) v + x.f_k (y.f) y.e (v) z_k - x.f (y.f_k) y.e (v) z_k|.                                            					\quad\cr}$$
We must show this (first) quotient is bounded by
$$\displaylines{\quad
(|x.f| + size(x) + x.e) /(|y.f| - (size(y) + y.e)) - 
\hfill\cr
\hfill{}  
((|x.f|/|y.f| + size(x)/|y.f|) + |x.f| size(y)/(|y.f||y.f|)).                                                					\quad\cr}$$
Putting this over a common denominator of 
$|y.f|^2 (|y.f| - (size(y) + y.e))$ 
and cancelling equal terms (in the numerator) we are left with a second quotient, whose numerator is 
$$x.e |y.f|^2 - (-|x.f| |y.f| y.e - size(x) |y.f| (size(y) + y.e) -
|x.f| size(y) (size(y) + y.e))$$
and we see that all terms in this numerator are positive.
Further, the terms in the numerators of the first and second quotients correspond in a natural way, and each term in the numerator of the second quotient  is greater than or equal to the absolute value of its corresponding term in the numerator of the first quotient.  

Finally, because the denominator in the second quotient is less than or equal to the absolute value of the denominator in the first quotient, we see that the absolute value of the first quotient is less than or equal to the second quotient, as desired. \qed
\bigskip
We present a couple of  variations on Proposition 6.8 which will be useful when we do round-off error.
\medskip
\noindent {\bf Proposition 6.9 (division of a double by an AffApprox):} If $x$ is a double and $y$ is the AffApprox $(y.f;\ y.f_0,\ yf_1,\ y.f_2;\ y.e)$, then $r(x / y) \supseteq r(x) / r(y)$, where
$$\displaylines{\quad
x / y = (x / y.f; -x \times y.f_0 /(y.f^2), 
-x.f \times y.f_1 /(y.f^2), -x.f \times y.f_2 /(y.f^2); 
\hfill\cr
\hfill{}  
(|x|  /(|y.f| - (size(y) + y.e)) - 
\hfill\cr
\hfill{}  
(|x|/|y.f|  + |x| size(y)/(|y.f||y.f|)))
                                              \quad\cr}$$
Of course, we have to avoid division by zero.  That is, we demand that $|y.f| > size(y) + y.e.$ \qed
\bigskip
\noindent {\bf Proposition 6.10 (division of an AffApprox by a double):} If $x$ is the AffApprox $(x.f;\ x.f_0,\ x.f_1,\ x.f_2;\ x.e)$ and $y$ is a double, then $r(x / y) \supseteq r(x) / r(y)$, where
$$
x / y = (x.f / y; x.f_0 / y, 
x.f_1 / y, x.f_2 / y; 
x.e/ |y| )
$$
Of course, we have to avoid division by zero.  That is, we demand that $|y| > 0.$ \qed
\bigskip
Finally, we do the square root.
\medskip
\noindent {\bf Proposition 6.11 (square root):} If $x$ is the AffApprox $(x.f;\ x.f_0,\ x.f_1,\ x.f_2;\ x.e)$ then $r(\sqrt x) \supseteq \sqrt {r(x)}$, where
$$\displaylines{\quad
\sqrt x = (\sqrt {x.f}; 
 {x.f_0 \over 2 \sqrt {x.f}}, 
 {x.f_1 \over 2 \sqrt {x.f}}, 
 {x.f_2 \over 2 \sqrt {x.f}};
\hfill\cr
\hfill{}  
 \sqrt {|x.f|} - ({size(x) \over 2 \sqrt {|x.f|}} + \sqrt {|x.f| - (size(x) + x.e)})
                                              )\quad\cr}$$
Of course, we have to avoid division by zero, and taking the real square root of a negative number.  In particular, when $|x.f| > size(x) + x.e$ we get the formula above, but if this does not hold, then we use the crude estimate $(0;0,0,0;\sqrt{|x.f| + size(x) + x.e}).$
\medskip
\noindent {\bf Proof:} 
As above, note that for any point $(z_0, z_1, z_2) \in A$ and any function $g \in r(x)$ we can find a complex number $u$ with 
$|u|   \le 1$, such that
$$ g(z_0, z_1, z_2) = x.f + (x.f_k z_k) + (x.e) u.$$
Also, because $|x.f| > size(x) + x.e.$ we see that the argument of
$x.f + (x.f_k z_k) + (x.e) u$ is within $\pi/2$ of the argument of 
$x.f$, and therefore, we can require that $\sqrt{g(z_0,z_1,z_2)}$ has argument within $\pi/4$ of the argument of $\sqrt{x.f}.$

We need to show that 
$$|\sqrt {x.f + x.f_k z_k + (x.e) u} - (\sqrt{x.f} + {x.f_k \over 2 \sqrt {x.f}}) z_k|$$
$$\le \sqrt{|x.f|} - ({size(x) \over 2 \sqrt {|x.f|}} + \sqrt {|x.f| - (size(x) + x.e)}\ )$$
Or, after multiplying both sides by $\sqrt {|x.f|}$,
$$\sqrt {x.f (x.f + x.f_k z_k + (x.e) u)} - (x.f + (x.f_k) z_k /2)$$
$$\le (|x.f| - size(x) /2) - \sqrt {|x.f|(|x.f| - (size(x) + x.e))}$$
The two  sides of the inequality are of the form $ A - B$ and $C - D$, and we
``simplify" by multiplying by  
${ A+ B \over A + B}$ and ${ C+ D \over C + D}$.   We now show that the (absolute value of the) left-hand numerator is less than or equal to the right-hand numerator.  Later, we will show that the (absolute value of the) left-hand denominator is larger than or equal to the right-hand denominator. 
The left-hand numerator is 
$$|x.f (x.f + x.f_k z_k + (x.e) u) - (x.f + (x.f_k) z_k /2)^2|$$
$$= |x.f^ 2 + x.f (x.f_k) z_k + x.f (x.e) u - x.f^2 - x.f (x.f_k) z_k - (x.f_k^2) z_k^2/4|$$
$$= |x.f (x.e) u - (x.f_k^2) z_k^2/4|$$
The right-hand numerator is 
$$(|x.f| - size(x)/2)^2 - |x.f| (|x.f| - (size(x) + x.e))$$
$$= |x.f|^2 - |x.f| size(x) + size(x)^2/4 - |x.f|^2 + |x.f| size(x) + |x.f| x.e$$
$$= |x.f| x.e + size(x)^2/4$$
So the left-hand numerator is indeed less than or equal to the right-hand numerator.

We now compare the denominators, but only after dividing each by
$\sqrt {|x.f|}$. 
The left-hand denominator is
$$|\sqrt {x.f + x.f_k z_k + (x.e) u} + 
(\sqrt {x.f} + {x.f_k z_k\over 2 \sqrt{x.f}})|$$
while the right-hand denominator is 
$$\sqrt{|x.f|} - {size(x) \over 2 \sqrt {|x.f|}} + \sqrt {|x.f| - (size(x) + x.e)}$$
The claim that the left-hand denominator is greater than or equal to the right-hand denominator is a bit complicated.  First, compare the $\sqrt{x.f}$ term and the $\sqrt{|x.f|}$ terms.  They are the same distance from the origin.  Next, note that as $z_k$ and $u$ take on all values, $x.f + x.f_k z_k + (x.e) u$ describes a disk centered at $x.f$ and whose radius is less than $ |\sqrt {x.f}|$.  Hence, $\sqrt {x.f + x.f_k z_k + (x.e) u}$ describes a convex symmetric (about the line formed by the origin and $x.f$) set centered at $\sqrt{x.f}$.  Further, $\sqrt{x.f} + \sqrt {x.f + x.f_k z_k + (x.e) u}$ describes a convex symmetric (about the line formed by the origin and $x.f$) set centered at $2 \sqrt{x.f}$.  In any case, it is easy enough to see that no points on this convex symmetric set get closer to the origin than $\sqrt{|x.f|}  + \sqrt {|x.f| - (size(x) + x.e)}$.

Finally, because $|{x.f_k z_k \over 2 \sqrt{x.f}}| \le {size(x) \over 2 \sqrt {|x.f|}},$ no points of $$\sqrt {x.f} + \sqrt {x.f + x.f_k z_k + (x.e) u} + 
{x.f_k z_k \over 2 \sqrt{x.f}}$$
can get closer to the origin than
$$\sqrt{|x.f|} + \sqrt {|x.f| - (size(x) + x.e)} - 
{size(x) \over 2 \sqrt {|x.f|}}$$ \qed
\bigskip \noindent {\bf Example 6.7 (Continuation of Example 5.4):}
We repeat the description of the sub-box under investigation.

The sub-box $Z(s01011)$ with  
$$s = 001000110001110111001111000101111111101111100111001111000001111011110111$$ 
\vfill\eject
\noindent is the region where
\medskip
$$\left(\matrix{
-1.381589027741\ldots  \le  Re(L')  \le  -1.379848991182\ldots\cr
-1.378124546093\ldots  \le  Re(D')  \le  -1.376574349753\ldots\cr
0.999893182771\ldots  \le Re(R')  \le  1.001274250703\ldots\cr
-2.535837191243\ldots  \le  Im(L')  \le  -2.534606799593\ldots\cr
2.535404997792\ldots  \le  Im(D')  \le  -2.534308843448\ldots\cr
-0.001953125000\ldots  \le  Im(R')  \le  0.000000000000\ldots\cr
}\right)$$
\bigskip
For this sub-box, we get
$$
f = \left[\matrix{
  \left(\matrix{
    -0.8677851121 + i 1.4607429651;\cr
    0.0000248810 - i 0.0003125810,\cr
    0.0000000000 + i 0.0000000000,\cr
    0.0000000000 + i 0.0000000000;\cr
    0.0000000289
  }\right)
 &
  \left(\matrix{
    0.0000000000 + i 0.0000000000;\cr
    0.0000000000 + i 0.0000000000,\cr
    0.0000000000 + i 0.0000000000,\cr
    0.0000000000 + i 0.0000000000;\cr
    0.0000000000
  }\right)
 \cr
  \left(\matrix{
    0.0000000000 + i 0.0000000000;\cr
    0.0000000000 + i 0.0000000000,\cr
    0.0000000000 + i 0.0000000000,\cr
    0.0000000000 + i 0.0000000000;\cr
    0.0000000000
  }\right)
 &
  \left(\matrix{
    -0.3006023265 - i 0.5060039953;\cr
    -0.0000909686 - i 0.0000593570,\cr
    0.0000000000 + i 0.0000000000,\cr
    0.0000000000 + i 0.0000000000;\cr
    0.0000000301
  }\right)
}\right]
$$
and
$$
w = \left[\matrix{
  \left(\matrix{
    -0.5845111829 + i 0.4773282853;\cr
    0.0000000000 + i 0.0000000000,\cr
    -0.0000296707 - i 0.0001657332,\cr
    -0.0004345111 - i 0.0001209539;\cr
    0.0000002590
  }\right)
 &
  \left(\matrix{
    -0.2840228472 + i 0.9825063583;\cr
    0.0000000000 + i 0.0000000000,\cr
    0.0000516606 - i 0.0001128245,\cr
    0.0005776611 - i 0.0001998632;\cr
    0.0000006462
  }\right)
 \cr
  \left(\matrix{
    -0.2832291572 + i 0.9833572297;\cr
    0.0000000000 + i 0.0000000000,\cr
    0.0000515806 - i 0.0001129408,\cr
    -0.0005778031 + i 0.0002005440;\cr
    0.0000002806
  }\right)
 &
  \left(\matrix{
    -0.5846352333 + i 0.4764792236;\cr
    0.0000000000 + i 0.0000000000,\cr
    -0.0000294917 - i 0.0001656653,\cr
    0.0004341392 + i 0.0001213070;\cr
    0.0000005286
  }\right)
}\right].$$
calculating
$g = f^{-1}wf^{-1}w^{-1}f^{-1}w^{-1}fw^{-1}f^{-1}w^{-1}f^{-1}wf^{-1}wfww$ gives
$$
g = \left[\matrix{
  \left(\matrix{
    -0.5764337542 + i 0.4752708071;\cr
    -0.0031657223 - i 0.0001436786,\cr
    -0.0017723577 + i 0.0000352928,\cr
    -0.0011623491 + i 0.0017516088;\cr
    0.0008229225
  }\right)
 &
  \left(\matrix{
    -0.2704033973 + i 0.9822741250;\cr
    -0.0045902952 - i 0.0019135041,\cr
    -0.0026219461 - i 0.0007506230,\cr
    -0.0002823450 + i 0.0033805602;\cr
    0.0008037640
  }\right)
 \cr
  \left(\matrix{
    -0.2861207992 + i 0.9766064999;\cr
    -0.0002777968 + i 0.0020330488,\cr
    0.0000837571 + i 0.0010241875,\cr
    0.0028322367 - i 0.0005972336;\cr
    0.0018172437
  }\right)
 &
  \left(\matrix{
    -0.5861133046 + i 0.4624368851;\cr
    -0.0021932627 + i 0.0040523411,\cr
    -0.0008612361 + i 0.0022394639,\cr
    0.0061581377 - i 0.0005862070;\cr
    0.0017738513
  }\right)
}\right].
$$
We then get
$$
length(g) = \left(\matrix{
    -1.3588762105 - i 2.4897230182;\cr
    0.0030210500 - i 0.0182284729,\cr
    0.0007938572 - i 0.0096614614,\cr
    -0.0122034521 + i 0.0074353043;\cr
    0.0080071969
  }\right)
$$
and
$$
{length(g) \over L'} =
  \left(\matrix{
    0.9825397896 - i 0.0008933519;\cr
    0.0053701602 + i 0.0037789019,\cr
    0.0028076072 + i 0.0018421952,\cr
    -0.0002400615 - i 0.0049443045;\cr
    0.0027802966
  }\right).
$$

This isn't quite good enough to kill the sub-box, since
$|length(g)/L'|$ can be high as $1.0001951323$.

When we subdivide $Z(s01011)$, we have to analyze two sub-boxes,
$Z(s010110)$ and $Z(s010111)$.
For $Z(s010110)$, the same calculation on the region
{\obeylines
$-1.381589027741073400 \leq Re(L') \leq -1.379848991182205200$
$-1.378124546093485700 \leq Re(D') \leq -1.376574349753672900$
$0.999893182771602220 \leq Re(R') \leq 1.001274250703607400$
$-2.535837191243490300 \leq Im(L') \leq -2.534606799593201600$
$-2.535404997792558600 \leq Im(D') \leq -2.534308843448505900$
$-0.001953125000000000 \leq Im(R') \leq -0.000976562500000000$
}
gives
$$
{length(g) \over L'} = 
  \left(\matrix{
    0.9814518667 + i 0.0008103446;\cr
    0.0053616729 + i 0.0037834001,\cr
    0.0028027236 + i 0.0018435245,\cr
    -0.0013175066 - i 0.0032448794;\cr
    0.0019033926
  }\right),
$$
and we can then bound $|{length(g) \over L'}| \leq 0.9967745579$,
which kills $Z(s010110)$.

On $Z(s010111)$, the calculation gives
$$
{length(g) \over L'} =
  \left(\matrix{
    0.9836225919 - i 0.0025990177;\cr
    0.0053786346 + i 0.0037743930,\cr
    0.0028124892 + i 0.0018408583,\cr
    -0.0013333182 - i 0.0032343347;\cr
    0.0019044429
  }\right)
$$
and $|{length(g) \over L'}| \leq 0.9989610507$,
which kills $Z(s010111)$.

\sectionskip\magnification=1200
\baselineskip = 0.2 true in
\def\qed{{\bf \ \ \ \ \vrule height6pt width5pt depth1pt}}

\noindent {\bf Chapter 7: Complex Numbers with Round-Off Error}
\medskip \noindent {\bf Remark 7.1:}
The theoretical method for proving Theorem 0.2 has been implemented on the computer programs {\it verify} and {\it fudging} given in the Appendices.  To make this computer-aided proof rigorous, we need to deal with round-off error in calculations.  

One approach to round-off error would be to use interval arithmetic packages to carry out all calculations with floating-point numbers (also called  ``doubles"), or to generate our own version of these packages.  However, it appears that this would be much too slow given the size of our collection of sub-boxes and conditions and killer words.  

To solve this problem of speed, we implement round-off error at a higher level of programing.  That is, we incorporate round-off error directly into AffApprox's.  This necessitates that we incorporate it into complex numbers as well.
\medskip \noindent {\bf Definition 7.2:}
In the next Chapter we work with AffApprox's.
In this section we show how to do standard operations on complex numbers while keeping track of round-off error.  There are two types of complex numbers to consider: 

1.)  An {\it XComplex} corresponds to a complex number that is represented exactly.  Thus, it simply consists of a real part and an imaginary part.

2.)  An {\it AComplex} corresponds to an ``interval'' that contains the complex number in question.  Thus, it consists of an XComplex and a floating-point number representing the error.  In particular, $(x;e)$ represents the set of complex numbers $\{w\ :\ |w - x| \le e\ \}$.
Following the notation  of Chapter 6, we could now denote this set of complex numbers by $r(x,e),$ but instead we suppress mention of the $r$ functions throughout the rest of this section.  It seems preferable to abuse notation in this fashion in the interests of simplicity.
\medskip \noindent {\bf Remark 7.3:}
In general,  our operations  act on XComplexes and produce AComplexes, or they act on AComplexes and produce AComplexes  .  In one case, the unary minus, an XComplex goes to an XComplex.  
In the calculations that follow the effect on the error is the whole point.
\medskip \noindent {\bf Conventions 7.4:}
We begin, by writing down our basic rules, which follow easily from the IEEE-754 double-precision standard for machine arithmetic.  (Actually,  the ``hypot" function $h(a,b)$, which computes by elaborate chicanery $\sqrt{a^2 + b^2},$ is not part of the IEEE standard, but  satisfies the appropriate standard according to the documentation provided (see [Kahan]).)  The operations here are on
double-precision floating-point real numbers (doubles) and we denote a true operation by the usual symbol and the associated machine operation by the same symbol in a circle, with two exceptions: a machine square root $\sqrt a$ is denoted $\root o \of a$ and the machine version of the hypot function is denoted $h_\circ$.  Perhaps a third exception is our occasional notation of true multiplication by the absence of a symbol.  The standard number used in error analysis is an ``EPS.''   It depends on the number of bits used to store floating-point numbers, and this can vary from machine to machine.  Finally, it should be noted that our analysis breaks down when there is ``underflow," so in the computer programs we ask the computer to inform us if underflow has occurred.  

As in Chapter 6, we now break with the usual numbering convention.
\bigskip
\noindent {\bf Basic Properties 7.0 (assuming no underflow):}
\medskip
$1 + k \times EPS = 1 \oplus (k \otimes EPS)$ when $k$ is an integer which is not huge in absolute value.

$2^k \times A = 2^k \otimes A$ when $k$ is an integer, $2^k \otimes A$ is not infinity.

\bigskip
$$|(a + b) - (a \oplus b)| \leq (EPS/2) |a + b|$$
$$|(a + b) - (a \oplus b)| \leq (EPS/2) |a \oplus b|$$
$$|(a - b) - (a \ominus b)| \leq (EPS/2) |a - b|$$
$$|(a - b) - (a \ominus b)| \leq (EPS/2) |a \ominus b|$$
$$|(a \times b) - (a \otimes b)| \leq (EPS/2) |a \times b|$$
$$|(a \times b) - (a \otimes b)| \leq (EPS/2) |a \otimes b|$$
$$|(a / b) - (a \oslash b)| \leq (EPS/2) |a / b|$$
$$|(a / b) - (a \oslash b)| \leq (EPS/2) |a \oslash b|$$
$$|\sqrt a - \root o \of a| \leq (EPS/2) |\sqrt a|$$
$$|\sqrt a - \root o \of a| \leq (EPS/2) |\root o \of a|$$
$$|h(a,b) - h_\circ(a,b)|\leq (EPS) |h(a,b)|$$
$$|h(a,b) - h_\circ(a,b)|\leq (EPS) |h_\circ(a,b)|$$
\medskip
From these formulas, we immediately compute the following.
$$(1 - EPS/2) |a + b| \leq |a \oplus b| 
\leq (1 + EPS/2) |a+b|$$
$$(1 - EPS/2) |a \oplus b| \leq |a + b| 
\leq (1 + EPS/2) |a \oplus b|$$
$$(1 - EPS/2) |a - b| \leq |a \ominus b| 
\leq (1 + EPS/2) |a - b|$$
$$(1 - EPS/2) |a \ominus b| \leq |a - b| 
\leq (1 + EPS/2) |a \ominus b|$$
$$(1 - EPS/2) |a \times b| \leq |a \otimes b| 
\leq (1 + EPS/2) |a \times b|$$
$$(1 - EPS/2) |a \otimes b| \leq |a \times b| 
\leq (1 + EPS/2) |a \otimes b|$$
$$(1 - EPS/2) |a / b| \leq |a \oslash b| 
\leq (1 + EPS/2) |a / b|$$
$$(1 - EPS/2) |a \oslash b| \leq |a / b| 
\leq (1 + EPS/2) |a \oslash b|$$
$$(1 - EPS/2) |\sqrt a| \leq |\root o \of a| 
\leq (1 + EPS/2) |\sqrt a|$$
$$(1 - EPS/2) |\root o \of a|\leq |\sqrt a| 
\leq (1 + EPS/2) |\root o \of a|$$
$$(1 - EPS) |h(a,b)| \leq |h_\circ (a,b)| 
\leq (1 + EPS) |h(a,b)|$$
$$(1 - EPS) |h_\circ (a,b)| \leq |h(a,b)| 
\leq (1 + EPS) |h_\circ (a,b)|$$
\bigskip
Of course, we can also get the following type of formula, which is sometimes convenient, for example, in the proof of Lemma 7.2.
\medskip
$$({1 \over {1 + {EPS \over 2}}}) |a \oplus b| \leq |a + b| 
\leq ({1 \over {1 - {EPS \over 2}}}) |a \oplus b|$$
\medskip
\noindent {\bf Proof:}  There is a finite set of numbers which are representable on the computer, and the IEEE standard states that the result of an operation is always the closest representable number to the true solution.  Ignoring technicalities,  a non-zero floating-point number is represented by a fixed number of bits of which the first determines the sign of the number, the next $m$ represent the exponent, and the remaining $n$ represent the mantissa of the number.  Because our non-zero numbers start with a 1, that means the $n$ mantissa bits actually represent the next $n$ binary digits after the 1.  That is, the mantissa is actually $1.b_1b_2b_3...b_{n}.$  With this divyying up of the $m + n + 1$ bits among exponent and mantissa,  $EPS$ would be $2^{-n}$ and $EPS/2$ would be $2^{-(n + 1)}.$ 

Given this set-up, properties of the form
$$|(a + b) - (a \oplus b)| \leq (EPS/2) |a + b|$$
follow immediately.  Then,
$$|(a + b) - (a \oplus b)| \leq (EPS/2) |a \oplus b|$$
follows because the true answer has ``exponent" which is less than or equal to the exponent of the machine answer. \qed
\bigskip
Before starting in with our Propositions, we prove  a couple of lemmas.
\bigskip
\noindent {\bf Lemma 7.0:} 
$$(1 - EPS) \otimes |a \oplus b| \le |a + b| \le (1 + EPS) \otimes |a \oplus b|$$
\centerline {\it Analagous formulas hold for $\ -,\ *,\ /,\ \sqrt.$}
\bigskip
{\bf Proof:}
Assume $a+b > 0$.  
If $(1+EPS) \otimes (a \oplus b) < (a+b)$ then the machine number 
$(1+EPS) \otimes (a \oplus b)$ is a better approximation to $a+b$ than $a \oplus b$, because $(a \oplus b) < (1+EPS) \otimes (a \oplus b)$.  This contradicts the IEEE standard.  The case $a+b < 0$ can be handled similarly, and the case $a+b = 0$ is trivial.  Similarly for the left-hand inequality. \qed
\bigskip
\noindent {\bf Lemma 7.1:} $$(1 + EPS/2)^a A \le (1 + k EPS) \otimes A$$ 
{\it where $A \ge 0$, and $a$ and $k$ are (not huge) integers, such that for $a$ even, $k = {a \over 2} + 1$ and 
for $a$ odd, $k = {a+1 \over 2} + 1$.}
\bigskip
{\bf Proof:} 
$$(1 + EPS/2)^a A\le (1 - EPS/2)(1 + k EPS) A \le (1 + k EPS) \otimes A$$ 
The first inequality holds if $a$ and $k$ are as in the Lemma, and the second inequality is a consequence of one of the formulas preceding Lemma 7.0 ($A \ge 0$). \qed
\bigskip 
We start the operations.  We will give proofs for most, the others should be straightforward to derive.  Note that for an XComplex, $x = (x.re,x.im)$, and for an AComplex, $x = (x.re,x.im,x.e).$
\bigskip
\noindent {\bf Proposition 7.1 (-X): } 

$-x = (-x.re,-x.im)$. \qed
\bigskip
\noindent {{\bf Proposition 7.2 (X + D)} (an XComplex and a double added together, which yields an AComplex):}

 $x + d = (x.re \oplus d, x.im;
(EPS/2)\otimes |x.re \oplus d|)$
\bigskip
{\bf Proof:}
The  error is bounded by  $$|(x.re + d) - (x.re \oplus d)|$$
$$\le (EPS/2) | x.re \oplus d|$$
$$  =  (EPS/2) \otimes | x.re \oplus d|$$ \qed
\bigskip
\noindent {{\bf Proposition 7.3 (X - D)} (a double subtracted from an XComplex, which yields an AComplex):}

 $x - d = (x.re \ominus d, x.im;
(EPS/2)\otimes |x.re \ominus d|)$ \qed
\bigskip
\noindent {{\bf Proposition 7.4 (X + X)} (an XComplex and an XComplex added together, which yields an AComplex):}

 $x + y = (x.re \oplus y.re, x.im \oplus y.im;
(EPS/2)\otimes ((1 + EPS) \otimes (|x.re \oplus y.re| \oplus |x.im \oplus y.im|)))$
\bigskip
{\bf Proof:}
The error is bounded by
$$|(x.re + y.re) - (x.re \oplus y.re)| + 
  |(x.im + y.im) - (x.im \oplus y.im)|$$
$$\le (EPS/2)(| x.re \oplus y.re| + 
   |x.im \oplus y.im|)$$
$$\le (EPS/2)((1 + EPS) \otimes (| x.re \oplus y.re| \oplus 
   |x.im \oplus y.im|))$$
$$= (EPS/2) \otimes ((1 + EPS) \otimes (| x.re \oplus y.re| \oplus 
   |x.im \oplus y.im|))$$
To go from line 2 to line 3 we used Lemma 7.0. \qed
\bigskip
\noindent {\bf Proposition 7.5 (X - X):}

 $x + y = (x.re \ominus y.re, x.im \ominus y.im;
(EPS/2)\otimes ((1 + EPS) \otimes (|x.re \ominus y.re| \oplus |x.im \ominus y.im|)))$ \qed
\bigskip
\noindent {{\bf Proposition 7.6 (A + A)} (an AComplex and an AComplex added together, which yields an AComplex):}

$x+y = (re,im;e)$ where 

$re = x.re \oplus y.re$

$im = x.im \oplus y.im$
 
$e = (1 + 2 EPS) \otimes ( ((EPS/2) \otimes (|re| \oplus |im|))
 \oplus (x.e \oplus y.e))$
\bigskip
{\bf Proof:}
The error is bounded by the sum of the contributions from the real part, the imaginary part, and the two individual errors:
$$|(x.re \oplus y.re) - (x.re + y.re)| + |(x.im \oplus y.im) - (x.im + y.im)| + (x.e + y.e).$$  
$$\le (EPS/2) |x.re \oplus y.re| + (EPS/2) |x.im \oplus y.im| + 
(1 + EPS/2) (x.e \oplus y.e)$$
$$\le  (1 + EPS/2) (EPS/2)(|x.re \oplus y.re| \oplus 
 |x.im \oplus y.im| ) + (1 + EPS/2) (x.e \oplus y.e)$$
$$= (1 + EPS/2) ((EPS/2)(|x.re \oplus y.re| \oplus 
 |x.im \oplus y.im| ) + (x.e \oplus y.e))$$
$$\le (1 + EPS/2)^2 (((EPS/2) (|x.re \oplus y.re| \oplus 
 |x.im \oplus y.im| )) \oplus  (x.e \oplus y.e))$$
$$\le   (1 + 2 EPS) \otimes (((EPS/2)\otimes (|x.re \oplus y.re| \oplus |x.im \oplus y.im| )) \oplus (x.e \oplus y.e))$$  \qed

The hierarchy for machine operations is the same as that for true operations, so one pair of parentheses is unnecessary and will often be omitted in what follows.
\bigskip
\noindent {{\bf Proposition 7.7 (A - A)}:}

$x+y = (re,im;e)$ where 

$re = x.re \ominus y.re$

$im = x.im \ominus y.im$
 
$e = (1 + 2 EPS) \otimes ( ((EPS/2) \otimes (|re| \oplus |im|))
 \oplus (x.e \oplus y.e))$ \qed
\bigskip
\noindent {\bf Proposition 7.8 (X $\times$ D):} 

 $x \times d = (re,im;e)$ where 

$re = x.re \otimes y$

$im = x.im \otimes y$
 
$e = (EPS/2) \otimes ((1 + EPS) \otimes (|re| \oplus |im|) )$
\bigskip
{\bf Proof:} The error is bounded by
$$|(x.re \times  y) - (x.re \otimes y)| + 
  |(x.im \times y) - (x.im \otimes y)|$$
$$\le (EPS/2)|x.re \otimes  y| + (EPS/2)|x.im \otimes y|$$
$$= (EPS/2)(|x.re \otimes  y| + |x.im \otimes y|)$$
$$\le (EPS/2) \otimes ((1 + EPS) \otimes (|x.re \otimes  y| \oplus 
  |x.im \otimes y|))$$ \qed
\bigskip
\noindent {\bf Proposition 7.9 (X / D):} 

 $x / d = (re,im;e)$where 

$re = x.re \oslash y$

$im = x.im \oslash y$
 
$e = (EPS/2) \otimes ((1 + EPS) \otimes (|re| \oplus |im|) )$ \qed
\bigskip
\noindent {\bf Proposition 7.10 (X $\times$ X):}

$x \times y = (re,im;e)$ where 

$re = re1 \ominus re2$, with $re1 = x.re \otimes y.re$ and $re2 = x.im \otimes y.im$

$im = im1 \oplus im2$, with $im1 = x.re \otimes y.im$ and $im2 = x.im \otimes y.re$
 
$e = EPS \otimes ((1 + 2EPS) \otimes  ((|re1| \oplus |re2|) \oplus  (|im1| \oplus |im2|)))$
\bigskip
{\bf Proof:}
The error is bounded by the sum of the contributions from the real part and the imaginary part:
$$\displaylines{\quad|(x.re \times  y.re - x.im \times y.im) - ((x.re \otimes y.re) \ominus (x.im \otimes y.im))| \hfill\cr
\hfill{}+ |(x.re \times  y.im + x.im \times y.re) - ((x.re \otimes y.re) \oplus (x.im \otimes y.im))|\quad\cr}$$  
We want to bound this by a machine formula.  Let's begin by bounding 
$$|(x.re \times  y.re - x.im \times y.im) - ((x.re \otimes y.re) \ominus (x.im \otimes y.im))| $$ 
by a machine formula.
$$|(x.re \times  y.re - x.im \times y.im) - ((x.re \otimes y.re) \ominus (x.im \otimes y.im))| $$
$$\displaylines{\quad\le |((x.re \times  y.re) - (x.im \times y.im)) - ((x.re \otimes y.re) - (x.im \otimes y.im))| \hfill\cr
\hfill{}+ |((x.re \otimes y.re) - (x.im \otimes y.im))-((x.re \otimes y.re) \ominus (x.im \otimes y.im))| \quad\cr}$$
$$\displaylines{\quad\le |(x.re \times  y.re) - (x.re \otimes y.re)| + |(x.im \times y.im)- (x.im \otimes y.im)| \hfill\cr
\hfill{}+ (EPS/2)|(x.re \otimes y.re) - (x.im \otimes y.im)| \quad\cr}$$
$$\displaylines{\quad\le (EPS/2)|(x.re \otimes y.re)| + (EPS/2)|(x.im \otimes y.im)| \hfill\cr
\hfill{}+ (EPS/2)(|x.re \otimes y.re| + |x.im \otimes y.im|)\quad\cr}$$
$$= (EPS/2)(2) (|x.re \otimes y.re| + |x.im \otimes y.im|)$$
$$\le EPS (1 + EPS/2) (|x.re \otimes y.re| \oplus |(x.im \otimes y.im|) $$

Almost the exact same calculation produces the analagous formula for the imaginary contribution, and we now combine the two to get a bound on the total error.  

$$\displaylines{\quad\le EPS (1 + EPS/2) (|x.re \otimes y.re| \oplus |x.im \otimes y.im|)  \hfill\cr
\hfill{}+
EPS (1 + EPS/2) (|x.re \otimes y.im| \oplus |x.im \otimes y.re| )\quad\cr}$$
$$\displaylines{\quad= EPS (1 + EPS/2) ((x.re \otimes y.re| \oplus |x.im \otimes y.im|)  \hfill\cr
\hfill{}+
(|x.re \otimes y.im| \oplus |x.im \otimes y.re|) )\quad\cr}$$
$$\displaylines{\quad\le EPS (1 + EPS/2)^2 ((|x.re \otimes y.re| \oplus |x.im \otimes y.im|)  \hfill\cr
\hfill{} \oplus
(|x.re \otimes y.im| \oplus |x.im \otimes y.re| ))\quad\cr}$$
$$\displaylines{\quad\le EPS \otimes ((1 + 2EPS)\otimes((|x.re \otimes y.re| \oplus |x.im \otimes y.im| ) \hfill\cr
\hfill{}\oplus
(|x.re \otimes y.im| \oplus |x.im \otimes y.re| )))\quad\cr}$$ \qed
\bigskip
\noindent {\bf Proposition 7.11 (D / X):}

$x  / y  = (re,im;e)$ where 

$re = (x  \otimes y.re)\oslash nrm$ where $nrm = y.re \otimes y.re \oplus y.im \otimes y.im$

$im = -(x  \otimes y.im) \oslash nrm$
 
$e = (2EPS) \otimes ((1 +  2EPS)\otimes(|re| \oplus |im|))$
\bigskip
{\bf Proof:}
The true version of $x/y$ is equal to  $(x \times y.re + i (-x \times y.im))/ (y.re^2 + y.im^2)$ and we need to compare this with the machine version to find the error.  Further, this error is less than or equal to the sum of the real error and the imaginary error. Thus, we start with the real calculation (as above, we use $nrm$ to represent the machine version of $y.re^2 + y.im^2$).

$$|{x \times y.re  \over y.re^2 + y.im^2} - ((x \otimes y.re)\oslash nrm)|$$
$$\le |(x \otimes y.re) \oslash nrm - {x \otimes y.re \over nrm}| +
|{x \otimes y.re \over nrm} - {x \times y.re \over nrm}| +
|{x \times y.re \over nrm} - {x \times y.re  \over y.re^2 + y.im^2}|$$

Before continuing, let's compare ${1 \over nrm}$ and ${1 \over y.re^2 + y.im^2}$ by developing a formula for comparing ${1 \over a^2 + b^2}$ and its associated ${1 \over nrm}$:
\bigskip
\noindent {\bf Lemma 7.2:} $$|{1 \over nrm} - {1 \over a^2 + b^2}| \le   (EPS + (EPS/2)^2) {1 \over nrm}$$
{\it where} $nrm = a\otimes a \oplus b \otimes b.$
\medskip
{\bf Proof:}
We compute that 
$$({1 \over 1 + EPS/2})^2 \times nrm
\le a^2 + b^2 
\le ({1 \over 1 - EPS/2})^2 \times nrm,$$ hence 
$${1 \over nrm} (1 - EPS/2)^2 
\le {1 \over a^2 + b^2}
\le {1 \over nrm} (1 + EPS/2)^2.$$  It then follows that 
$$|{1 \over nrm} - {1 \over a^2 + b^2}| \le 
    {1 \over nrm} (1 + EPS/2)^2 - {1 \over nrm}$$
 $$= {1 \over nrm} ((1 + EPS/2)^2 - 1) =
  (EPS + (EPS/2)^2){1 \over nrm} $$  \qed
\bigskip
Getting back to our main calculation (with $nrm = y.re \otimes y.re \oplus y.im \otimes y.im$),
$$|(x \otimes y.re) \oslash nrm - {x \otimes y.re \over nrm}| +
|{x \otimes y.re \over nrm} - {x \times y.re \over nrm}| +
|{x \times y.re \over nrm} - {x \times y.re  \over y.re^2 + y.im^2}|$$
$$\le(EPS/2){|x \otimes y.re| \over nrm} +
        (EPS/2){|x \otimes   y.re| \over nrm} +
(EPS + (EPS/2)^2) {|x \times   y.re| \over nrm}$$
$$=(EPS/2)({1 \over nrm})(2 |x \otimes y.re| +
(2 + EPS/2) \times |x \times y.re| )$$
$$\le (EPS/2)({1 \over nrm})(2 |x \otimes y.re| +
(2 + EPS/2)(1 + EPS/2) \times |x \otimes y.re| )$$
$$= (EPS/2)({1 \over nrm})(|x \otimes y.re|)
(2 + (2 + EPS/2)(1+EPS/2))$$
$$\le (EPS/2)(4 + 3EPS/2 + (EPS/2)^2)(|x \otimes y.re|)({1 \over nrm})$$
$$\le (EPS/2)(4 + 3EPS/2 + (EPS/2)^2)(1+EPS/2)(|x \otimes y.re| \oslash nrm)$$
$$\le (2EPS)(1 + 3EPS/8 + (EPS/4)^2)(1+EPS/2)(|(x \otimes y.re \oslash nrm)|)$$
We also get the analagous formula for the imaginary contribution for the error, so our total error is bounded by
$$(2EPS)(1 + 3EPS/8 + (EPS/4)^2)(1+EPS/2)((|(x \otimes y.re)\oslash nrm|)+(|(x \otimes y.im)\oslash nrm|))$$
$$\le (2EPS)(1 + 3EPS/8 + (EPS/4)^2)(1+EPS/2)^2((|(x \otimes y.re)\oslash nrm|) \oplus (|(x \otimes y.im)\oslash nrm|))$$
$$\le (2EPS)(1 - EPS/2)(1+ 2EPS)((|(x \otimes y.re)\oslash nrm|) \oplus (|(x \otimes y.im)\oslash nrm|))$$
$$\le (2EPS) \otimes ((1+ 2EPS) \otimes ((|(x \otimes y.re)\oslash nrm|) \oplus (|(x \otimes y.im)\oslash nrm|)))$$

Here we used the fact that 
$$(1 + 3EPS/8 + (EPS/4)^2)(1+EPS/2)^2 \le (1 - EPS/2)(1+ 2EPS)$$ \qed
\bigskip

\noindent {\bf Proposition 7.12 (X / X):}

$x  / y  = (re,im;e)$ where 

$re = (x.re \otimes y.re \oplus x.im \otimes y.im) \oslash nrm$ where $nrm = y.re \otimes y.re \oplus y.im \otimes y.im$

$im = (x.im \otimes y.re \ominus x.re \otimes y.im) \oslash nrm$
 
$e = (5EPS/2) \otimes ((1 +  3EPS) \otimes A)$  where
$$A = 
((|x.re \otimes y.re| \oplus |x.im \otimes y.im|)
 \oplus 
 (|x.im \otimes y.re| \oplus |x.re \otimes y.im|))
\oslash nrm
$$
\bigskip
{\bf Proof:}  

It will be useful to begin by comparing $(a \times b + c \times d) / (e \times e + f \times f)$ with $(a \otimes b \oplus c \otimes d) \oslash nrm.$
\medskip
\noindent {\bf Lemma 7.3:}
$$|(a \otimes b \oplus c \otimes d) \oslash nrm
- (a \times b + c \times d) / e \times e + f \times f)|$$
$$\le (EPS/2) (1+EPS/2)(5 + 3EPS/2 + (EPS/2)^2)(|a \otimes b| \oplus |c \otimes d|) / nrm$$

{\bf Proof:}
At one point in the proof, we will use the fact that $|A \oplus B| \le |A| \oplus |B|$.  This fact is easily seen by considering the various possibilities for the sign of $A$ versus the sign of $B$.  At a different point in the proof, Lemma 7.2 with $a,b$ replaced by $e,f$ will be used. 

$$|(a \otimes b \oplus c \otimes d) \oslash nrm
- (a \times b + c \times d) / e \times e + f \times f)|$$
$$\displaylines{\quad\le |(a \otimes b \oplus c \otimes d) \oslash nrm
- (a \otimes b \oplus c \otimes d) / nrm| \hfill\cr
\hfill{}+
|(a \otimes b \oplus c \otimes d) / nrm
- (a \times b + c \times d) / (e \times e + f \times f)|\quad\cr}$$
$$\displaylines{\quad\le (EPS/2) |(a \otimes b \oplus c \otimes d) / nrm| +
 |(a \otimes b \oplus c \otimes d) / nrm - (a \otimes b + c \otimes d) / nrm|
\hfill\cr
\hfill{}+|(a \otimes b + c \otimes d) / nrm -
(a \times b + c \times d) / (e \times e + f \times f)|\quad\cr}$$
$$\displaylines{\quad\le (EPS/2) |a \otimes b \oplus c \otimes d| / nrm + 
|(a \otimes b \oplus c \otimes d) - (a \otimes b + c \otimes d)| / nrm \hfill\cr
\hfill{}+
|(a \otimes b + c \otimes d) / nrm -
 (a \times b + c \times d) / nrm| \hfill\cr
\hfill{}+
|(a \times b + c \times d) / nrm - 
(a \times b + c \times d) / (e \times e + f \times f)|\quad\cr}$$
$$\displaylines{\quad\le (EPS/2) |a \otimes b \oplus c \otimes d| / nrm + 
(EPS/2) |a \otimes b \oplus c \otimes d| / nrm 
\hfill\cr
\hfill{}+
|(a \otimes b + c \otimes d) / nrm -
 (a \times b + c \times d) / nrm| \hfill\cr
\hfill{}+
|(a \times b + c \times d) / nrm - 
(a \times b + c \times d) / (e \times e + f \times f)|\quad\cr}$$
$$\displaylines{\quad\le (EPS/2) |a \otimes b \oplus c \otimes d| / nrm + 
(EPS/2) |a \otimes b \oplus c \otimes d| / nrm 
\hfill\cr
\hfill{}+
|(a \otimes b - a \times b) / nrm -
 ( c \otimes d - c \times d) / nrm| 
\hfill\cr
\hfill{}+
(EPS + (EPS/2)^2)|(a \times b + c \times d)| / nrm\quad\cr}$$
$$\displaylines{\quad\le (EPS) |a \otimes b \oplus c \otimes d| / nrm +
(EPS/2)|a \otimes b| / nrm +
(EPS/2)|c \otimes d|/ nrm  \hfill\cr
\hfill{}+
(EPS + (EPS/2)^2)(1+EPS/2)|(a \otimes b + c \otimes d)| / nrm\quad\cr}$$
$$\displaylines{\quad\le (EPS) |a \otimes b \oplus c \otimes d| / nrm +
(EPS/2)|a \otimes b| / nrm +
(EPS/2)|c \otimes d| / nrm  \hfill\cr
\hfill{}+
(EPS + (EPS/2)^2)(1+EPS/2)(|a \otimes b| + |c \otimes d|) / nrm\quad\cr}$$
$$\displaylines{\quad\le (EPS) |a \otimes b \oplus c \otimes d| / nrm 
\hfill\cr
\hfill{}+
(EPS/2 + (EPS + (EPS/2)^2)(1+EPS/2)) (|a \otimes b| + |c \otimes d|) / nrm\quad\cr}$$
$$\displaylines{\quad\le (EPS) |a \otimes b \oplus c \otimes d| / nrm 
\hfill\cr
\hfill{}+
(EPS/2)(1 + (2 + EPS/2)(1+EPS/2)) (|a \otimes b| + |c \otimes d|) / nrm\quad\cr}$$
$$\displaylines{\quad\le (EPS) (|a \otimes b| \oplus |c \otimes d| )/ nrm 
\hfill\cr
\hfill{}+
(EPS/2)(3 + 3EPS/2 + (EPS/2)^2) (|a \otimes b| + |c \otimes d|) / nrm\quad\cr}$$
$$\displaylines{\quad\le (EPS) (1+EPS/2)(|a \otimes b| \oplus |c \otimes d|) / nrm 
\hfill\cr
\hfill{}+
(EPS/2)(3 + 3EPS/2 + (EPS/2)^2)(1+EPS/2) (|a \otimes b| \oplus |c \otimes d|) / nrm\quad\cr}$$
$$\le (EPS/2)(1+EPS/2)(2 + (3 + 3EPS/2 + (EPS/2)^2)) (|a \otimes b| \oplus |c \otimes d|) / nrm$$
$$\le (EPS/2) (1+EPS/2)(5 + 3EPS/2 + (EPS/2)^2)(|a \otimes b| \oplus |c \otimes d|) / nrm$$ \qed
\bigskip

We now use this lemma to get the error term for $X/X$. Of course, this error is less than the sum of the real and imaginary errors.  Also, we let $$A=
((|x.re \otimes y.re| \oplus |x.im\otimes y.im|) 
 \oplus
(|x.im \otimes y.re| \oplus |(-x.re)\otimes y.im|)) \oslash nrm$$
\medskip
$$\displaylines{\quad |(x.re \otimes y.re \oplus x.im\otimes y.im) \oslash nrm
- (x.re \times y.re + x.im\times y.im) / (y.re \times y.re + y.im \times y.im)|
+\hfill\cr\hfill{}
|(x.im \otimes y.re \oplus (-x.re)\otimes y.im) \oslash nrm
- (x.im \times y.re + (-x.re)\times y.im) / (y.re \times y.re + y.im \times y.im)|
\quad\cr}$$
$$\displaylines{\quad\le
(EPS/2) (1+EPS/2)(5 + 3EPS/2 + (EPS/2)^2)
((|x.re \otimes y.re| \oplus |x.im\otimes y.im|) 
\hfill\cr
\hfill{} +
(|x.im \otimes y.re| \oplus |x.re\otimes y.im|)) / nrm
\quad\cr}$$
$$\displaylines{\quad\le(EPS/2) (1+EPS/2)^3 (5 + 3EPS/2 + (EPS/2)^2)
(((|x.re \otimes y.re| \oplus |x.im\otimes y.im|) 
\hfill\cr
\hfill{} \oplus
(|x.im \otimes y.re| \oplus |(-x.re)\otimes y.im|)) \oslash nrm)
\quad\cr}$$
$$=(EPS/2) (1+EPS/2)^3 (5 + 3EPS/2 + (EPS/2)^2)(A)$$
$$\le(1 - (EPS/2))^2(5EPS/2) (1 + 3EPS)(A)$$
$$\le(1 - (EPS/2))(5EPS/2) ((1 + 3EPS) \otimes A)$$
$$\le(5EPS/2) \otimes ((1 + 3EPS) \otimes A)$$
Here we used the fact that
$$(1+EPS/2)^3 (5 + 3EPS/2 + (EPS/2)^2) \le (1-EPS/2)^2(5)(1+3EPS)$$ \qed
\bigskip
\noindent {\bf Proposition 7.13 (A / A):}

$x  / y  = (re,im;e)$ where 

$re = (x.re \otimes y.re \oplus x.im \otimes y.im) \oslash nrm$ where $nrm = y.re \otimes y.re \oplus y.im \otimes y.im$

$im = (x.im \otimes y.re \ominus x.re \otimes y.im) \oslash nrm$

$e = (1 + 4EPS)\otimes 
(
((5 EPS/2) \otimes A \oplus (1+103EPS) \otimes B)
\oslash nrm)$ where
$$A = (|x.re \otimes y.re| \oplus |x.im\otimes y.im|) 
 \oplus
(|x.im \otimes y.re| \oplus | x.re \otimes y.im|)$$
$$B = x.e\otimes(|y.re|\oplus|y.im|)
\oplus
   (|x.re|\oplus|x.im|)\otimes y.e$$
 We are assuming (the computer will complain if this is not true)  $y.e < 100 EPS \otimes |y|$, or, more accurately, 
$$y.e^2< ((10000 EPS) \otimes EPS)\otimes nrm.$$
\bigskip
{\bf Proof:}

The error is bounded by $|x/y - x \oslash y| + ({x \over y}).e$ where $$({x \over y}).e = {|x| + x.e \over |y| - y.e} - {|x| \over |y|}$$
$$={|x||y| + |y|x.e - |x||y| + |x| y.e \over |y|(|y| - y.e)} $$
$$={|y|x.e + |x| y.e \over |y|(|y| - y.e)} $$
$$\le (1 + 101 EPS){|y|x.e + |x| y.e \over |y|^2} $$
We have used the fact that our assertion $y.e^2< (10000 EPS \otimes EPS)\otimes nrm$ implies that
$${1 \over |y|(|y| - y.e)} \le (1 + 101 EPS){1 \over |y|^2},$$  
which we derive in the next seven lines.

$$y.e^2< (10000 EPS \otimes EPS)\otimes nrm
 \le (1 + (EPS/2))(10000 EPS \otimes EPS) \times nrm$$
$$\le (1 + (EPS/2))^3(10000 EPS \otimes EPS) \times 
(y.re^2 + y.im^2)$$
$$\le (1 + (EPS/2))^4(10000 EPS \times EPS) \times 
(y.re^2 + y.im^2)$$
$$= (1 + (EPS/2))^4(100EPS)^2 \times |y|^2$$

This implies that $y.e < A \times |y|$, where
$A = (1 + (EPS/2))^2 (100 EPS)$. Now, noting that ${1 \over 1-A} \le 1 + 101 EPS$), we see that 
$${1 \over |y|(|y| - y.e)} \le {1 \over |y|(|y| - A|y|)}
={1 \over |y|^2(1 - A)}\le (1 + 101 EPS){1 \over |y|^2}$$

Resuming with our main proof, and setting
$$A= (|x.re \otimes y.re| \oplus |x.im\otimes y.im|) 
\oplus
(|x.im \otimes y.re| \oplus | x.re \otimes y.im|)$$
$$B= x.e\otimes(|y.re|\oplus|y.im|)
\oplus(|x.re|\oplus|x.im|)\otimes y.e$$
we have that
$$|x/y - x \oslash y| + ({x \over y}).e$$
$$\le |x/y - x \oslash y| + (1 + 101 EPS){|y|x.e + |x| y.e \over |y|^2}$$
$$\displaylines{\quad\le (EPS/2) (1+EPS/2)(5 + 3EPS/2 + (EPS/2)^2)
((|x.re \otimes y.re| \oplus |x.im\otimes y.im|) 
\hfill\cr
\hfill{} +
(|x.im \otimes y.re| \oplus | x.re \otimes y.im|)) / nrm
\hfill\cr
\hfill{} +
( 1 + 101 EPS)
(x.e\times (|y.re|+|y.im|)+(|x.re|+|x.im|)\times y.e)/(|y|^2)
\quad\cr}$$
$$\displaylines{\quad\le(5EPS/2) (1+EPS/2)^2(1 + (3EPS/10) + (EPS^2/20))
((|x.re \otimes y.re| \oplus |x.im\otimes y.im|) 
\hfill\cr
\hfill{} \oplus
(|x.im \otimes y.re| \oplus | x.re \otimes y.im|)) / nrm
\hfill\cr
\hfill{} +
( 1 + 101 EPS)(1+EPS/2)
(x.e\times(|y.re|\oplus|y.im|)
\hfill\cr
\hfill{}+
(|x.re|\oplus|x.im|)\times y.e)/(({1 \over (1 + EPS/2)})^2 nrm)
\quad\cr}$$
$$\displaylines{\quad\le(5EPS/2) (1+EPS/2)^3
((|x.re \otimes y.re| \oplus |x.im\otimes y.im|) 
\hfill\cr
\hfill{} \oplus
(|x.im \otimes y.re| \oplus | x.re \otimes y.im|)) / nrm
\hfill\cr
\hfill{} +
( 1 + 101 EPS)(1+EPS/2)^3 (1 + EPS/2)^2
(x.e\otimes(|y.re|\oplus|y.im|)
\hfill\cr
\hfill{}\oplus
(|x.re|\oplus|x.im|)\otimes y.e)/nrm
\quad\cr}$$
$$= (1+EPS/2)^3(
(5EPS/2) (A) +
( 1 + 101 EPS)(1 + EPS/2)^2 (B)
)/nrm $$
$$(1+EPS/2)^3(
(5EPS/2) A +
( 1 + 103EPS) B
)/nrm $$
$$(1+EPS/2)^6((
(5EPS/2) \otimes A \oplus
( 1 + 103EPS)\otimes B
) \oslash nrm) $$
$$(1+4EPS) \otimes ((
(5EPS/2) \otimes A \oplus
( 1 + 103EPS)\otimes B
) \oslash nrm) $$
Here we have used the fact that
$$(1+EPS/2)^2(1 + (3EPS/10) + (EPS^2/20)) \le (1+(EPS/2)^3)$$ \qed
\bigskip
\noindent {\bf Proposition 7.14 ($\sqrt{X}$):}

Let $s_o = \root o \of {(|x.re| \oplus h_o (x.re,x.im)) \otimes 0.5}$
 and $d_o = (x.im \oslash s) \otimes 0.5$, then
$\sqrt{x}  = (re,im;e)$ where 

$re = s_o$ if $x.re > 0.0$ and $re = d_o$ otherwise,

$im = d_o$ if $x.re > 0.0$ and $im = s_o$ otherwise,

$e = EPS \otimes ((1 + 4 EPS) \otimes (1.25 \otimes s_o \oplus 1.75 \otimes |d_o|))$
\bigskip
{\bf Proof:}

This will be a little nasty.  Let's begin by analyzing $e_s$, which is the difference between the true calculation of $s$ and the machine calculation of $s$, that is $e_s = |s - s_o|.$  
First, we bound $s.$
$$s = \sqrt{(|x.re| + h(x.re,x.im)) * 0.5}$$
$$\le (1+EPS)^{1/2} \sqrt{(|x.re| + h_o(x.re,x.im)) * 0.5}$$
$$\le (1+EPS)^{1/2} (1+EPS/2)^{1/2}
\sqrt{(|x.re| \oplus h_o(x.re,x.im)) * 0.5}$$
$$\le (1+EPS)^{1/2} (1+EPS/2)^{1/2}(1+EPS/2)
\root o \of {(|x.re| \oplus h_o(x.re,x.im)) * 0.5}$$
$$= (1+EPS)^{1/2} (1+EPS/2)^{3/2} s_o$$
By a power series expansion, we see that
$$(1+EPS)^{1/2} (1+EPS/2)^{3/2}$$
$$=(1 + {1 \over 2} EPS - {1 \over 8} EPS^2 + ...) +
     (1 + {3 \over 2} EPS/2 + {3 \over 8} (EPS/2)^2 + ...)$$
$$=(1 + {1 \over 2} EPS - {1 \over 8} EPS^2 + ...) +
     (1 + {3 \over 4} EPS + {3 \over 32} (EPS)^2 + ...)$$
$$=(1 + {5 \over 4} EPS + {11 \over 32} EPS^2 + ...)$$
So that,
$$s \le (1 + {5 \over 4} EPS + {11 \over 32} EPS^2 + ...) s_o$$
Similarly, 
$$s \ge (1 - {5 \over 4} EPS) s_o$$

Thus, we can bound the $s$ error, 
$$e_s = |s - s_o|$$
$$\le ((1 + {5 \over 4} EPS + {11 \over 32} EPS^2 + ...) - 1) s_o$$
$$= ({5 \over 4} EPS + {11 \over 32} EPS^2 + ...) s_o$$
\medskip
Next, we analyze $e_d$, which is the absolute value of the difference between the true calculation of $d$ and the machine calculation of $d$.  That is, $e_d = |d - d_o|.$
$$e_d = |x.im/(2s) - x.im \oslash (2 s_o)|$$
$$\le |x.im \oslash (2 s_o) - x.im / (2 s_o)| + 
|x.im/(2s_o) - x.im/(2s) |$$
$$\le (EPS/2) |x.im / (2 s_o)| + |{x.im \over 2} {s - s_o \over s s_o} |$$
$$\le (EPS/2) |x .im/ (2 s_o)| + |{x.im \over 2} {1 \over s s_o}
((5/4)EPS + (11/32)EPS^2+...) s_o|$$
$$\le (EPS/2) |x.im / (2 s_o)| + |{x.im \over 2}{1 \over s}
((5/4)EPS + (11/32)EPS^2+...)|$$
$$\le (EPS/2) |x.im / (2 s_o)| + |{x.im \over 2}
{1 \over s_o (1 - (5/4)EPS)}((5/4)EPS + (11/32)EPS^2+...)|$$
$$= (EPS/2) |x .im/ (2 s_o)|
(1 + {(5/2) + (11/16)EPS+...) \over (1 - (5/4)EPS)})$$
$$= (EPS/2) {(7/2) + (-9/16) EPS + ... \over (1 - (5/4)EPS)}
|x .im/ (2 s_o)|$$
$$\le (EPS/2) (1+EPS/2)
 {7/2 \over (1 - (5/4)EPS)}
|x .im \oslash (2 s_o)|$$
$$= (EPS/2) (1+EPS/2)
 {7/2 \over (1 - (5/4)EPS)}
|d_o|$$
\bigskip
Finally, we can bound the overall error $e = e_s + e_d$.
$$e_s + e_d$$
$$\le ({5 \over 4} EPS + {11 \over 32} EPS^2 + ...) s_o +
 (EPS/2) (1+EPS/2)
 {7/2 \over (1 - (5/4)EPS)}|d_o|$$
$$\le ( EPS + {11 \over 40} EPS^2 + ...) ({5 \over 4} s_o) +
 EPS (1+EPS/2)
 {1 \over (1 - (5/4)EPS)}|{7 \over 4} d_o|$$
$$\le EPS (1+EPS/2)
 {1 \over (1 - (5/4)EPS)} ({5 \over 4} s_o) +
 EPS (1+EPS/2)
 {1 \over (1 - (5/4)EPS)}|{7 \over 4} d_o|$$
$$\le EPS (1+EPS/2)
 {1 \over (1 - (5/4)EPS)}
({5 \over 4} s_o + |{7 \over 4} d_o|)$$
$$\le EPS (1+EPS/2)^3
 {1 \over (1 - (5/4)EPS)}
({5 \over 4} \otimes s_o \oplus |{7 \over 4} \otimes d_o|)$$
$$\le EPS (1 - (EPS/2))(1 + 4 EPS) 
({5 \over 4} \otimes s_o \oplus |{7 \over 4} \otimes d_o|)$$
$$\le EPS \otimes ((1 + 4 EPS) \otimes
({5 \over 4} \otimes s_o \oplus |{7 \over 4} \otimes d_o|))$$ \qed
\bigskip
Finally, we develop a couple of formulas for the absolute value of an XComplex.
\bigskip
\noindent {\bf Formula 7.0 (absUB(X)):}

If $x$ is an XComplex, then
we get an upper bound on the absolute value of $x$ as follows.

$$|x| = h(x.re,x.im) \le (1 + EPS) h_{\circ} (x.re,x.im)$$
$$\le (1 - EPS/2) (1 + 2EPS)h_{\circ} (x.re,x.im)$$
$$\le (1 + 2EPS) \otimes h_{\circ} (x.re,x.im)$$
Thus, we define 
$$absUB(x) = (1 + 2EPS) \otimes h_{\circ} (x.re,x.im)$$ \qed
\bigskip
\noindent {\bf Formula 7.1 (absLB(X)):}

If $x$ is an XComplex, then
we get a lower bound on the absolute value of $x$ as follows.

$$|x| = h(x.re,x.im) \ge (1 - EPS) h_{\circ} (x.re,x.im)$$
$$\ge (1 + EPS/2) (1 - 2EPS) h_{\circ} (x.re,x.im)$$
$$\ge (1 - 2EPS) \otimes h_{\circ} (x.re,x.im)$$
Thus, we define 
$$absLB(x) = (1 - 2EPS) \otimes h_{\circ} (x.re,x.im)$$ \qed
\bigskip
Finally, in several places in the programs {\it verify} and {\it fudging} we perform a standard operation on a pair of doubles and must take into account round-off error.  This is easy if we use Lemma 7.0.

For example, in {\it inequalityHolds} we want to show that $wh \times wh > absUB(along),$ where $wh = absLB(whirle).$   
By Lemma 7.0, we know that $(1 - EPS) \otimes (wh \otimes wh) \le wh \times wh$ and we simply test that 
$(1 - EPS) \otimes (wh \otimes wh) \ge absUB(along).$

Similar situations occur in the functions {\it horizon} and {\it larger-angle} in {\it fudging}.

A slightly more complicated version of this occurs in the computer calculation of $pos[i]$ and $size[i],$ that is, the center and size of a sub-box.  Prior to multiplication by $scale[i] = 2 ^{(5 - i)/6},$  the calculations of $pos$ and $size$ are exact.  However, multiplication by $scale$ introduces round-off error.  For the center of the box we will have the computer use $pos[i] \otimes scale[i]$ with the realization that this is not necessarily $pos[i] \times scale[i]$.   Thus, we have to choose appropriate sizes to ensure that the machine sub-box contains the true sub-box.  

Notationally, this is annoying, because we typically use a computer command like $pos[i] = pos[i] \otimes scale[i],$ while in an exposition, we need to avoid that.  We will denote the true center of the box by $p[i]$ and the machine center of the box by $p_0[i],$ and the true and machine sizes will be denoted $s[i]$ and $s_0[i].$  We will let $pos[i]$ and $size[i]$ be the position and size (true and machine are the same) before multiplication by $scale[i].$

Let $p[i] = pos[i] \times scale[i],\ p_0[i] = pos[i] \otimes scale[i],$ and 
$s[i]= size[i] \times scale[i].$  We must select $s_0[i]$ so that 
$p_0[i] + s_0[i] \ge p[i] + s[i].$  (Here, taking $+$ on the left-hand side is correct, because the need for machine calculation there is incorporated at other points in the programs.)  So, we must find $s_0[i]$ such that $s_0[i] \ge (p[i] - p_0[i]) + s[i].$
$$(p[i] - p_0[i]) + s[i]. \le (EPS/2) |p_0[i]| + size[i] \times scale[i]$$
$$\le (EPS/2) |p_0[i]| + (1 + EPS/2) (size[i] \otimes scale[i])$$
$$\le (1 + EPS/2) ((EPS/2) |p_0[i]| +  (size[i] \otimes scale[i]))$$
$$\le (1 + EPS/2)^2 ((EPS/2) |p_0[i]| \oplus  (size[i] \otimes scale[i]))$$
$$\le (1 + 2EPS) \otimes ((EPS/2) |p_0[i]| \oplus  (size[i] \otimes scale[i]))$$
Thus we take 
$$s_0 [i]= (1 + 2EPS) \otimes ((EPS/2) |p_0[i]| \oplus  (size[i] \otimes scale[i]))$$
This also works to give $p_0[i] - s_0[i] \le p[i] - s[i].$  \qed

\sectionskip\magnification=1200
\baselineskip = 0.2 true in
\def\qed{{\bf \ \ \ \ \vrule height6pt width5pt depth1pt}}

\noindent {\bf Chapter 8: AffApprox's with Round-Off Error}
\medskip
In Chapter 6, we saw how to do calculations with AffApprox's.  In this chapter, we incorporate round-off error into these calculations.  \medskip \noindent {\bf Convention 8.1:}
Recall that an AffApprox $x$ is a five-tuple
$(x.f; x.f_0,x.f_1,x.f_2; x.err)$
consisting of four complex numbers $(x.f, x.f_0,x.f_1, x.f_2)$ and one real number $x.err.$  In Chapter 6, the real number was denoted $x.e,$ but it seems preferable to use $x.err$ in this Chapter.  Also, we will suppress mention of the $r$ functions used in Chapter 6 (see Definition 6.4).  As such, statements in this section will have to be translated (implicitly) to account for this suppression.
\medskip \noindent {\bf Remark 8.2:}
One approach to round-off error for AffApprox's would be to replace the four complex numbers by four AComplex numbers complete with their round-off errors, and similarly for the one real number.  We will not do this because it would necessitate keeping track of five separate round-off-error terms when we do AffApprox calculations.

Instead, we will 
replace the four complex numbers by four XComplex numbers  (``exact" complex numbers) and push all the round-off error into the $.err$ term.  In particular, in doing an AffApprox calculation, our subsidiary calculations will generally be on XComplex numbers and produce an AComplex number whose $.e$ term will be plucked off and forced into the $.err$ term of the final AffApprox.
\medskip \noindent  {\bf Convention 8.3:}
In what follows, we will use Basic Properties 7.0 and Lemmas 7.0 and 7.1.  Also, the Propositions in Chapter 6 will be utilized;  as such, the numbering of the Propositions is the same in both Chapters (for example, Propositon 6.7 corresponds to Proposition 8.7).
\medskip
\noindent {\bf -X:}
\medskip
\noindent {\bf Proposition 8.1:}
If $x = (x.f; x.f_0,x.f_1,x.f_2; x.err)$, then 
$$-x = (-x.f; -x.f_0,-x.f_1,-x.f_2; x.err)$$ \qed
\medskip
\noindent {\bf X + Y:}
\medskip
We analyze the addition of the AffApproxs 
$x = (x.f; x.f_0,x.f_1,x.f_2; x.err)$ and 
$y = (y.f; y.f_0,y.f_1,y.f_2; y.err).$  To get the first term in $x+y$ we add the XComplex numbers $x.f$ and $y.f;$ which produces the AComplex number $r\_f = x.f + y.f,$ and then we pluck off the XComplex part, $r\_f.z.$  The round-off error part $r\_f.e$ will be foisted into the overall error term $r\_error$ for $x+y.$ Similarly for the next three terms in $x + y.$

Abstractly, the overall error term $r\_error$ comes from adding the round-off error contributions $r\_f.e,\  r\_f_0.e,\  r\_f_1.e, 
\ r\_f_2.e$ and the AffApprox error contributions $x.err,\  y.err$.  Of course, we have to produce a machine version.
\medskip
\noindent {\bf Proposition 8.2:} If 
$x = (x.f; x.f_0,x.f_1,x.f_2; x.err)$ and 
$y = (y.f; y.f_0,y.f_1,y.f_2; y.err),$ then
$$x + y = (r\_f.z; r\_f_0.z, r\_f_1.z, r\_f_2.z; r\_error)$$
where 
$$r\_error = (1 + 3EPS) \otimes ((x.e \oplus y.e) \oplus ((r\_f.e \oplus r\_f_0.e) \oplus (r\_f_1.e \oplus r\_f_2.e)))$$

\noindent {\bf Proof:}
The error is given by 
$$(x.e + y.e) + ((r\_f.e + r\_f_0.e) + (r\_f_1.e + r\_f_2.e))$$
$$\le (1+EPS/2) (x.e \oplus y.e) + (1+EPS/2)((r\_f.e \oplus r\_f_0.e) + (r\_f_1.e \oplus r\_f_2.e))$$
$$\le (1+EPS/2)^2 (x.e \oplus y.e) + (1+EPS/2)^2((r\_f.e \oplus r\_f_0.e) \oplus (r\_f_1.e \oplus r\_f_2.e))$$
$$\le (1+EPS/2)^3 ((x.e \oplus y.e) \oplus ((r\_f.e \oplus r\_f_0.e) \oplus (r\_f_1.e \oplus r\_f_2.e)))$$
$$\le (1 + 3EPS) \otimes ((x.e \oplus y.e) \oplus ((r\_f.e \oplus r\_f_0.e) \oplus (r\_f_1.e \oplus r\_f_2.e)))$$
To get the last line we used Lemma 1 in Section 5. \qed
\bigskip
\noindent {\bf X - Y:}
\medskip
\noindent {\bf Proposition 8.3:} If 
$x = (x.f; x.f_0,x.f_1,x.f_2; x.err)$ and 
$y = (y.f; y.f_0,y.f_1,y.f_2; y.err),$ then
$$x - y = (r\_f.z; r\_f_0.z, r\_f_1.z, r\_f_2.z; r\_error)$$
where
$$r\_f = x.f - y.f$$
$$r\_f_k = x.f_k - y.f_k$$
$$r\_error = (1 + 3EPS) \otimes ((x.e \oplus y.e) \oplus ((r\_f.e \oplus r\_f_0.e) \oplus (r\_f_1.e \oplus r\_f_2.e)))$$ \qed
\medskip
\noindent {\bf X + D}
\medskip
Here, we add the AffApprox $x = (x.f; x.f_0,x.f_1,x.f_2; x.err)$  to the double $y$.  The only terms that change are the first and the last.
\medskip
\noindent {\bf Proposition 8.4:} If 
$x = (x.f; x.f_0,x.f_1,x.f_2; x.err)$ and 
$y$ is a double, then
$$x + y = (r\_f.z; r\_f_0.z, r\_f_1.z, r\_f_2.z; r\_error)$$
where
$$r\_f = x.f + y$$
$$r\_f_k = x.f_k$$
$$r\_error = (1 + EPS) \otimes (x.err \oplus r\_f.e)$$
\medskip
\noindent {\bf Proof:}
The error is given by 
$$x.err + r\_f.e$$
$$\le (1 + EPS) \otimes (x.err \oplus r\_f.e)$$
by Lemma 0 in Section 5. \qed
\medskip
\noindent {\bf X - D}
\medskip
\noindent {\bf Proposition 8.5:} If 
$x = (x.f; x.f_0,x.f_1,x.f_2; x.err)$ and 
$y$ is a double, then
$$x - y = (r\_f.z; r\_f_0.z, r\_f_1.z, r\_f_2.z; r\_error)$$
where
$$r\_f = x.f - y$$
$$r\_f_k = x.f_k$$
$$r\_f.z = x.f - y$$
$$r\_error = (1 + EPS) \otimes (x.err \oplus r\_f.e)$$ \qed
\medskip
\noindent {\bf X * Y}
\medskip
We multiply the AffApproxs $x,y$ while pushing all error into the $.err$ term.

We will use the functions (see Formulas 7.0 and 7.1, at the end of Chapter 7) $absUB= (1 + 2EPS) \otimes hypot_o(x.re,x.im)$ and $absLB(x) = (1 - 2EPS) \otimes hypot_o(x.re,x.im)$.

When $x$ is an AffApprox, we define $dist(x)$ to be $$(1 + 2EPS) \otimes (absUB(x.f_0) \oplus (absUB(x.f_1) \oplus absUB(x.f_2))).$$  This is the machine representation of the sum of the absolute values of the linear terms in the AffApprox $x$ (the proof is straightforward).
\medskip
\noindent {\bf Proposition 8.6:} If 
$x = (x.f; x.f_0,x.f_1,x.f_2; x.err)$ and 
$y = (y.f; y.f_0,y.f_1,y.f_2; y.err),$ then
$$x * y = (r\_f.z; r\_f_0.z, r\_f_1.z, r\_f_2.z; r\_error)$$
where

$$r\_f = x.f \times y.f$$
$$r\_f_k = x.f \times y.f_k + x.f_k \times y.f$$
\centerline {Then, $x \times y = (r\_f.z; r\_f_0.z,r\_f_1.z,r\_f_2.z; r\_error)$, where}
$$r\_error = 
(1 + 3EPS) \otimes (A \oplus (B \oplus C))$$
with
$$A = (dist(x) \oplus x.e)\otimes (dist(y) \oplus y.e)$$
$$B = absUB(x.f) \otimes y.e \oplus absUB(y.f) \otimes x.e$$
$$C = (r\_f.e \oplus r\_f_0.e) \oplus (r\_f_1.e \oplus r\_f_2.e)$$
\noindent {\bf Proof:}
We add the non-round-off error term for $x \times y$ to the various round-off error terms that accumulated.
$$\displaylines{\quad 
((dist(x) + x.e)\times (dist(y) + y.e)) +
 ((absUB(x.f) \times y.e 
\hfill\cr
\hfill{} + 
absUB(y.f) \times x.e) +
 (r\_f.e + r\_f_0.e) + (r\_f_1.e + r\_f_2.e))
                                                \quad\cr}$$
$$\displaylines{\quad \le 
(1+EPS/2)^2[(dist(x) \oplus x.e)\times (dist(y) \oplus y.e)] +
 (1+EPS/2)(absUB(x.f) \otimes y.e 
\hfill\cr
\hfill{} + 
absUB(y.f) \otimes x.e) +
(1+EPS/2)((r\_f.e \oplus r\_f_0.e) + (r\_f_1.e \oplus r\_f_2.e))
                                                \quad\cr}$$
$$\displaylines{\quad \le 
(1+EPS/2)^3[(dist(x) \oplus x.e)\otimes (dist(y) \oplus y.e)] +
 (1+EPS/2)^2\{(absUB(x.f) \otimes y.e 
\hfill\cr
\hfill{} \oplus 
absUB(y.f) \otimes x.e) +
((r\_f.e \oplus r\_f_0.e) \oplus (r\_f_1.e \oplus r\_f_2.e))\}
                                                \quad\cr}$$
$$\le (1+EPS/2)^3 A + (1+EPS/2)^2 (B + C) $$
$$\le (1+EPS/2)^3 A + (1+EPS/2)^3 (B + C) $$
$$\le (1+EPS/2)^4 (A \oplus (B + C) )$$
$$\le (1+3EPS) \otimes (A \oplus (B + C) )$$ \qed
\bigskip
\noindent {\bf X * D:}
\medskip
\noindent {\bf Proposition 8.7:} If 
$x = (x.f; x.f_0,x.f_1,x.f_2; x.err)$ and 
$y$ is a double, then
$$x \times y = (r\_f.z; r\_f_0.z, r\_f_1.z, r\_f_2.z; r\_error)$$
where
$$r\_f = x.f \times y$$
$$r\_f_k = x.f_k \times y$$
$$r\_error = (1+3EPS) \otimes ( (x.e \otimes |y|) \oplus
( (r\_f.e \oplus r\_f_0.e) \oplus (r\_f_1.e \oplus r\_f_2.e)) )$$
\medskip
\noindent {\bf Proof:}

$$(x.e \times |y| ) + ( (r\_f.e + r\_f_0.e) + (r\_f_1.e + r\_f_2.e))$$
$$\le (1+EPS/2) (x.e \otimes |y| ) + (1+EPS/2)
( (r\_f.e \oplus r\_f_0.e) + (r\_f_1.e \oplus r\_f_2.e))$$
$$\le (1+EPS/2)^2 (x.e \otimes |y| ) + (1+EPS/2)^2
( (r\_f.e \oplus r\_f_0.e) \oplus (r\_f_1.e \oplus r\_f_2.e))$$
$$\le (1+EPS/2)^3 ((x.e \otimes |y| ) \oplus 
( (r\_f.e \oplus r\_f_0.e) \oplus (r\_f_1.e \oplus r\_f_2.e)))$$
$$\le (1 + 3EPS) \otimes ((x.e \otimes |y| ) \oplus 
( (r\_f.e \oplus r\_f_0.e) \oplus (r\_f_1.e \oplus r\_f_2.e)))$$ \qed
\bigskip
\noindent {\bf X/Y:}
\medskip
For convenience, let  $ax = absUB(x.f),\ ay = absLB(y.f).\ $ 
 \medskip
\noindent {\bf Proposition 8.8:} If 
$x = (x.f; x.f_0,x.f_1,x.f_2; x.err)$ and 
$y = (y.f; y.f_0,y.f_1,y.f_2; y.err),$ then
$$x / y = (r\_f.z; r\_f_0.z, r\_f_1.z, r\_f_2.z; r\_error)$$
where
$$r\_f = x.f / y.f$$
$$r\_f_k = (x.f_k \times y.f - x.f \times y.f_k)/ (y.f \times y.f)$$
$$r.error = (1 + 3EPS) \otimes (
((1 + 3EPS) \otimes A \ominus (1 - 3EPS) \otimes B )\oplus C ) $$
with
$$A = (ax \oplus (dist(x) \oplus x.e)) \oslash D$$
$$B = (ax \oslash ay \oplus dist(x) \oslash ay) \oplus  ((dist(y) \otimes ax) \oslash (ay \otimes ay))$$
$$C = (r\_f.e \oplus r\_f_0.e) \oplus (r\_f_1.e \oplus r\_f_2.e)$$
$$D = ay \ominus (1 + EPS) \otimes (dist(y) \oplus y.e)$$
Of course, we do have to be concerned about division by zero, so the program will complain if $ay$ is not greater than $(1 + EPS) \otimes (dist(y) \oplus y.e)$, that is, if $D$ is not greater than zero.
\medskip
\noindent {\bf Proof:}

As usual, we add the round-off errors to the old AffApprox error, 
taking into account round-off error.  Let's work on it bit by bit.

$$(ax + dist(x) + x.e)/(ay - (dist(y) + y.e))$$
$$\le (ax + (1 + EPS/2)(dist(x) \oplus x.e))/
         (ay - (1 + EPS) \otimes (dist(y) \oplus y.e))$$
$$\le (1 + EPS/2)^2 (ax \oplus (dist(x) \oplus x.e))/
         (ay - (1 + EPS) \otimes (dist(y) \oplus y.e))$$
$$\le (1 + EPS/2)^2 (ax \oplus (dist(x) \oplus x.e))/
         ({1 \over 1 + EPS/2}) (ay \ominus (1 + EPS) \otimes (dist(y) \oplus y.e))$$
$$\le (1 + EPS/2)4 (ax \oplus (dist(x) \oplus x.e)) \oslash
         (ay \ominus (1 + EPS) \otimes (dist(y) \oplus y.e))$$
$$= (1 + EPS/2)^4 A$$
$$\le   (1 + 3EPS) \otimes A$$

The next term, being subtracted, requires opposite inequalities.
$$(ax/ay + dist(x)/ay) + dist(y) \times ax / (ay \times ay)$$
$$\ge (1-EPS/2)(ax \oslash ay + dist(x) \oslash ay) + 
(1-EPS/2) (dist(y) \otimes ax )/ ({1 \over 1 - EPS/2})(ay \otimes ay)$$
$$\ge ((1-EPS/2)^2 (ax \oslash ay \oplus dist(x) \oslash ay) + 
(1 - EPS/2)^3 (dist(y) \otimes ax) \oslash (ay \otimes ay)$$
$$\ge (1 - EPS/2)^3 [(ax \oslash ay \oplus dist(x) \oslash ay) +  (dist(y) \otimes ax) \oslash (ay \otimes ay))$$
$$\ge ((1-EPS/2)^4 ((ax \oslash ay \oplus dist(x) \oslash ay) \oplus  ((dist(y) \otimes ax) \oslash (ay \otimes ay)))$$
$$\ge (1 + EPS/2)(1 + 3EPS) (B)$$
$$\ge (1 - 3EPS) \otimes B$$

Finally, we do the round-off terms.  
$$((r\_f.e + r\_f_0.e) + (r\_f_1.e + r\_f_2.e))$$
$$\le (1 + EPS/2)((r\_f.e \oplus r\_f_0.e) + (r\_f_1.e \oplus r\_f_2.e))$$
$$\le (1 + EPS/2)^2 ((r\_f.e \oplus r\_f_0.e) \oplus (r\_f_1.e \oplus r\_f_2.e))$$
$$= (1 + EPS/2)^2 C$$
 Now, we put these three pieces together.

$$\displaylines{\quad 
(ax + dist(x) + x.e)/(ay - (dist(y) + y.e))
\hfill\cr
\hfill{} - 
((ax/ay + dist(x)/ay) + dist(y) \times ax / (ay \times ay))
+ ((r\_f.e + r\_f_0.e) + (r\_f_1.e + r\_f_2.e))
\quad\cr}$$
$$\le 
(1 + 3EPS) \otimes A - (1 - 3EPS) \otimes B + (1 + EPS/2)^2 C $$
$$\le 
(1+EPS/2) ((1 + 3EPS) \otimes A \ominus (1 - 3EPS) \otimes B)
 + (1 + EPS/2)^2 C $$
$$\le (1+EPS/2)^2 (
((1 + 3EPS) \otimes A \ominus (1 - 3EPS) \otimes B)
 + C
) $$
$$\le (1+EPS/2)^3 (
((1 + 3EPS) \otimes A \ominus (1 - 3EPS) \otimes B)
 \oplus C
) $$
$$\le (1+3EPS) \otimes (
((1 + 3EPS) \otimes A \ominus (1 - 3EPS) \otimes B)
 \oplus C
) $$ \qed
\bigskip
\noindent {\bf D/X:}
\medskip
We are dividing a double $x$ by an AffApprox $y$.
For convenience, let $ax = |x|,\ ay = absLB(y.f).\ $ 
 \medskip
\noindent {\bf Proposition 8.9:} If 
$x$  is a double,  and 
$y = (y.f; y.f_0,y.f_1,y.f_2; y.err),$ then
$$x / y = (r\_f.z; r\_f_0.z, r\_f_1.z, r\_f_2.z; r\_error)$$
where
 $$r\_f = x / y.f$$
$$r\_f_k = -(x \times y.f_k)/ (y.f \times y.f)$$
$$r\_error = (1+3EPS) \otimes ( 
((1 + 2EPS) \otimes (ax  \oslash D) \ominus (1 - 3EPS) \otimes B)
\oplus C
)$$
$$B = ax \oslash ay  \oplus (dist(y) \otimes ax \oslash (ay \otimes ay))$$
$$C = (r\_f.e \oplus r\_f_0.e) \oplus (r\_f_1.e \oplus r\_f_2.e)$$
$$D = ay \ominus (1 + EPS) \otimes (dist(y) \oplus y.e)$$
Again, we have to be concerned about division by zero, so the program will complain if $ay$ is not greater than $(1 + EPS) \otimes (dist(y) \oplus y.e)$, that is if $D$ is not greater than zero.
\medskip
\noindent {\bf Proof:}

Let's start with the pieces.

$$ax / (ay - (dist(y) + y.e))$$
$$\le ax / (ay - (1+EPS) \otimes (dist(y) \oplus y.e))$$
$$\le ax / ({1 \over 1 + EPS/2}) (ay \ominus (1+EPS) \otimes (dist(y) \oplus y.e))$$
$$\le (1 + EPS/2)^2 ax \oslash D$$
$$\le (1 + 2EPS) \otimes (ax \oslash D)$$

The next term, being subtracted, requires opposite inequalities.
$$ax/ay + dist(y) \times ax / (ay \times ay)$$
$$\ge (1 - EPS/2) ax \oslash ay + 
(1 - EPS/2) dist(y) \otimes ax / ({1 \over 1 - EPS/2})(ay \otimes ay)$$
$$\ge (1 - EPS/2) ax \oslash ay + 
(1 - EPS/2)^3 dist(y) \otimes ax \oslash (ay \otimes ay)$$
$$\ge (1 - EPS/2)^3 (ax \oslash ay  +  dist(y) \otimes ax \oslash (ay \otimes ay))$$
$$\ge (1 - EPS/2)^4 [ax \oslash ay  \oplus  (dist(y) \otimes ax \oslash (ay \otimes ay))]$$
$$\ge (1 + EPS/2)(1 - 3EPS) B$$
$$\ge (1 - 3EPS) \otimes B$$

Finally, we do the round-off terms as before.  
$$((r\_f.e + r\_f_0.e) + (r\_f_1.e + r\_f_2.e))$$
$$\le (1 + EPS/2)^2 C$$

Putting the pieces together, we have:
$$\displaylines{\quad \le 
ax/(ay - (dist(y) + y.e))
\hfill\cr
\hfill{} - 
(ax/ay + dist(y) \times ax / (ay \times ay))
+ ((r\_f.e + r\_f_0.e) + (r\_f_1.e + r\_f_2.e))
\quad\cr}$$
$$\le (1 + 2EPS) \otimes (ax  \oslash D)
- (1 - 3EPS) \otimes B
+ (1 + EPS/2)^2 C$$
$$\le (1+EPS/2)((1 + 2EPS) \otimes (ax  \oslash D)
\ominus (1 - 3EPS) \otimes B)
+ (1 + EPS/2)^2 C$$
$$\le (1+EPS/2)^2(((1 + 2EPS) \otimes (ax  \oslash D)
\ominus (1 - 3EPS) \otimes B)
+ C)$$
$$\le (1+EPS/2)^3(((1 + 2EPS) \otimes (ax  \oslash D)
\ominus (1 - 3EPS) \otimes B)
\oplus C)$$
$$\le (1+3EPS) \otimes (((1 + 2EPS) \otimes (ax  \oslash D)
\ominus (1 - 3EPS) \otimes B)
\oplus C)$$ \qed
\bigskip
\noindent {\bf X/D:}
\medskip
We are dividing an AffApprox $x$ by a double $y$ and the computer will object if $y = 0$. 
\medskip
\noindent {\bf Proposition 8.10:} If 
$x = (x.f; x.f_0,x.f_1,x.f_2; x.err)$ and 
$y$ is a double, then
$$x / y = (r\_f.z; r\_f_0.z, r\_f_1.z, r\_f_2.z; r\_error)$$
where
  $$r\_f = x.f / y$$
$$r\_f_k = x.f_k / y$$
$$r\_error  = 
(1+3EPS) \otimes ( 
(x.e\oslash |y|)
\oplus 
[(r\_f.e \oplus r\_f_0.e) \oplus (r\_f_1.e \oplus r\_f_2.e)]
                                               )$$
\medskip
\noindent {\bf Proof:}

This is easy.
$$x.e/|y| + ((r\_f.e + r\_f_0.e) + (r\_f_1.e + r\_f_2.e))$$
$$\le (1 + EPS/2) x.e \oslash |y| + 
(1 + EPS/2)^2
[(r\_f.e \oplus r\_f_0.e) \oplus (r\_f_1.e \oplus r\_f_2.e)]$$
$$\le (1 + EPS/2)^2 (x.e \oslash |y| + 
[(r\_f.e \oplus r\_f_0.e) \oplus (r\_f_1.e \oplus r\_f_2.e)])$$
$$\le (1 + EPS/2)^3 (x.e \oslash |y| \oplus 
[(r\_f.e \oplus r\_f_0.e) \oplus (r\_f_1.e \oplus r\_f_2.e)])$$
$$\le (1 + 3EPS) \otimes (x.e \oslash |y| \oplus 
[(r\_f.e \oplus r\_f_0.e) \oplus (r\_f_1.e \oplus r\_f_2.e)])$$ \qed
\bigskip
\noindent  ${\bf \sqrt X:}$
\medskip
Here, $x$ is an AffApprox and we let $ax = absUB(x.f)$.  There are two cases to consider depending on whether or not
$D = ax \ominus (1 + EPS) \otimes (dist(x) \oplus x.e)$ is or is not greater than zero.
\medskip
\noindent {\bf Proposition 8.11a:} If 
$x = (x.f; x.f_0,x.f_1,x.f_2; x.err)$  and
$D$ is not greater than zero, then we use the crude over-estimate $\sqrt x = (0;0,0,0;(1 + 2EPS) \otimes
\root o \of {(ax \oplus (xdist \oplus x.e))}\ )$
\medskip
\noindent {\bf Proof:}

$$\sqrt {ax + xdist + x.e}$$
$$\le \sqrt {ax + (1 + EPS/2) (xdist \oplus x.e)}$$
$$\le \sqrt {(1 + EPS/2)^2 (ax \oplus (xdist \oplus x.e))}$$
$$= (1 + EPS/2) \sqrt { (ax \oplus (xdist \oplus x.e))}$$
$$\le  (1 + EPS/2)^2 \root o \of { (ax \oplus (xdist \oplus x.e))}$$
$$\le  (1 + 2EPS) \otimes \root o \of { (ax \oplus (xdist \oplus x.e))}$$ \qed
\medskip
\noindent {\bf Proposition 8.11b:} If 
$x = (x.f; x.f_0,x.f_1,x.f_2; x.err)$  and
$D$ is greater than zero,
then, $$\sqrt{x} = (r\_f.z; r\_f_0.z,r\_f_1.z,r\_f_2.z; r\_error)$$ where
  $$r\_f = \sqrt {x.f}$$
 $$t = r\_f + r\_f$$
$$r\_f_k = {\rm AComplex} (x.f_k.re,x.f_k.im,0) / t$$
(Simply put, $r\_f_k = x.f_k / (2\sqrt{x.f})$. The reason we have to fuss to define $r\_f_k$ is because $\sqrt{x.f}$ is an AComplex.)
$$\displaylines{\quad
r\_error = 
(1+3EPS) \otimes ( 
\hfill\cr
\hfill{} \{
(1+EPS) \otimes \root o \of {ax}
\hfill\cr
\hfill{}\ominus
(1 - 3EPS) \otimes [dist(x)\oslash(2 \times \root o \of {ax})
\oplus \root o \of {D}\ ]\}
\hfill\cr
\hfill{}\oplus 
((r\_f.e \oplus r\_f_0.e) \oplus (r\_f_1.e \oplus r\_f_2.e))
 \hfill\cr
\hfill{}                                              )\quad\cr}$$
\medskip
\noindent {\bf Proof:}

Let's work on the pieces.

$$\sqrt {ax} \le (1 + EPS) \otimes \root o \of {ax}$$

Next,
$$dist(x)/(2 \sqrt {ax}) + \sqrt {ax - (dist(x) + x.e)}$$
$$\ge  (1 - EPS/2)^2 dist(x) \oslash (2 \root o \of {ax}) + \sqrt {ax - (1+EPS) \otimes (dist(x) \oplus x.e)}$$
$$\displaylines{\quad
\ge  (1 - EPS/2)^2 dist(x) \oslash (2 \root o \of {ax}) 
\hfill\cr
\hfill{} 
+ (1 - EPS/2)^{1/2}\sqrt {ax \ominus (1+EPS) \otimes (dist(x) \oplus x.e)}
                                                \quad\cr}$$
$$\ge (1 - EPS/2)^2 dist(x) \oslash (2 \root o \of {ax}) 
+ (1 - EPS/2)^{3/2}\root o \of {D}$$
$$\ge (1 - EPS/2)^2 [dist(x) \oslash (2 \root o \of {ax}) + \root o \of {D}\ ]$$
$$\ge(1 - EPS/2)^3 [dist(x) \oslash (2 \root o \of {ax}) \oplus \root o \of {D}\ ]$$
$$\ge (1 + EPS/2) (1 - 3EPS) [dist(x) \oslash (2 \root o \of {ax}) \oplus \root o \of {D}\ ]$$
$$\ge (1 - 3EPS) \otimes [dist(x) \oslash (2 \root o \of {ax}) \oplus \root o \of {D}\ ]$$

Adding in the usual term, we get as our error bound
$$\sqrt {ax} - (dist(x)/(2 \sqrt {ax}) + \sqrt {ax - (dist(x) + x.e)}\ )
+ ((r\_f.e + r\_f_0.e) + (r\_f_1.e + r\_f_2.e))$$
$$\displaylines{\quad
\le (1 + EPS) \otimes \root o \of {ax} 
\hfill\cr
\hfill{} 
- (1 - 3EPS) \otimes [dist(x) \oslash (2 \root o \of {ax}) \oplus 
\root o \of {D}\ ]
\hfill\cr
\hfill{} 
+ (1 + EPS/2)^2 ((r\_f.e \oplus r\_f_0.e) \oplus (r\_f_1.e \oplus r\_f_2.e))
                                                \quad\cr}$$
$$\displaylines{\quad
\le (1 + EPS/2) \{(1 + EPS) \otimes \root o \of {ax} 
\hfill\cr
\hfill{} 
\ominus (1 - 3EPS) \otimes [dist(x) \oslash (2 \root o \of {ax}) \oplus \root o \of {D}\ ]\}
\hfill\cr
\hfill{} 
+ (1 + EPS/2)^2 ((r\_f.e \oplus r\_f_0.e) \oplus (r\_f_1.e \oplus r\_f_2.e))
                                                \quad\cr}$$
$$\displaylines{\quad
\le (1 + EPS/2)^3 (\{(1 + EPS) \otimes \root o \of {ax} \hfill\cr
\hfill{} 
\ominus (1 - 3EPS) \otimes [dist(x) \oslash (2 \root o \of {ax}) \oplus \root o \of {D}\ ]\}
\hfill\cr
\hfill{} 
\oplus ((r\_f.e \oplus r\_f_0.e) \oplus (r\_f_1.e \oplus r\_f_2.e))
                                                )\quad\cr}$$
$$\displaylines{\quad
\le (1 + 3EPS) \otimes (
\{(1 + EPS) \otimes \root o \of {ax} 
\hfill\cr
\hfill{} 
\ominus (1 - 3EPS) \otimes [dist(x) \oslash (2 \root o \of {ax}) \hfill\cr
\hfill{} 
\oplus \root o \of {D}\ ]\}
\oplus ((r\_f.e \oplus r\_f_0.e) \oplus (r\_f_1.e \oplus r\_f_2.e))
                                                )\quad\cr}$$ \qed

\sectionskip\magnification=1200
\baselineskip = 0.2 true in
\centerline {\bf References}
\bigskip
\noindent\hang [B]  A. Beardon, {\it The Geometry of Discrete Groups}, Graduate Texts in mathematics 91, Springer-Verlag, 
New York, 1983.
\medskip
\noindent\hang [Be]  J. Berge, Heegaard, a program to understand Heegaard splittings of
3-manifolds, available free from the author.
\medskip
\noindent\hang [F]  W. Fenchel,  Elementary Geometry in Hyperbolic Space, de Gruyter Studies Math., 11, de Gruyter, Berlin (1989).
\medskip
\noindent\hang [G] D. Gabai, On the Geometric and Topological Rigidity of Hyperbolic 3-Manifolds,  to appear, {\it J. Amer. Math. Soc.\/}
\medskip
\noindent \hang [GM1] F. Gehring and G. Martin, Inequalities of M\"obius Transformations and Discrete Groups, {\it j. Reine Agnew. math.\/} {\bf 418} (1991) 31 - 76.
\medskip
\noindent \hang [GM2] F. Gehring and G. Martin, Precisely invariant collars and the volume of hyperbolic 3-folds, preprint.
\medskip
\noindent\hang [HW1] C. Hodgson and J. Weeks,  personal communication.
\medskip
\noindent\hang [HW2] C. Hodgson and J. Weeks, Symmetries, Isometries and Length Spectra of Closed Hyperbolic 3-Manifolds, 
{\it Exp. Math.\/} {\bf 3} (1994), 261-274.
\medskip
\noindent\hang [JR] K. Jones and A. Reid, 
\medskip
\noindent\hang [K]  W. Kahane, Interval Arithmetic Options in the Proposed IEEE Floating Point Arithmetic Standard; {\it Interval Mathematics 1980},  Ed. Karl L. E. Nickel, 99-128.
\medskip
\noindent\hang [Math] Mathematica 2.0, A system for doing mathematics by computer, 
Wolfram Research Inc., Champaign Il.
\medskip
\noindent\hang [M1] G. R. Meyerhoff, A Lower Bound for the Volume of Hyperbolic 
3-Manifolds, {\it Canadian J. Math. \/} {\bf 39} (1987) 1038-1056.
\medskip
\noindent\hang [M2] G. R. Meyerhoff,  Sphere-Packing and Volume in Hyperbolic 3-Space, {\it Comment. Math. Helvetici\/} {\bf 61} (1986) 271-278.
\medskip
\noindent\hang [Mo] G. D. Mostow, Quasi-Conformal Mappings in n-Space and the Rigidity of Hyperbolic Space Forms, {\it IHES Publ. Math. \/} {\bf 34} (1968) 53-104.
\medskip
\noindent\hang [T] N. Thurston,  Finding Killer Words, in preparation. 
\medskip
\noindent \hang [Wa] F. Waldhausen, On Irreducible 3-Manifolds which Are Sufficiently Large, {\it Annals of Math. \/} {\bf 87} (1968) 56-88.
\medskip
\noindent\hang [W]  J. Weeks, SNAPPEA, available from the Geometry Center, geom.umn.edu.

\bye